\documentclass[11pt]{amsart}
\usepackage{amssymb,amsbsy,latexsym,amsmath,bbm,epsfig,psfrag,amsthm,mathrsfs,stmaryrd}
\usepackage[bookmarks=true]{hyperref}   % Hyperref in DVI and PDF
\usepackage[swedish,english]{babel}
\usepackage{graphicx,tikz,esint}
\usepackage[T1]{fontenc}
\usepackage[latin1]{inputenc}
\usepackage{amsfonts}
\usepackage{amsmath}
\usepackage{amssymb}
\usepackage{epsf}
\usepackage{graphics}
\usepackage{graphicx}
\usepackage{latexsym}
\usepackage{psfrag}
\usepackage{relsize}
\usepackage{verbatim}
\usepackage{hyperref}
\usepackage{pgfplots}
\usepackage{float}
\usetikzlibrary{positioning}
\usetikzlibrary{arrows}
\usetikzlibrary{decorations.pathreplacing}

\theoremstyle{plain}

\theoremstyle{definition}
\newtheorem{Lem}{Lemma}
\numberwithin{Lem}{section}
\newtheorem{Prop}{Proposition}
\numberwithin{Prop}{section}
\newtheorem{Thm}{Theorem}
\numberwithin{Thm}{section}

\numberwithin{Cor}{section}

\numberwithin{Con}{section}

\newtheorem{Def}{Definition}
\numberwithin{Def}{section}

\numberwithin{hyp}{section}

\numberwithin{conj}{section}
\newtheorem{ex}{Example}
\numberwithin{ex}{section}

\theoremstyle{remark}
\newtheorem{rem}{\bf{Remark}}
\numberwithin{rem}{section}

\numberwithin{equation}{section}

\DeclareMathOperator*{\wslim}{w^\ast-lim}
\DeclareMathOperator*{\wlim}{w-lim}
\DeclareMathOperator*{\esslim}{\text{ess lim}}
\DeclareMathOperator*{\esup}{\text{ess sup}}
\DeclareMathOperator*{\esslimsup}{\text{ess limsup}}
\DeclareMathOperator*{\essliminf}{\text{ess liminf}}
\DeclareMathOperator{\pv}{p.v.}

\newcommand{\dv}{\partial}
\newcommand{\Om}{\Omega}

\newcommand{\eps}{\varepsilon}
\newcommand{\R}{{\mathbb R}}
\newcommand{\C}{{\mathbb C}}

\newcommand{\N}{{\mathbb N}}

\newcommand{\Di}{\mathbb{D}}

\newcommand{\LL}{\mathcal{L}}

\newcommand{\A}{\mathscr{A}}
\newcommand{\Aa}{\mathbb{A}}

\newcommand{\gp}{\text{\tiny{$\triangle$}}}

\newcommand{\di}{\text{div}\,}

\newcommand{\ri}{\,\lrcorner\,}

\begin{document}
%\begin{flushright}

%December 8, 2008 (first part) 
%January 9, 2009 (second part) 
%\end{flushright}
\vspace{.4cm}

%\title{\bf \sffamily On Fuglede's flux extensions and the pointwise definition of linear partial differential operators}
\title[]
{\bf \sffamily On Fuglede's flux extensions and the pointwise definition of linear partial differential operators}
\author[Erik Duse]{Erik Duse}

\address{Erik Duse \\ Department of mathematics, KTH \\SE-100 44 Stockholm, Sweden} \email{duse@kth.se}

%\keywords{}
%\subjclass[2010]{42B35, 42B20, 42B37}

%\date{}

\dedicatory{Dedicated to  Bent Fuglede in admiration}

\thanks{ED was supported by the Swedish Research Council (VR), grant no. 2019-04152 and the Knut and Alice Wallenberg Foundation grant KAW 2015.0270.}

\maketitle

\begin{abstract}
In this work we provide a survey of Fuglede's flux extensions of first order partial differential operators, a concept largely forgotten today. A long the way we also survey the classical weak and strong extensions of PDE operators and the works of Friedrichs and Hörmander. We give several applications of this theory showing its usefulness, as well as connecting it to more recent developments in connection to various sharp versions of the divergence theorem. In particular, we use it to prove a generalization of Morera's theorem valid for general first order operators. Using this theory we also prove a new local limit formula for the maximal extension of a first order operator. We initiate a study of this limit and connect it to the wave cone of the operator, a concept that first arose in the theory of compensated compactness. Hopefully, this will contribute to a rivival of Fuglede's beautiful ideas. 
\end{abstract}

\tableofcontents

\newpage

\section{\sffamily Introduction}

%============NEW SUBSECTION=================================================================

\subsection{\sffamily Motivation and background}

Let $\Om \subset \R^n$ be a domain and let $E$ and  $F$ be a finite dimensional euclidean vector spaces. Let $\mathscr{A}: C^\infty_0(\Om,E)\to C^\infty_0(\Om,F)$
be a first order partial differential system with smooth coefficients given in a coordinate system by 
\begin{align}\label{eq:OpA}
\mathscr{A}u(x)=\sum_{j=1}^nA_j(x)\dv_ju(x)+B(x)u(x)
\end{align}
where $B,A_j\in C^\infty(\overline{\Om},\LL(E,F))$ for $j=1,2,...,n.$ 

A very important step in the history of the analysis of PDE:s was the realisation that in order prove solvability of partial differential equations one should first the extend the domain of definition of the operator $\A$ and allow it to act on functions that are not necessarily smooth. The question of whether or not these generalized solutions are then in fact smooth becomes a second question, that of the regularity theory of the PDE at hand. It turns out that in many cases solutions are in fact not smooth, and that by staying in very narrow confines of smooth functions, one would in retrospect never have been able to prove any general solvability result. 

How then should one go about extending an operator $\A$ to act on not necessarily continuously differentiable functions? The arguably most important and general way has been the use of integration by parts and duality, underlying the theory of weak derivatives and distribution theory. In the case of the linear PDE operator defined by \eqref{eq:OpA}, integration by parts gives for any $u\in C^\infty(\overline{\Om},E)$ and any $\phi\in C^\infty_0(\Om,F)$ 
\begin{align}\label{eq:OpAweak}
\int_{\Om}\langle \A u(x),\phi(x)\rangle_F dx +\int_{\Om}\langle u(x),\A^\ast \phi(x)\rangle_E dx=0,
\end{align}
where 
\begin{align}\label{eq:OpAdj}
\mathscr{A}^\ast v(x)=-\sum_{j=1}^nA_j(x)^\ast \dv_jv(x)+\bigg(B(x)^\ast -\sum_{j=1}^n \dv_j A_j(x)^\ast\bigg)v(x)
\end{align}
is the formal adjoint of $\A$ and $\langle \cdot,\cdot\rangle_E$ and $\langle \cdot,\cdot\rangle_F$ are the euclidean inner products on $E$ and $F$ respectively. In the \emph{weak extension} of $\A$ in various function spaces we insist that \eqref{eq:OpAweak} holds also for functions $u$ which are not necessarily differentiable. 
 
 Another way of creating an extension $\widetilde{\A}: X\to Y$ of $\A$ to functions $u$ in some function space $X$ such that the space of smooth functions $C^\infty(\Om,E)\cap X$ is dense in $X$, is to consider sequences $\{u_j\}_j\subset C^\infty(\Om,E)\cap X$ and $\{\mathscr{A}u_j\}_j\subset C^\infty(\Om,F)\cap Y$ such that both $\lim_{j\to \infty}u_j=u\in X$ and $\lim_{j\to \infty}\A u_j=v\in Y$ in the topologies of $X$ and $Y$ respectively, and then to \emph{define} $\A u:=v$. This leads to the notion of \emph{strong extensions} (the precise definition is deferred to Section \ref{sec:WeakStrong}).

 These are all very classical topics in the analysis of partial differential operators, going back to the work of K. O. Friedrichs, L. Schwartz, L. Hörmander, any many others and their importance, especially the notion of weak extension, can hardly be overestimated. However, something in common with all these extensions is that they do not provide an explicit way to construct the extended operator via a limit as in the smooth case. It is therefore not a priori obvious which information of the original function $u$ that is actually encoded in $\A u$. In addition, it is also hard to understand how the operator $\mathscr{A}$ acts on these weakly differentiable functions in a pointwise way directly from the definitions and why one actually gets a larger domain of definition. There is however another way to extend a first order partial differential operator system, due to B. Fuglede called \emph{flux extension}, which has received much less attention than it deserves. It is the purpose of the this paper to rectify the situation and hopefully revive the interest in this notion. Among other things, we will use flux extensions to prove a point-wise formula for closed extensions of $\A$ on $L^p(\Om,E)$. 

 We will now briefly outline the idea behind flux extensions and how they can be applied. They can be viewed as a vast extension of Cauchy's integral theorem in complex analysis. The precise definitions are given in Section \ref{sec:Flux}. Firstly, we recall that the \emph{principal symbol} of $\A\in C^1(\Om, \LL(\R^n, \LL(E,F)))$ is given by
\begin{align}\label{eq:PSymbol}
\Aa(x)=\sum_{j=1}^n A_j(x)\otimes e_j
\end{align} 
 and we write $\Aa(x,\xi):=\Aa(x)(\xi)\in C^\infty(\Om,\LL(E,F))$. Write $\mathscr{A}$ on the divergence form 
\begin{align}\label{eq:DivA}
\mathscr{A}u(x)=\sum_{j=1}^n\dv_j(A_j(x)u(x))+(B(x)-\text{div}\, \mathbb{A}(x))u(x), \quad \text{div}\, \mathbb{A}(x)=\sum_{j=1}^n\dv_jA_j(x)
\end{align}
and use Stokes theorem to get for every smooth domain $U\Subset \Om$ and every $u\in C^\infty(\overline{\Om},E)$

\begin{align}\label{eq:FluxInt}
\int_{U}\mathscr{A}u(x)dx=\int_U(B(x)-\text{div}\, \mathbb{A}(x))u(x)dx+\int_{\dv U}\mathbb{A}(y,\nu(y))u(y)d\sigma(y),
\end{align}
where 
\begin{align*}
\mathbb{A}(x,\xi)u(x)=\sum_{j=1}^nA_j(x)\xi_ju(x)
\end{align*}
and $\nu$ is the outward pointing unit normal vector field of $\dv U$. The observation is now that the right hand side of \eqref{eq:FluxInt} is well-defined in many cases even though $u$ is not necessarily smooth. Thus if there exists a $v\in L^{p}(\Om, F)$ such that the relation 
\begin{align*}
\int_{U}v(x)dx=\int_U(B(x)-\text{div}\, \mathbb{A}(x))u(x)dx+\int_{\dv U}\mathbb{A}(y,\nu(y))u(y)d\sigma(y),
\end{align*}
holds for at least almost every smooth domain $U$, in a sense that was made rigorous by B. Fuglede in \cite{F1} using the notion of \emph{moduli of a family of Green domains}, we can take this $v$ to be the extension of $\mathscr{A}$ acting on $u$. By abuse of notation we still call the extended operator $\mathscr{A}$. Using the Lebesgue differentiation theorem we see that 
\begin{align}\label{eq:LimitA}
\mathscr{A}u(x)=(B(x)-\text{div}\,\Aa(x))u(x)+\lim_{r\to 0^+}\frac{1}{\vert B_r(x)\vert}\int_{\dv B_r(x)}\mathbb{A}(y,\nu(y))u(y)d\sigma(y), 
\end{align}
defining the operator $\mathscr{A}$ pointwise, through a limit \footnote{Here the limit is the \emph{essential limit}, i.e., it holds for all $r$ outside a null set.}. This idea is not new, and was used for $\overline{\dv}$. by D. Pompieu who called it \emph{areolar derivative}, see \cite{Pomp}.

However, it does not seem to have been used in a systematic way for general first order operators. Moreover it is a very natural extension of 1-dimensional definition of derivative. Indeed the fundamental theorem of calculus (of which Stokes theorem can been regarded as a higher dimensional generalization) applied to an absolutely continuous function implies
\begin{align*}
f'(x)=\lim_{h\to 0}\frac{1}{h}\int_{x}^{x+h}f'(x)dx=\lim_{h\to 0}\frac{f(x+h)-f(x)}{h}. 
\end{align*}

Also from a physics point of view, \eqref{eq:LimitA} is very natural as the local form of a conservation or balance law, see \cite[ch. 1]{Da}.

Using a new notion of differentiability called harmonic differential, see Definition \ref{def:HarmonicDiff} , we will see how this limit can be computed and when it exists in a point-wise sense. Indeed, assume that $u$ has one-sided directional derivatives in all directions, i.e., that 
\begin{align*}
u'(x;\xi)=\lim_{\eps\to 0^+}\frac{u(x+\eps \xi)-u(x)}{\eps}
\end{align*}
exists for all $\xi \in S^{n-1}$. Note that we {\bf do not} require that $u'(x;\xi)$ is linear in $\xi$. Assume for simplicity that $\mathscr{A}$ is a homogenous constant coefficient operator and that $u$ is continuous. Then after an affine change of coordinates 
\begin{align*}
\mathscr{A}u(x)&=\lim_{r\to 0^+}\frac{1}{\vert B_r(x)\vert}\int_{\dv B_r(x)}\mathbb{A}(\nu(y))u(y)d\sigma(y)\\&=\lim_{r\to 0^+}\frac{1}{\vert B_r(x)\vert}\int_{\dv B_r(x)}\mathbb{A}(\nu(y))(u(y)-u(x))d\sigma(y)\\
&=\lim_{r\to 0^+}\frac{1}{\omega_n}\int_{S^{n-1}}\mathbb{A}(\xi)\frac{u(x+r\xi)-u(x)}{r}d\sigma(\xi)\\
&=\frac{1}{\omega_n}\int_{S^{n-1}}\mathbb{A}(\xi)u'(x;\xi)d\sigma(\xi),
\end{align*}
where $\omega_n$ is the Lebesgue measure of $B_1(0)$, and where we used the very important \emph{cancellation property} that 
\begin{align*}
\int_{S^{n-1}}\mathbb{A}(\xi)d\sigma(\xi)=0.
\end{align*}
Thus, we see that $\mathscr{A}u(x)$ can be given a meaning even at points where $u$ is not-differentiable, since a function $u$ may possess one-sided directional derivatives without being differentiable. It is the purpose of the present work to explore this in greater detail as well as providing a survey of Fuglede's flux extensions.

%============NEW SUBSECTION=================================================================

\subsection{\sffamily Overview of the paper}

In Section \ref{sec:Stokes} we present different versions of Stokes theorem, relying on very recent results, where we base our presentation on the survey papers \cite{CT,MM}. Section \ref{sec:WeakStrong} contains an overview of different closed extensions of differential operators on $L^p$ spaces, building primarily on the work of Kurt Friedrich and Lars Hörmander. Section \ref{sec:Flux} contains a survey of Fugelede's theory of flux extensions and constitute the heart of the matter of this paper. The section ends with applications of the theory, in particular a striking generalization of Cauchy's integral theorem and Morera's theorem to general first order operators, showing perhaps surprisingly that these theorems have little to do with complex analysis. In Section \ref{sec:Point} we show how flux extensions provide a way to define the maximal extension of a differential operator on $L^p$ pointwise through a limit integral formula over the unit sphere. Related to the point-wise limit we introduce a new notion of maximal operator associated to a differential operator $\A$ which we believe can useful also for other purposes. Furthermore, we show how the formulas we derive are completely analogous to those for singular integrals. In Section \ref{sec:Limit} we study the limit formula from Section \ref{sec:Point} in greater detail and show how to express the limit for elliptic operators using the weak total derivative. In particular we give a limit formula for the weak total derivative and relate it to harmonic extensions. For non-elliptic homogeneous constant coefficient operators we show using cancellation properties, how the limit formula can be well-defined even outside the Lebesgue set of the function. The section finishes by relating the limit formula to the concept of the wave cone of the operator, and how properties of the wave cones affect the limit. The last section is also the one which contains new results.

%============NEW SUBSECTION=================================================================
 
\subsection{\sffamily Notation and conventions}

In this paper $E$ and $F$ will always denote real finite dimensional euclidean vectors space. The inner products on $E$ and $F$ will be denoted by $\langle \cdot,\cdot\rangle_E$ and $\langle \cdot,\cdot\rangle_F$ respectively, or just $\langle \cdot, \cdot\rangle$ when the vector space is clear from the context. In addition $\LL(E,F)$ denote the space of linear transformations from $E$ to $F$, and if $E=F$ we only write $\LL(E)$. If $X\in \LL(E,F)$ then 
\begin{align*}
\Vert X\Vert=\sup_{\substack {v\in E\\ \neq 0}}\frac{\vert Xv\vert}{\vert v\vert}
\end{align*}
denotes the operator norm and 
\begin{align*}
\vert X\vert=\sqrt{\text{tr}(X^\ast X)}
\end{align*}
the Hilbert-Schmidt norm. We let $\mathcal{L}^n$ denote the Lebsegue measure on $\R^n$, and for a Lebesgue measurable set $A$ we write $\vert A\vert:=\mathcal{L}^n(A)$. In addition, $dx:=d \mathcal{L}^n(x)$. 
$\mathcal{H}^k$ will always denote the $k$-dimensional Hausdorff measure and when we write integrals on the form 
\begin{align*}
\int_\Gamma f(x)d\sigma(x)
\end{align*}
over some $n-1$-Hausdorff dimensional set $\Gamma$, we let $\sigma:=\mathcal{H}^{n-1}\lfloor \Gamma$. $B_r(x)$ will denote the open ball $\{y\in \R^n: \vert y-x\vert<r\}$ and $S^{n-1}=\{x\in \R^n: \vert x\vert=1\}$ the unit sphere in $\R^n$. Furthermore we set $\omega_n=\vert B_1(0)\vert$ and $\sigma_{n-1}=\mathcal{H}^{n-1}(S^{n-1})$. By essential limit $\esslim_{\eps\to 0^+}f_\eps(x)=f(x)$ we mean that for some $\delta>0$ there exists a measurable set $N\subset (0,\delta)$ with $\mathcal{L}^1(N)=0$ so that 
\begin{align*}
\lim_{\substack{\eps \to 0^+\\ \eps \in (0,\delta) \setminus N}}f_\eps(x)=f(x). 
\end{align*}
A first order operator $\A$ given by
\begin{align*}
\mathscr{A}u(x)=\sum_{j=1}^nA_j(x)\dv_ju(x)+B(x)u(x)
\end{align*}
where $B,A_j\in C^\infty(\Om,\LL(E,F))$ for $j=1,2,...,n$ is said to satisfy the \emph{standard assumptions}. Whenever, the standard assumptions are mentioned in the rest of the paper they refer to these assumptions on the coefficients. Finally, if we write $\A$ without any subscript like $\A_s$ of $\A_f$, $\A$ will always mean the weak extension $\A_w$.

%============NEW SECTION=================================================================

\section{\sffamily Stokes Theorem and first oder partial differential operators}\label{sec:Stokes}

%============NEW SUBSECTION=================================================================

\subsection{\sffamily Stokes theorem}

Define the matrix field 
\begin{align}\label{eq:MatrixField}
\mathcal{M}_\A(x,u(x))=\sum_{j=1}^n\Aa(x,e_j)u(x)\otimes e_j,
\end{align}
so that $\mathcal{M}_\A(x,u(x))\in \LL(\R^n,F)$. Then, for $u\in C^1(\Om,E)$

\begin{align*}
\text{div}\, \mathcal{M}_\A(x,u(x))&=\sum_{k=1}^n\dv_k\bigg(\sum_{j=1}^n\Aa(x,e_j)u(x)\otimes e_j\bigg)(e_k)\\&=\sum_{j=1}^n\dv_j(\Aa(x,e_j)u(x))=\sum_{j=1}^n(\dv_j\Aa(x,e_j))u(x))+\sum_{j=1}^n\Aa(x,e_j)\dv_ju(x)\\
&=\text{div}\, \Aa(x)u(x)-B(x)u(x)+\A u(x).
\end{align*}

In Fuglede's theory of flux extensions it is essential that Stokes theorem is allowed to hold on domains that are not necessarily $C^1$. This is due to the fact that $C^1$-domains have undesirable partition properties. We will therefore recall versions of Stokes theorem applicable to non-$C^1$ domains. If for example $\Om$ is a manifold with corners as, in the definition given in \cite[p. 415]{Lee} (see also the much more general version in \cite[Thm. 4.7]{Sauv}), then for any $u\in C^1(\Om,E)\cap C(\overline{\Om},F)$
\begin{align}\label{eq:Stokes!}
\int_{\Om}\mathscr{A}u(x)dx=\int_\Om(B(x)-\text{div}\, \mathbb{A}(x))u(x)dx+\int_{\dv \Om}\mathbb{A}(y,\nu(y))u(y)d\sigma(y),
\end{align}
In particular any polygonal domain is a manifold with corners.

For the purpose of flux extensions Fuglede defines Green domains as follows:
\begin{Def}[Green domains]
A bounded open set $U\subset \R^n$ is called a \emph{Green set} if there exists: 
\begin{itemize}
\item[]
\item[] A finite measure $\sigma\geq 0$, defined for all Borel sets $A\subset \dv U$. 
\item[]
\item[] A vector valued Baire function $\nu(x)=(\nu_1(x),...,\nu_n(x))$ defined $\sigma$-a.e. on $\dv U$ and of unit length with the property that 
\begin{align*}
\int_U\dv_ju(x)dx=\int_{\dv U}\nu_j(x)u(x)d\sigma(x). 
\end{align*}
for any $u\in C^1(\overline{U})$. 
\item[]
\end{itemize}
A connected Green set is called a Green domain. 
\end{Def}

As is shown in \cite{F2}, if $U$ is a Green domain, then $\sigma$ and $\nu$ are unique. In particular any domain so that $\dv U$ is a Lipschitz manifold is a Green domain with $\sigma=\mathcal{H}^{n-1}\lfloor \dv U$. In \cite[Thm. 1, p. 22]{F2} Fuglede gave the following implicit characterization of Green domains.
\begin{Thm}\label{thm:Green}
A bounded domain $U\subset \R^n$ is a Green domain if and only if there exist constants $M_1,M_2,...,M_n$ depending only on $U$ such that for every $u\in C^1(\overline{U})$
\begin{align*}
\bigg\vert \int_U\dv_ju(x)dx\bigg\vert \leq M_j\max_{x\in \dv U}\vert u(x)\vert. 
\end{align*}
\end{Thm}

While the notion of flux extensions does not require Stokes theorem to hold on more general domains with Lipschitz boundaries, for applications it is necessary to have a more general version of Stokes theorem available. Different essentially sharp versions of Stokes theorem for the divergence operator are given in \cite{MM,CT, MM2,CT2,CCT,Sil}. To state them we need some preliminary notions from geometric measure theory taken from \cite{AFP}. 

\begin{Def}
A set of density $\alpha\in [0,1]$ of a set $A\subset \R^n$ is defined by 
\begin{align*}
A^\alpha:=\bigg\{y\in \R^n: \lim_{r\to 0^+}\frac{\vert B_r(y)\cap A\vert}{\vert B_r(y)\vert}=\alpha\bigg\},
\end{align*}
where $\vert X \vert$ denotes the $n$-dimensional Lebesgue measure of measurable subset $X\subset \R^n$. In particular, $A^0$ is the measure theoretic exterior of a set $A$ and $A^1$ is the measure theoretic interior of set $A$. Furthermore the set 
\begin{align*}
\dv_s A=\R^n\setminus(A^0\cup A^1) 
\end{align*}
is called the \emph{measure theoretic} or \emph{essential boundary} of $A$. 
\end{Def}

\begin{Def}[Sets of finite perimeter]
Let $U\subset \R^n$ be a  Lebesgue measurable set. Let  $\Om \subset \R^n$ be open. The \emph{perimeter} of $U$ in $\Om$ is defined according to 
\begin{align}
P(U,\Om):=\sup\bigg\{ \int_{U}\text{div}\, \varphi(x)dx: \varphi\in C^1_0(\Om,\R^n), \,\,\, \Vert \varphi\Vert_{L^\infty}=1 \bigg\}. 
\end{align}
If $P(U,\Om)<+\infty$ we say that $U$ is a set of finite perimeter in $\Om$. In particular, $U$ has finite perimeter if and only if the characteristic function $\chi_U\in \text{BV}(\Om)$, where $\text{BV}(\Om)$ denotes the space of functions of bounded variations and $\vert \nabla \chi_U\vert(\Om)=P(U,\Om)$, where $\vert \nabla \chi_U\vert(\Om)$ denotes the total variation of the finite Radon measure $\nabla \chi_U$. 
\end{Def}

In particular any $C^1$-domain $U$ such that $\mathcal{H}^{n-1}(\dv U\cap \Om)<+\infty$ is a set of finite perimeter and the Gauss-Green formula
\begin{align}
 \int_{U}\text{div}\, \varphi(x)dx=-\int_{\dv U\cap \Om}\langle \nu_U(x),u(x)\rangle d\sigma(x)
\end{align}
holds where $\nu_{U}$ is the \emph{inner} unit normal of $\dv U$. Moreover, it follows from Theorem \ref{thm:Green} that any Green domain is a set of finite perimeter.

\begin{Def}[Reduced boundary and measure theoretic normal]
Let $U\subset \R^n$ be a Lebesgue measurable set in $\R^n$ and let $\Om$ be the largest open set such that $U$ is locally of finite perimeter in $\Om$. The \emph{reduced boundary} $\dv_\ast U\subset \dv U$ is the set of all points $x\in \text{supp}(\vert \nabla \chi_U\vert)$ such that the limit 
\begin{align*}
\nu_U(x):=\lim_{r\to 0^+}\frac{\nabla \chi_U(B_r(x))}{\vert \nabla \chi_U\vert(B_r(x))}.
\end{align*}
exists and satisfies $\vert \nu_U(x)\vert=1$. The function $\nu_U: \dv_\ast \Om\to S^{n-1}$ is called the \emph{measure theoretic inner unit normal} to $U$.
\end{Def}

By theorems of De Giorgi and Federer, see \cite[Thm. 359, p. 157, Thm. 3.61, p. 158]{AFP}, a set of finite perimeter $U$ in $\Om$ satisfies 
\begin{align*}
\dv_\ast U\subset \dv_s U, \quad \mathcal{H}^{n-1}(\Om \setminus (U^0\cup \dv_\ast U\cup U^1))=0, 
\end{align*}
and 
\begin{align*}
\vert \nabla \chi_U\vert=\mathcal{H}^{n-1}\lfloor \dv_\ast U.  
\end{align*}

Given an open set $\Om \subset \R^n$ and an aperture parameter $\kappa\in (0,+\infty)$ define the \emph{non-tangential approach regions}
\begin{align*}
\Gamma_\kappa(x):=\{y\in \Om: \vert y-x\vert<(1+\kappa)\text{dist}(y,\dv \Om)\}. 
\end{align*} 

\begin{Def}[Nontangentially accessible boundary]
Let $\Om\subset \R^n$ be a domain. The nontangentially accessible boundary of $\Om$ is defined as
\begin{align*}
\dv_{nta}\Om:=\{x\in \dv\Om: x\in \overline{\Gamma_\kappa(x)} \text{ for each }\kappa>0\}. 
\end{align*}
\end{Def}

Associated to a family of non-tangential approach regions is the \emph{non-tangential maximal function} of a measurable $u:\Om \to E$ at $\dv \Om$ defined according to 
\begin{align*}
(\mathcal{N}_\kappa u)(x):=\sup_{y\in \Gamma_\kappa(x)}\vert u(y)\vert
\end{align*} 
for $x\in \dv \Om$. We note that for some $x\in \dv \Om$ we may have $\Gamma_\kappa(x)=\varnothing$, for example if $\dv \Om$ has an outward pointing cusp at $x$. In this case $(\mathcal{N}_\kappa u)(x)$ is undefined. More generally, let $A\subset \Om$ be an arbitrary Lesbesgue measurable subset. Then the non-tangential maximal function with respect to $A$ is 
\begin{align*}
(\mathcal{N}_\kappa^A u)(x):=\sup_{y\in \Gamma_\kappa(x)\cap A}\vert u(y)\vert
\end{align*} 
for each $x\in \dv \Om$. 

\begin{Def}[Non-tangential limit]
Let $\Om \subset \R^n$ be open and let $x\in \dv \Om$ be such that $\Gamma_k(x)\neq \varnothing$. Then a measurable function $u:\Om \to E$ is said to have the $\kappa$-nontangetial limit $a$ at $x$ if for every $\eps>0$ there exists an $r>0$ such that $\vert u(y)-a\vert<\eps$ for $\mathscr{L}^n$-a.e. $y\in \Gamma_\kappa(x)\cap B_r(x)$. Whenever $u$ has a $\kappa$-nontangential limit at $x$, it will be denoted by $u\vert_{\dv \Om}^{\kappa-\text{n.t.}}(x).$. If the limit exits for every $\kappa'>0$ we simply say that $u$ has a \emph{nontangetial limit} denoted by $u\vert_{\dv \Om}^{\text{n.t.}}(x).$
\end{Def}

\begin{Def} 
Let $\Sigma\in \R^n$ be a closed set. $\Sigma$ is a called \emph{Ahlfors-David regular} if there exists constants $c,C\in (0,\infty)$ such that for each $x\in \Sigma$, $r\in (0,2\text{diam}(\Sigma))$ and $\rho>0$,
\begin{align*}
cr^{n-1}\leq \mathcal{H}^{n-1}(B(x,r)\cap) \Sigma),\quad  \mathcal{H}^{n-1}(B(x,\rho)\cap) \Sigma)\leq C\rho^{n-1}. 
\end{align*}
\end{Def}

Whenever $\dv \Om$ is Ahlfors-David regular, $\Om$ is a set of finite perimeter and the measure $\sigma=\mathcal{H}^{n-1}\lfloor \dv \Om$ is a \emph{doubling measure} and 
\begin{align*}
\sigma(\dv_s\Om \setminus \dv_{nta}\Om)=0. 
\end{align*}
For a proof see \cite{MM2}.

With all these notions defined we can now state the following version of the divergence theorem from \cite[Thm. 2.5]{MM}  and \cite[Thm. 1.2.1, p. 18]{MM2} (which for simplicity we state in a slightly more restrictive form for bounded domains)
\begin{Thm}\label{thm:Div1}
Let $\Om\subset \R^n$ be a bounded domain and assume that $\dv \Om$ is Ahlfors-David regular. Let $u\in L^1(\Om,\R^n)$ such that $\text{div}\, u\in L^1(\Om)$ in the sense of distributions. Assume further that $u\vert_{\dv \Om}^{\kappa-\text{n.t}}$ exists at $\sigma$-a.e. point on $\dv_{nta}\Om$ and $\mathcal{N}^\kappa u\in L^1(\dv \Om,\sigma;\R^n)$. Then for any $\kappa'>0$ $\kappa$-nontangential limit at $x$ exits $\sigma$-a.e on $\dv_{nta}\Om$ and is independent of $\kappa'$. This nontangential trace belongs to $L^1(\dv_\ast \Om,\sigma;\R^n)$. Furthermore, 
\begin{align*}
\int_{\Om}\text{div}\, u(x)dx=\int_{\dv_\ast \Om}\langle \nu(x),u\vert_{\dv \Om}^{\text{n.t.}}(x)\rangle d\sigma(x)
\end{align*}
where $\nu$ is the outward pointing measure theoretic unit normal. 
\end{Thm}

\begin{Def}[Divergence measure fields]
A vector field $u\in L^p(\Om,\R^n)$, $1\leq p\leq \infty$ is a called a \emph{divergence measure field} if $\text{div}\, F\in \mathcal{M}_b(\Om)$, where $\mathcal{M}_b(\Om)$ denotes the space of  Radon measures on $\Om$ with finite total variation. The space of all divergence measure fields is denoted by $\text{div}\,\mathcal{M}^p_b(\Om)$.
\end{Def}

For divergence measure fields the following version of the divergence theorem was proven independently by \cite{CT} and \cite{Sil} 

\begin{Thm}\label{thm:Div2}
Let $\Om$ be a set of finite perimeter and let $u\in \text{div}\,\mathcal{M}^\infty_b(\Om)$. Then for every $\phi\in C^1_c(\R^n)$ 
\begin{align*}
\int_{\Om^1}\phi(x)d \text{div}\, u(x)+\int_{\Om^1}\langle \nabla \phi(x),u(x)\rangle dx=-\int_{\dv^\ast \Om}\phi(x)\langle \text{Tr}_i u(x),\nu(x)\rangle d\sigma. 
\end{align*}
Here $\langle \text{Tr}_i u(x),\nu(x)\rangle$ is the interior normal trace as defined in \cite{CT}.
\end{Thm}

For divergence measure fields in $\mathcal{M}^p_b(\Om)$ with $p<\infty$ one instead have the following version of the divergence theorem.

\begin{Thm}\label{thm:Div4}
Let $\Om$ be open and let $U\subset \Om$ be an open set with Lipschitz boundary. Let $u\in \text{div}\,\mathcal{M}^p_b(\Om)$. Then for every $\phi\in C^0(\Om)$ such that $\nabla \phi\in L^{q}(\Om,\R^n)$, $1/p+1/q=1$, there exists a null set $\mathscr{N}\subset \R$ such that for every sequence $\{\eps_k\}_k\nsubseteq \mathscr{N}$ with $\eps_k\to 0^+$ as $k\to \infty$ 
\begin{align*}
&\int_{U}\phi(x)d \text{div}\, u(x)+\int_{U}\langle \nabla \phi(x),u(x)\rangle dx\\&=\lim_{k\to \infty}-\int_{\dv U}\frac{\phi(f_\eps(y))}{\vert \nabla \rho(f_\eps(y))}\langle u(f_\eps(y)),\nabla \rho(f_\eps(y)))\rangle \det(Df_\eps(y))d\sigma(y), 
\end{align*}
where $f_{\eps_k}: U \to U_{\eps_{k}}$ is a bilipschitz mapping to smooth domains $U_{\eps_k}$ and $\rho$ is a regularized distance to $U$. 
\end{Thm}
Theorem \ref{thm:Div4} is proven in \cite{CCT}. In fact an even more general version is given without assuming that $U$ is Lipschitz. See also the nice survey \cite{CT}. 

For an open set $\Om \subset \R^n$ denote by $\mathcal{M}(\Om)$ the space of Radon measures on $\Om$ and $\mathcal{E}'(\Om)$ the space of compactly supported distributions on $\Om$. By $\mathcal{E}'(\Om)+\mathcal{M}(\R^n)$ denote the subspace in $\mathcal{D}'(\Om)$ of all distributions that can be written on the form $\eta+\mu$ with $\eta\in \mathcal{E}'(\Om)$ and $\mu \in \mathcal{M}(\Om)$. In addition, $\mathcal{E}'(\Om)+\mathcal{M}(\R^n)\subset C^\infty_b(\Om)^\ast$, where $C^\infty_b(\Om)$ denotes the space of smooth and bounded functions on $\Om$, and $C^\infty_b(\Om)^\ast$ denotes the algebraic dual of $C^\infty_b(\Om)$, i.e., the space of all linear functionals on $C^\infty_b(\Om)$, not necessarily continuous. In \cite{MM2} an even more general version of Theorem \ref{thm:Div1} is proven. 

\begin{Thm}\label{thm:Div3}
Let $\Om\subset \R^n$ be a bounded domain and assume that $\dv \Om$ is Ahlfors-David regular and such that $\sigma=\mathcal{H}^{n-1}\lfloor \dv\Om$ is a doubling measure on $\dv \Om$. In particular $\Om$ is a set of finite perimeter and its geometric measure theoretic outward pointing unit normal $\nu$ is define $\sigma$-a.e on $\dv_\ast \Om$. Fix a $\kappa\in (0,\infty)$ and assume that the distributional vector field $u\in \mathcal{D}'(\Om,\R^n)$ satisfy the following conditions:
\begin{enumerate}
\item There exists a compact set $K\subset \Om$ such that $u\in L^1_{loc}(\Om\setminus K,\R^n)$ and $\mathcal{N}^{\Om \setminus K}_\kappa u(\dv\Om,\sigma;\R^n)\in L^1(\dv \Om,\R^n)$.
\item[]
\item The pointwise nontangential boundary trace $u^{\kappa-n.t.}\vert_{\dv \Om}$ exists $\sigma$-a.e. on $\dv_{nt}\Om$ and $\langle \nu, u^{\kappa-n.t.}\vert_{\dv \Om}\rangle \in L^1(\dv_s \Om,\sigma)$.
\item[]
\item $\text{div}\, u\in \mathcal{E}'(\Om)+\mathcal{M}(\Om)$.
\end{enumerate}
Then for any $\kappa'>0$ $u^{\kappa'-n.t}\vert_{\dv \Om}$ exists $\sigma$-a.e. on $\dv_{nta}\Om$ and is independent of $\kappa'$. Furthermore, the divergence theorem holds in the following form 
\begin{align*}
(\text{div}\,u,1)=\int_{\dv_\ast \Om}\langle \nu(x),u\vert_{\dv \Om}^{\text{n.t.}}(x)\rangle d\sigma(x),
\end{align*}
where $(\cdot,\cdot)$ denotes the duality pairing between $(C_b^\infty(\Om))^\ast$ and $C_b^\infty(\Om)$.
\end{Thm}
Theorem \ref{thm:Div3} is proven in \cite[Thm. 1.4.1, p. 38]{MM2}. The main difference between Theorem \ref{thm:Div2} and Theorem \ref{thm:Div3} besides the fact that $\text{div}\, u$ lies in different spaces is that they rely on different notions of traces. We now go back to the version of Stokes theorem for first order operator $\A$. Recall the identity 
\begin{align*}
\A u(x)&=\sum_{j=1}^nA_j(x)\dv_ju(x)+B(x)u(x)=\sum_{j=1}^n\dv_j(A_j(x)u(x))+(B(x)-\text{div}\,\Aa(x))\\&=\text{div}\,(\mathcal{M}_\A(x,u(x)))+(B(x)-\text{div}\,\Aa(x))u(x),
\end{align*}
where $\mathcal{M}_\A(x,u(x))$ is the matrix field defined by \eqref{eq:MatrixField}, and $\text{div}\,\mathcal{M}_\A(x,u(x))$ means the row-wise divergence. This identity was derived under the assumption that $u\in C^1(\Om,E)$. When $u\notin C^1(\Om,E)$ we will take the right hand-side as the definition of $\A$. By applying Theorem \ref{thm:Div3} component-wise for the $m$ rows $\mathcal{M}_\A(x,u(x))_1,...,\mathcal{M}_\A(x,u(x))_m$ of $\mathcal{M}_\A(x,u(x))$ (with $m=\text{dim}(F)$) and under the same assumptions about the non-tangential trace of the components of $u$ on $\dv \Om$ and the same assumptions about $\Om$ the divergence theorem implies for $u\in L^1(\Om,E)$
\begin{align}\label{eq:DualityDiv}
\sum_{j=1}^m(\text{div}\,\mathcal{M}_\A(x,u(x))_j,1 )f_j+\int_{\Om}(B(x)-\text{div}\,\Aa(x))u(x)dx=\int_{\dv \Om}\Aa(x,\nu(x))u\vert_{\dv_\ast \Om}^{\text{n.t.}}(x)d\sigma(x)
\end{align}
where $\{f_j\}_j$ is an ON-basis for $F$. We can identify 
\begin{align*}
\sum_{j=1}^m(\text{div}\,\mathcal{M}_\A(\cdot ,u)_j ,1)f_j+(B-\text{div}\,\Aa)u
\end{align*}
with an element in $\mathcal{E}'(\Om,F)+\mathcal{M}(\Om,F)\subset \mathcal{D}'(\Om,F)\cong F\otimes \mathcal{D}'(\Om)$ of $F$-valued distributions. An element in $v\in \mathcal{D}'(\Om,F)$ besides giving an $\R$-valued duality by acting on $C^\infty_c(\Om,F)$ also gives rise to an $F$-valued duality by acting on $C^\infty_c(\Om)$ according to 
\begin{align*}
(v,\phi)_F=\sum_{j=1}^\text{dim(F)}(f_j\otimes v_j)(\phi)=\sum_{j=1}^\text{dim(F)}f_j(v_j,\phi),
\end{align*}
where all $v_j\in \mathcal{D}'(\Om)$. The same applies when an element in $\mathcal{E}'(\Om,F)+\mathcal{M}(\Om,F)\subset \mathcal{D}'(\Om,F)\cong F\otimes \mathcal{D}'(\Om)$ acts on an element $\phi\in C^\infty_b(\Om)$. Thus, we interpret \eqref{eq:DualityDiv} as the $F$-valued duality pairing between $\mathcal{E}'(\Om,F)+\mathcal{M}(\Om,F)$ and $C^\infty_b(\Om)$ and \eqref{eq:DualityDiv} can be written
\begin{align}\label{eq:StokesA}
(\A u, 1)_F:=(\text{div}\, \mathcal{M}_\A(\cdot,u)+(B-\text{div}\,\Aa)u,1)_F=\int_{\dv \Om}\Aa(x,\nu(x))u\vert_{\dv_\ast \Om}^{\text{n.t.}}(x)d\sigma(x). 
\end{align}
The same as above can of course be done using Theorem \ref{thm:Div4} instead. In this case we say that $u\in \mathcal{M}^p_{\A,b}(\Om,E)$ if $\mathcal{M}_\A(\cdot ,u)\in F\otimes \mathcal{M}^p_b(\Om)$.

We conclude this section by mentioning Stokes theorem on smooth vector bundles $\mathcal{E}$ and $\mathcal{F}$ over a Riemannian base manifold $M$ with fibers $E_x$ and $F_x$. Consider a first order PDE operator
$\A: C^\infty(M,\mathcal{E})\to C^\infty(M,\mathcal{F})$, where by abuse of notation $C^{\infty}(M,\mathcal{E}/\mathcal{F})$ denotes smooth \emph{sections} rather than maps. For any $x\in M$ the principal symbol $\Aa(x,\nu(x)^\flat)\in \LL(T^\ast _xM, \LL(E_x,F_x))$ and the Green's formula
\begin{align*}
\int_{\dv \Om}\langle \A u(x),v(x)\rangle_{F_x} d\text{vol}_g(x)+\int_{\dv \Om}\langle u(x),\A v(x)\rangle_{E_x} d\text{vol}_g(x)=\int_{\dv \Om}\langle \Aa(x,\nu(x)^\sharp) u(x),v(x)\rangle_{F_x} d\text{vol}_{g_{\dv \Om}}(x)
\end{align*}
holds for any smooth domain $\Om \subset M$, and where $\text{vol}_g$ and $\text{vol}_{g_{\dv \Om}}$ are the induced Riemannian volume forms, and $\nu(x)^\flat$ is the outward pointing unit covector field and $\flat: T_x M\to T^\ast_xM$ is the musical isomorphism. In particular the integrals are real valued. However, the integral 
\begin{align}\label{eq:BadInt}
\int_{\dv \Om} \A u(x)d\text{vol}_g
\end{align}
is ill-defined because we are summing vectors over different fibers $F_x$ at different points of $\mathcal{F}$ which is not defined. We could use a local trivialization to define \eqref{eq:BadInt}, this however would not give rise to an intrinsic definition and would depend on the choice of local trivialization. If both $\mathcal{E}$ and $\mathcal{F}$ admit an isometric embedding into some euclidean space $\R^n$, then the integral \eqref{eq:BadInt} can be defined, but again the value of $\eqref{eq:BadInt}$ will depend on the choice of isometric embedding. Since we are only going to study flux extension on euclidean domains, these issues will not arise. It is however an interesting problem whether or not flux extensions can be defined on vector bundles.

%============NEW SUBSECTION=================================================================

\subsection{\sffamily Balance laws and a converse of Stokes}

In this section we will discuss balance laws and a converse of Stokes theorem. We will follow the discussion in Section 1.1 of \cite{Da}. Let $\Om\subset \R^n$ be a domain. A \emph{balance law} on $\Om$ postulates that the \emph{production} of an \emph{extensive quantity} in any (Lipschitz) subdomain $U\Subset \Om$ is balanced by the \emph{flux} of this quantity through the boundary $\dv U$ of $U$. The defining property of an extensive quantity is that both its production and flux are additive over disjoint subsets. The production of an extensive quantity over a domain is modelled by a vector valued Radon measure $\mathscr{P}$ on $\Om$ so that the production on a subset $U$ is given by $\mathscr{P}(U)$. The flux on the other hand is modelled by a countably additive set function $\mathscr{F}_U$ defined on Borel subsets of $\dv U$ such that the flux in and out of $\dv U$ is given by $\mathscr{F}_U(\dv U).$ A balance law is then simply the requirement that 
\begin{align}\label{eq:Balance}
\mathscr{P}(U)=\mathscr{F}_U(\dv U)
\end{align}
for every (Lipschitz) $U\Subset \Om$. From this point of view we see that Stokes theorem applied to a first order linear PDE operator can be view as a balance law by letting the production be given by 
\begin{align*}
\mathscr{P}_{\A u}(U)=\int_{U}\A u(x)-(B(x)-\di \Aa(x))u(x)dx
\end{align*}
and the flux be given by 
\begin{align*}
\mathscr{F}_{\Aa u}(\dv U)=\int_{\dv U}\Aa(x,\nu(x))u(x)d\sigma(y). 
\end{align*}
We will later see in Section \ref{sec:Flux} that balance laws and flux extensions are intimately related. 

On the other hand, given a balance law of the form \eqref{eq:Balance} what can be said about the relation between production and flux? Are they necessarily related via Stokes theorem for some suitable operator $\A$? Surprisingly it turns out the requirement coming from a balance law is quite restrictive and leads to a converse of Stokes theorem.

\begin{Thm}[Thm. 1.2.1 in \cite{Da}]
\label{thm:ConvStokes}
Consider a balance law \eqref{eq:Balance} on $\Om$ where $\mathscr{P}$ is a signed Radon measure. Assume furthermore that the flux is induced by a \emph{density flux function} $f_{\dv U}\in L^1(\dv U,\sigma)$ such that 
\begin{align*}
\mathscr{F}_U(\Gamma)=\int_{\Gamma}f_{\dv U}(x)d\sigma(x). 
\end{align*}
for any Borel subset $\Gamma \subset \dv U$. Assume further that $\vert f_{\dv U}(x)\vert \leq C$ for all $x\in \dv U$ and all Lipschitz domains $U\Subset \Om$. Then, \newline
(i) For each $\nu\in S^{n-1}$ there exists a measurable function $a_{\nu}\in \Om$ with the following properties: For any Lipschitz domain $U\Subset \Om$ and any point $x\in \dv U$ where outward pointing unit normal exits and equals $\nu$. Assume further that $x$ is a Lebesgue point of $f_{\dv U}$ relative to $\sigma=\mathscr{H}^{n-1}\lfloor \dv U$, and that the upper derivative of $\vert \mathscr{P}\vert$ at $x$ with respect to the Lesbegue measure on $\R^n$ is finite. Then 
\begin{align}
f_{\dv U}(x)=a_{\nu}(x).
\end{align} 
\newline (ii) There exists a vector field $V\in L^\infty(\Om,\R^n)$ such that for any fixed $\nu\in S^{n-1}$,
\begin{align}
a_{\nu}(x)=\langle \nu(x),V(x)\rangle, \quad \text{a.e. on $\Om$.}
\end{align} 
\newline ({iii}) The vector field $V$ satisfies the field equation 
\begin{align}
\di V(x)=\mathscr{P}
\end{align} 
in the sense of distributions on $\Om$. 
\end{Thm}

\begin{rem}
For the readers convenience, the upper derivative of a measure $\nu$ with respect to another measure $\mu$ is given by 
\begin{align*}
D^+_\mu \nu(x)=\limsup_{r\to 0^+}\frac{\nu(\overline{B_r(x)})}{\mu(\overline{B_r(x)})}. 
\end{align*}
\end{rem}

For a proof of Theorem \ref{thm:ConvStokes} we refer the reader to the proof of Theorem 1.2.1 in \cite{Da}.

%============NEW SUBSECTION=================================================================

\subsection{\sffamily Traces}

In relation to both Stokes theorem and Flugede's flux extensions it is important to know when functions belonging to some function space have traces, i.e., generalized notions of boundary values, and in what sense these traces exist. To that end, we will gather some known results in the literature that will be of useful later.

We begin by giving a proof of Seeley's trace theorem for first order operators presenting the nice proof of M. Taylor in \cite{TaylorN1}. 

\begin{Def}\label{def:Bmap}
Let $\Om\subset \R^n$ be a smooth domain. Define the boundary multiplier operator $\mathcal{B}_{\Aa}: C^\infty(\dv \Om,F)\to C^\infty(\dv \Om,E)$ through 
\begin{align*}
\mathcal{B}_{\Aa}v(x)=\Aa(x,\nu(x))^\ast v(x), 
\end{align*}
where $\nu(x)$ is the outward pointing unit normal of $\dv \Om$ at $x$. 
\end{Def}

\begin{Lem}
Assume that $\Om$ is $C^1$ domain and assume that the principle symbol $\Aa(x,\nu(x))$ of the operator $\A$ is injective for all $x\in \dv \Om$. Then the boundary multiplier operator $\mathcal{B}_{\Aa}$ extends to a surjective bounded map 
\begin{align*}
\mathcal{B}_{\Aa}: L^p(\dv \Om,F)\to L^p(\dv \Om,E)
\end{align*}
for all $1<p<\infty$. 
\end{Lem}

\begin{proof}
Since for each $x\in \dv \Om$, $\Aa(x,\nu(x)):E\to F$ is injective, $\Aa(x,\nu(x))^\ast :F\to E$ is surjective. Furthermore, since $\Om$ is $C^1$, $\nu:\dv \Om \to \R^n$ is continuous, thus $\Aa(x,\nu(x))^\ast \in C(\dv \Om, \LL(F,E))$. By subjectivity, $\Aa(x,\nu(x))^\ast$ has a right inverse that can be taken to be the Moore-Penrose pseudoinverse, $(\Aa(x,\nu(x))^\ast)^+$.  For any $v\in L^p(\dv \Om)$, define $u(x)=(\Aa(x,\nu(x))^\ast)^+v(x)$. Then $\Aa(x,\nu(x))^\ast u(x)=\Aa(x,\nu(x))^\ast(\Aa(x,\nu(x))^\ast)^+v(x)=v(x)$. Furthermore, $(\Aa(x,\nu(x))^\ast)^+$ is continuous. Thus $u\in L^p(\dv \Om,E)$.  
\end{proof}

\begin{Thm}[Seeley's trace theorem]
\label{thm:Trace}
Assume that $\Om$ is $C^1$ domain.  Let $\A$ be a first order partial differential operator whose symbol $\Aa(x,\xi)$ is injective for all $\xi\in S^{n-1}$ and $x\in \overline{\Om}$ and let
\begin{align*}
\text{dom}_{p,\Om}(\A)=\{\mathcal{D}'(\Om,E): u\in L^p(\Om,E),\,\, \A u\in L^p(\Om,F)\}.
\end{align*}
Assume furthermore that the matrix coefficients of $\A$ satisfies  $A_j\in C^1(\overline{\Om},\LL(E,F))$ for $j=1,2,...,n$ and $B\in C^0(\overline{\Om},\LL(E,F))$. Then there exists a bounded map $\text{Tr}:\text{dom}_{p,\Om}(\A)\to W^{-1/p,p}(\dv \Om, E)$ such that if $u\in C^1(\overline{\Om},E)$, 
\begin{align*}
\text{Tr}(u)=u\vert_{\dv \Om}
\end{align*}
and the Green identity 
\begin{align*}
(\text{Tr}(u),\mathcal{B}_{\Aa}\text{Tr}(v))=\int_{\Om}\langle \A u,v\rangle_F dx+\int_{\Om}\langle u, \A^\ast v\rangle_E dx
\end{align*}
holds for all $v\in W^{1,q}(\Om,F)$. 
 \end{Thm}

\begin{proof}
The Gauss-Green theorem implies that for all $u\in C^\infty(\overline{\Om},E)$ and all $v\in C^\infty(\overline{\Om},F)$
\begin{align*}
\int_{\Om}\langle \A u(x),v(x)\rangle_Fdx+\int_{\Om}\langle u(x),\A^\ast v(x)\rangle_Edx&=\int_{\dv \Om}\langle \Aa(x,\nu(x))u(x),v(x)\rangle_F d\sigma(x)\\
&=\int_{\dv \Om}\langle u(x),\Aa(x,\nu(x))^\ast v(x)\rangle d\sigma(x)\\
&=\int_{\dv \Om}\langle u(x),\mathcal{B}_{\Aa}v(x)\rangle d\sigma(x).
\end{align*}
Abstractly we can write this as
\begin{align}\label{eq:Duality}
(\A u,v)_{L^2(\Om,F)}+(u,\A^\ast v)_{L^2(\Om,E)}=(u,\mathcal{B}_{\Aa}v)_{L^2(\dv \Om,E)}
\end{align}

By an approximation argument this extends to $u\in W^{1,p}(\Om,E)$ and $v\in W^{1,q}(\Om,F)$, where we use that $\text{Tr}:W^{1,p}(\Om,V)\mapsto W^{1-1/p,p}(\dv \Om,V)\subset L^p(\dv \Om,V)$ for any real euclidean vector space $V$ (see \cite[Thm. 18.2]{DiB}). In addition the boundary multiplication operator extends to a bounded map $\mathcal{B}_{\Aa}: L^q(\dv \Om, F)\to  L^q(\dv \Om, E)$. Thus, 
\begin{align}\label{def:Trace}
(u,\mathcal{B}_{\Aa}v):=\int_{\dv \Om}\langle u(x),\mathcal{B}_{\Aa}v(x)\rangle_E d\sigma(x).
\end{align}
gives a well-defined duality $W^{1-1/p,p}(\dv \Om,E)\times W^{1-1/q,q}(\dv \Om,E)\to \R$. We will now use the duality \eqref{eq:Duality} to define $\text{Tr}(u)$ for any $u$ in $\text{dom}_{p,\Om}(\A)$. By \cite[Thm. 18.2]{DiB}), the trace map $\text{Tr}: W^{1,p}(\Om,E)\to W^{1-1/p,p}(\dv \Om,E)$ is surjective. Hence for any $\varphi\in W^{1-1/p,p}(\dv \Om,E)$ we can find a $w\in W^{1,p}(\Om,E)$ such that $\text{Tr}(w)=\varphi$. Using that the boundary multiplication operator is surjective there exists a $\phi\in W^{1-1/q,q}(\dv \Om,F)$ such that $\varphi(x)=\Aa^\ast(\nu(x))\phi(x)$. We now \emph{define} the trace $\psi=\text{Tr}(u)\in W^{-1/p,p}(\dv \Om,E)$ of $u\in \text{dom}_{p,\Om}(\A_w)$ by
\begin{align*}
( \psi, \mathcal{B}_{\Aa}v):=(\A u,v)-(u,\A^\ast v)
\end{align*}
for all $v \in W^{1,q}(\dv \Om,F)$ such that the $\text{Tr}(v)=\phi$ in the Sobolev sense. To give a well-defined duality we need to check that this definition is independent of the choice of $v$. Therefore let $v_1,v_2\in W^{1,q}(\dv \Om,F)$ be such that $\text{Tr}(v_1)=\text{Tr}(v_2)$. We need to show that 
\begin{align*}
(\A u,v_1)-(u,\A^\ast v_1)=(\A u,v_2)-(u,\A^\ast v_2).
\end{align*}
or equivalently that $w=v_1-v_2$ satisfy 
\begin{align*}
(\A u,w)-(u,\A^\ast _w)=0
\end{align*}
as long as $u\in \text{dom}_{p,\Om}(\A_w)$. This is clear however since it is true for any $w\in C^{\infty}_0(\Om,F)$ and since $C^{\infty}_0(\Om,F)$ is dense in $W^{1,q}_0(\Om,F)$ the general case follows by approximation. Finally, we show that $\text{Tr}: \text{dom}_{p,\Om}(\A)\to W^{1-1/p,p}(\dv \Om,E)$ is bounded. Since $\Om$ is a Sobolev extension domain there exists a bounded extension operator $E: W^{1-1/q,q}(\dv \Om,F)\to W^{1,q}(\Om,F)$. Consequently, for any $\varphi\in W^{1-1/p,p}(\dv \Om,E)$ 
\begin{align*}
\vert (\text{Tr}(u),\mathcal{B}_{\Aa}v)\vert\leq \vert (\A u,E\phi)\vert +\vert (u,\A^\ast E\phi )\vert 
\end{align*}
\end{proof}

There is another more constructive way to define the trace of $u\in \text{dom}_{p,\Om}(\A)$ also found in \cite{TaylorN1}. This view point also fits naturally with the divergence theorem, Theorem \ref{thm:Div4} as well as Rauch's continuity theorem for traces, Theorem \ref{thm:RauchCont} below. Let $V$ be a smooth vector field in $\R^n$, such that $V$ intersects $\dv \Om$ transversely and that 
\begin{align*}
\langle \nu(x),V(x)\rangle <0, \quad x\in \dv \Om
\end{align*}
where $\nu(x)$ is the outward pointing unit normal of $\dv \Om$ at $x$. Let $\Phi_t$ denote the flow generated by $V$ and consider the one-parameter family of domains 
\begin{align*}
\overline{\Om_t}=\Phi_t(\overline{\Om}), \quad t\geq 0,
\end{align*}
and so for all $t\in (0,\delta)$ for $\delta>0$ sufficiently small, we have $\overline{\Om}_{t_2}\Subset \Om_{t_1}$ whenever $t_2>t_1$ and $\Om_0=\Om$. By local elliptic regularity, $u\in W^{1,p}_{loc}(\Om,E)$ which implies that $u\vert_{\Om_t}\in W^{1,p}(\Om_t,E)$. Thus $u\vert_{\dv \Om_t}\in W^{1-1/p,p}$ for all $t>0$. In particular for each $v\in W^{1,p}(\Om,F)$ and each $t>0$
\begin{align}\label{eq:GreenTrace}
\int_{\Om_t}\langle \A u(x),v(x)\rangle_F dx+\int_{\Om_t}\langle u(x),\A^\ast v(x)\rangle_Edx&=\int_{\dv \Om_t}\langle u(x),\Aa(x,\nu(x))^\ast v(x)\rangle d\sigma(x).
\end{align}
Let 
\begin{align*}
u_t=\Phi_t^\ast (u\vert_{\Om_t})\in , \quad \psi_t=\Phi_t^\ast (u_t\vert_{\dv \Om})\in W^{1-1/p,p}(\dv \Om,E).
\end{align*}
By \eqref{eq:GreenTrace}
\begin{align*}
\Vert \psi_t\Vert_{W^{-1/p,p}}\leq C(\Vert u\Vert_{L^p}+\Vert \A u\Vert_{L^p})
\end{align*}
for $t>0$, and with $C$ independent of $t$. Furthermore, 
\begin{align*}
\lim_{t\to 0+}&\int_{\Om_t}\langle \A u(x),v(x)\rangle_F dx+\int_{\Om_t}\langle u(x),\A^\ast v(x)\rangle_Edx\\&=\int_{\Om}\langle \A u(x),v(x)\rangle_F dx+\int_{\Om}\langle u(x),\A^\ast v(x)\rangle_Edx=(\psi,\mathcal{B}_{\Aa}\phi)
\end{align*}
by the definition of trace according to \eqref{def:Trace} for each $\phi$ with $v\vert_{\dv \Om}=\phi$. This shows that $\psi_t$ converges weak${}^*$ in $W^{-1/p,1}$ to $\text{Tr}(u)$ as $t\to 0^+$.

\begin{rem}
Note that the trace map $\text{Tr}$ on $\text{dom}_{p,\Om}(\A)$ is not surjective in general. This is for example the case when $\A$ is a Dirac type operator. In this case the image of $\text{Tr}$ can be characterized by the spectrum of an adapted operator of $\A$ on the boundary. For details see \cite{BärBal}. On the other hand in the case when $\A =D$ is the total derivative then the trace map is indeed surjective, see \cite[Thm. 18.2]{DiB}. More generally this holds true whenever $\A$ is a $\C$-elliptic operator. For details see the discussion in Section \ref{subsec:StokesAll}.
\end{rem}

In fact the assumptions of Theorem \ref{thm:Trace} can be relaxed significantly. This is done in \cite{R1}. Before we state Rauch's theorem we define:
\begin{align}
\label{eq:Space1}
\mathscr{K}_\A&:=\{u\in L^2(\Om,E): \A u \in W^{1,2}(\Om,E)'\}\\
\mathscr{H}_\A&:=\text{dom}_{2,\Om}(\A)=\{u\in L^2(\Om,E): \A u \in L^{2}(\Om,E)\}\label{eq:Space2}
\end{align}
equipped with the graph norms 
\begin{align}
\label{eq:Norm1}
\Vert u\Vert_{\mathscr{K}_\A}&=\Vert u\Vert_2+\Vert \A u\Vert_{W^{1,2}(\Om,F)'}\\
\label{eq:Norm1}
\Vert u\Vert_{\mathscr{H}_\A}&=\Vert u\Vert_2+\Vert \A u\Vert_{L^2(\Om,F)}
\end{align}
where $W^{1,2}(\Om,E)'$ is the dual space of $W^{1,2}(\Om,E)$. Note that $\mathscr{K}_\A$ is only a slightly smaller space than $L^2(\Om,E)$ since automatically $\A u\in \dot W^{1,2}(\Om, F)'$, where 
\begin{align*}
 \dot W^{1,2}(\Om, F)=\{u\in \mathcal{D}'(\Om,E): \A u\in L^2(\Om,F)\}
\end{align*}
is the homogeneous Sobolev space. It is shown in \cite[Prop. 1]{R1} that if $\Om$ is a $C^1$ domain, then $\mathscr{K}_\A$ and $\mathscr{H}_\A$ are Hilbert spaces and $C^1(\overline{\Om},E)$ is dense in both.

\begin{Thm}[Rauch]
\label{thm:Rauch1}
Let $\Om \subset \R^n$ be a $C^1$-domain and consider the first order operator 
\begin{align*}
\A u(x)=\sum_{j=1}^nA_j(x)\dv_ju(x)+B(x)u(x),
\end{align*}
where $A_j\in \text{Lip}(\overline{\Om}, \LL(E,F)$ for $j=1,2,...,n$ and $B\in L^\infty(\Om,\LL(E,F))$. Then the map 
\begin{align*}
C^1(\overline{\Om},E)\ni u\mapsto \Aa(x,\nu(x))u(x)
\end{align*}
extends to a continuous linear map
\begin{align*}
\text{Tr}: \mathscr{K}_\A \to W^{-1/2,2}(\dv \Om,\sigma;E)
\end{align*}
and Green's identity 
\begin{align}\label{eq:GreenRauch}
\int_{\dv \Om}\langle \Aa(x,\nu(x))u(x),v(x)\rangle d\sigma(x)=\int_{\Om}\langle \A u,v\rangle_F dx+\int_{\Om}\langle u, \A^\ast v\rangle_E dx
\end{align}
holds. Furthermore the bilinear map 
\begin{align*}
C^1(\overline{\Om},E)\times C^1(\overline{\Om},F)\ni (u,v)\mapsto \langle \Aa(x,\nu(x)u(x),v(x)\rangle
\end{align*}
extends to a continuous bilinear map $\mathscr{H}_\A\times \mathscr{H}_{\A^\ast}\to \text{Lip}(\dv \Om)'$ and again \eqref{eq:GreenRauch} holds for any $(u,v)\in \mathscr{H}_\A\times \mathscr{H}_{\A^\ast}$. 
\end{Thm}

In Theorem \ref{thm:Rauch1} the boundary integral in Green's identity is understood the action of the distribution $\langle \Aa(\nu)u,v \rangle$ acting on the Lipschitz continuous function $1$ on $\dv U$, i.e., 
\begin{align*}
\int_{\dv U}\langle \Aa(\nu(x))u(x),v(x) \rangle d\sigma(x):=(\langle \Aa(\nu)u,v \rangle,1).
\end{align*}
Also note that $\text{Lip}(\dv \Om)'\subset \mathcal{D}'(\dv \Om)$ does not have a useful elementary characterization.

Finally, Rauch also proved in \cite{R1} the following continuity result with respect to traces of near by surfaces. 
\begin{Thm}\label{thm:RauchCont}
Let $\Om$ be a $C^1$-domain and let $\Om_\eps\Subset \Om$ be monotonically increasing sequence of $C^1$-domains generated be a one-parameter family of diffeomorphism $\{\rho_t\}_{t\geq 0}$, such that $\Om_t=\rho_t(\Om)$, $\rho_0(\Om)=\Om$ and $\lim_{t\to 0^+}\rho_\eps=\text{id}$ in $C^1$ and $\{\dv \Om_\eps\}_\eps$ form a $C^1$-foliation close to $\dv \Om$. Then for any $u \in \mathscr{K}_\A(\Om)$ the map 
\begin{align*}
[0,\delta)\ni t\mapsto \Aa(\nu)u\vert_{\dv \Om_t}\in W^{-1/2}(\dv \Om_t,\sigma;E)
\end{align*}
for some $\delta>0$ is continuous. Similarly, if $u\in \mathscr{H}_\A(\Om)$, $v\in  \mathscr{H}_{\A^\ast}(\Om)$ then 
\begin{align*}
[0,\delta)\ni t\mapsto \langle \Aa(\nu)u,v\rangle \vert_{\dv \Om_t}\in Lip(\dv \Om_t,\sigma)'
\end{align*}
is continuous. 
\end{Thm}

%============NEW SUBSECTION=================================================================

\section{\sffamily Weak and strong extensions of partial differential operators}\label{sec:WeakStrong}

Here we follow the exposition in \cite{F2}, see also Hörmander in \cite{H3} for which we refer for a more complete treatment. For simplicity of exposition and relevance to the present paper we will only consider first order operators. 

Let $X$ and $Y$ be Banach spaces. The cartesian product $X\times Y$ is also a Banach space equipped with the norm $\Vert (x,y)\Vert_{X\times Y}=\Vert x\Vert_X+\Vert y\Vert_Y$.  We denote by $X^\ast$ and $Y^\ast$ the dual spaces of $X$ and $Y$ respectively. Let $V$ be a linear subspace of some Banach space $Z$. The annihilator space $V^0$ of $V$ is defined as all functionals $z'\in Z^\ast$ such that
\begin{align*}
\langle z',v\rangle=0
\end{align*}
 for all $z\in Z^\ast$ and all $v\in V$. Let $\overline{V}$ denote the closure of $V$ in $Z$. Then in particular $V^0=(\overline{V})^0$ and $V$ is everywhere dense in $Z$ if and only if $V^0=\{0\}$.

Let $V\subset X$ be a linear subspace not necessarily closed. Let $T:V\to Y$ be a linear operator and let 
\begin{align*}
\text{gr}(T):=\{(u,Tu)\in V\times Y\}\subset X\times Y
\end{align*}
denote its graph. We will call $V$ the domain of $T$ and denote it by $\text{dom}(T)$. An operator $T$ is called closed if its graph is a closed subspace of $X\times Y$. If $T_1:V_1\subset X\to Y$ and $T_2:V_2\subset X\to Y$ are operators we say that $T_1\subsetneq T_2$ if $\text{gr}(T_1)\subset \text{gr}(T_2)$ as sets and we call $T_1$ the restriction of $T_2$, (obviously we have $\text{dom}(T_1)\subset \text{dom}(T_2)$) and $T_2$ an extension of $T_1$. An operator $T$ is called pre-closed if admits at least one closed extension $\overline{T}$. The domain $\text{dom}(T)$ of a closed operator $T$ equipped with the graph norm $\Vert x\Vert_T:=\Vert x\Vert_X+\Vert Tx\Vert_Y$ is complete (see \cite[Prop. 1.4, p. 6]{KSch}). Furthermore, a densely defined linear operator $T:X\to Y$ is closable if and only if for any sequence $\{x_n\}_n\subset \text{dom}(T)\subset X$ such that $\lim_{n\to \infty}x_n=0$ in $X$ and $\lim_{n\to \infty}Tx_n=y$ in $Y$, then $y=0$. For a proof of this statement see \cite[Prop. 1.5, p. 6]{KSch}. We will now consider the case when $X=L^p(\Om,E)$ and $Y=L^p(\Om,F)$ where $\Om \subset \R^n$ is a domain and $E$ and $F$ are finite dimensional real euclidean vector spaces equipped with an euclidean inner product $\langle \cdot,\cdot\rangle$ or finite dimensional complex hermitian vector spaces equipped with an hermitian inner product $(\cdot,\cdot)$. Moreover, using the inner products we will always identify $E^\ast\cong E$ and $F^\ast \cong F$. We let $C^\infty_0(\Om,E)$ denote the space of compactly supported smooth functions with values in $E$ and we let $C^\infty(\overline{\Om},E)$ denote the space of smooth functions with values in $E$ such that every $u\subset C^\infty(\overline{\Om},E)$ is the restriction of some smooth function $\overline{u}\in C^\infty(U,E)$ where $\Om \Subset U$ for some open set $U$. Let $\mathscr{A}$ be a first order partial differential operator of the form \eqref{eq:OpA}. Let $\mathscr{A}^\ast$ denote its formal adjoint given by \eqref{eq:OpAdj}.
Then for any smooth $U\Subset \Om$ and any $u\in C^1(\Om,E)$ and $\phi \in C^1(\Om,F)$ the Lagrange-Green identity 
\begin{align*}
\int_U\langle \mathscr{A}u(x),\phi(x)\rangle_Fdx-\int_U\langle u(x),\mathscr{A}^\ast \phi(x)\rangle_Edx&=\int_{\dv U}\langle \mathbb{A}(y,\nu(y))u(y),\phi(y)\rangle_F d\sigma(y)
\end{align*}
holds, which again is a consequence of partial integration and ultimately Stokes theorem. Moreover, if either $u$ or $\phi$ has compact support in $\Om$, then 
\begin{align*}
\int_U\langle \mathscr{A}u(x),\phi (x)\rangle_Fdx=\int_U\langle u(x),\mathscr{A}^\ast \phi(x)\rangle_Edx
\end{align*}

\begin{Def}[Minimal extensions]
The \emph{minimal extension} $\mathscr{A}_0$ of an operator $\mathscr{A}$ is defined to be the closure the graph $\{(u,\A u): u\in C_0^\infty(\Om,E)\}$ in $ L^p(\Om,E)\times L^p(\Om,F)$. The minimal extension is also called the \emph{minimal operator}. We write $\text{dom}_{p,\Om}(\A_0)$ for the domain of $\A_0$, and $\text{gr}_{p,\Om}(\A_0)$ for the graph of the extension.
\end{Def}

\begin{ex}
Let $\A=D$, the total derivative, so that $Du(x)\in E\otimes \R^n=F$, for $u\in C^\infty(\Om,E)$. Then $\text{dom}_{p,\Om}(D_0)=W^{1,p}_{0}(\Om,E)$.
\end{ex}

\begin{Def}[Weak extension]
The weak extension $\mathscr{A}_w$ of an operator $\mathscr{A}$ is defined to be the set of all pairs $(u,v)\in L^p(\Om,E)\times L^p(\Om,F)$ such that the Lagrange-Green formula 
\begin{align*}
\int_U\langle v(x),\phi(x)\rangle_Fdx=\int_U\langle u(x),\mathscr{A}^\ast \phi(x)\rangle_Edx
\end{align*}
for all $\phi\in C^1_0(\Om,E)$. The weak extension is also called the \emph{maximal operator} (not to be confused with the Hardy-Littlewood maximal operator). We write $\text{dom}_{p,\Om}(\A_w)$ for the domain of $\A_w$, and $\text{gr}_{p,\Om}(\A_w)$ for the graph of the extension.
 \end{Def}

\begin{ex}
As in the case of the minimal extension we take $\A=D$, the total derivative. Then $\text{dom}_{p,\Om}(D_w)=W^{1,p}(\Om,E)=\{u\in \mathcal{D}'(\Om,E): u\in L^p(\Om,E), \,\,\, Du\in L^p(\Om, E\otimes \R^n)\}$.
\end{ex}

The notion of weak extension is of course also intimately related to the very important notion of weak solution of $\mathscr{A}u=0$, namely:
$u\in \text{dom}(\mathscr{A}_w)\subset L^p(\Om,E)$ is a \emph{weak solution} if for every $\phi\in C^1_0(\Om,E)$ we have 
\begin{align*}
\int_U\langle u(x),\mathscr{A}^\ast \phi(x)\rangle_Edx=0.
\end{align*} 

Moreover, there is a direct connection between the minimal and maximal extension. Namely that $(\A_w)^\ast=\A_0$, where $(\A_w)^\ast$ is the adjoint of $\A_w$, so that in particular we have for all $u\in \text{dom}(\A_0)$ and all $v\in \text{dom}((\A_w)^\ast)$
\begin{align*}
(\A_0 u,v)=(u,(\A_w)^\ast v). 
\end{align*}

\begin{Def}[Locally strong extension]
Let $V=\{u\in C^1(\Om): u\in L^p(\Om,E), \mathscr{A}u\in L^p(\Om,F)\}$. The \emph{locally strong extension} $\mathscr{A}_{ls}$ of an operator $\mathscr{A}$ is the closure of $\mathscr{A}$ with respect to convergence in $L^p_{loc}(\Om,E)\times L^p_{loc}(\Om,F)$, i.e. the graph of $\mathscr{A}_{ls}$ consists of all pairs $(u,v)\in L^p(\Om,E)\times L^p(\Om,F)$ such that for every open set $U\Subset \Om$ there exists sequence $\{u_\eps\}_\eps\subset C^\infty_0(\Om,E)$ such that 
\begin{align*}
\lim_{\eps\to 0}\Vert u-u_\eps\Vert_{L^p(U,E)}=0,\quad \lim_{\eps\to 0}\Vert v-\mathscr{A}u_\eps\Vert_{L^p(U,F)}=0
\end{align*}
We write $\text{dom}_{p,\Om}(\A_{ls})$ for the domain of $\A_{ls}$, and $\text{gr}_{p,\Om}(\A_{ls})$ for the graph of the extension.
\end{Def}

\begin{rem}
In the literature there does not seem to be a consistent standard terminology for various different types of strong extensions. The reader should be aware that sometimes locally strong extensions are just called strong extensions as for example in \cite{KF}.
\end{rem}

\begin{Def}[Strong extension]
Let $V=\{u\in C^1(\Om): u\in L^p(\Om,E), \mathscr{A}u\in L^p(\Om,F)\}$. The \emph{strong extension} $\mathscr{A}_s$ of an operator $\mathscr{A}$ is the closure of $\mathscr{A}$ with respect to the domain $V$, i.e., $(u,v)\in \text{gr}(\mathscr{A}_s)$ if there exists a sequence $\{u_n\}_n\subset V$ such that $\lim_{n\to \infty}u_n=u$ and such that $\lim_{n\to \infty}\mathscr{A}u_n=v\in L^p(\Om,F)$. We write $\text{dom}_{p,\Om}(\A_s)$ for the domain of $\A_s$ and $\text{gr}_{p,\Om}(\A_s)$ for the graph of the extension. 
\end{Def}

\begin{ex}
Let again $\A=D$ be the total derivative acting on $C^\infty(\Om,E)$. It is shown in \cite[Thm. 2, Ch. 4.2, p. 125]{EG} that $D_s=D_w$ and that consequently $\text{dom}_{\Om,p}(D_s)=W^{1,p}(\Om,E)$ for any domain $\Om$ and any $1\leq p<\infty$. 
\end{ex}

\begin{Def}[Very strong extension]
Let $V=\{u\in C^1(\overline{\Om}): u\in L^p(\Om,E), \mathscr{A}u\in L^p(\Om,F)\}$. The \emph{very strong extension} $\mathscr{A}_S$ of an operator $\mathscr{A}$ is the closure of $\mathscr{A}$ with respect to the domain $V$, i.e., $(u,v)\in \text{gr}(\mathscr{A}_S)$ if there exists a sequence $\{u_n\}_n\subset V$ such that $\lim_{n\to \infty}u_n=u$ and such that $\lim_{n\to \infty}\mathscr{A}u_n=v\in L^p(\Om,F)$. We write $\text{dom}_{p,\Om}(\A_S)$ for the domain of $\A_S$, and $\text{gr}_{p,\Om}(\A_S)$ for the graph of the extension.
\end{Def}

Clearly, we have 
\begin{align*}
\text{gr}(\A_0) \subset \text{gr}(\A_S)\subset \text{gr}(\A_s)\subset \text{gr}(\A_{ls})\subset \text{gr}(\A_w).
\end{align*}
One can ask under which conditions we have 
\begin{align*}
\A_S=\A_s=\A_{ls}=\A_w\,\,\, \text{or}\,\,\, \A_{ls}=\A_w, \,\,\, \text{or}\,\,\, \A_s=\A_w \,\,\, \text{or}\,\,\,  \A_S=\A_w?
\end{align*}
Here equality means that the graphs are equal. The answer to these question depends both on the type of operator $\A$ and properties of the domain $\Om$. Firstly, in the case of the total derivative and Sobolev spaces we have by \cite[Thm. 3, Ch. 4.2, p. 127]{EG} that if $\Om$ is a bounded domain whose boundary is locally the graph of a Lipschitz function and $1\leq p<\infty$, then $D_S=D_w$, that is for any $u\in W^{1,p}(\Om,E)$ there exists a sequence $\{u_j\}_j\subset C^\infty(\overline{\Om},E)\cap W^{1,p}(\Om,E)$ such that $\lim_{j\to \infty}\Vert u_j-u\Vert_{W^{1,p}(\Om,E)}=0$.

Secondly, Hörmander proved in \cite{H1} that for any constant coefficient scalar operator there exists a domain $\Om\subset \R^n$ such that $\A_S\neq \A_w$. Since the example is very illuminating we will give it here adapted to the case when $\A$ is a constant coefficient homogeneous first order system. 
\begin{ex}
We may assume that $A_n=\Aa(e_n)\neq 0$. Let $\widetilde{\Om}=B_2(0)\subset \R^n$. Let $\gamma=B_1(0)\cap \{x\in \R^n:x_n=0\}$ and set $\Om=\widetilde{\Om}\setminus \gamma$. We note that $L^p(\widetilde{\Om},E)\vert_{\Om}=L^p(\Om,E)$ since $\vert \gamma\vert=0$. Furthermore, $C^\infty(\overline{\widetilde{\Om}},E)\vert_{\Om}=C^\infty(\overline{\Om},E)$ since any $u\in C^\infty(\overline{\Om},E)$ extends to $C^\infty(\overline{\widetilde{\Om}},E)$. On the other hand $C^\infty(\widetilde{\Om},E)\vert_{\Om}\neq C^\infty(\Om,E)$. To see this, let $\varphi\in C^\infty_0(\widetilde{\Om})$ be nonnegative, have support in $B_{1/2}(0)$ and be equal to 1 in a neighbourhood of $0$. Let
\begin{equation*}
  \phi(x) = \left\{
    \begin{array}{rl}
      \varphi(x) & \text{if } x_n>0 ,\\
      0 & \text{if } x_n < 0.
     \end{array} \right.
\end{equation*}
Since by assumption $A_n\neq 0$, $\text{ker}(A_n)\neq E$. Therefore we may choose $u\in C^\infty_0(\widetilde{\Om},E)$ such that $u(x)$ is a constant vector on $B_1(0)$ not in $\text{ker}(A_n)$ and define $u_\phi(x)=\phi(x)u(x)$. Then $u_\phi \in C^\infty_0(\widetilde{\Om},E)$ but $u_\phi \notin C^\infty_0(\Om,E)$, in fact $u_\phi\in \text{dom}_{p,\Om}(\A_w)$. Let $u_\phi^+(y,0)=\lim_{x_n\to 0^+}u_\phi(x)$ and $u_\phi^-(y,0)=\lim_{x_n\to 0^-}u_\phi(x)$ for $(y,0)\in \gamma$. We now choose $w\in C^\infty_0(\widetilde{\Om},F)$ to be equal to the constant vector $A_n u(x)$ on $B_1(0)$. Then, by Stokes theorem 
\begin{align*}
&\int_{\widetilde{\Om}}\langle \A u_{\phi}(x),w(x)\rangle dx-\int_{\widetilde{\Om}}\langle  u_{\phi}(x),\A^\ast w(x)\rangle dx=\int_{\gamma}\langle \Aa(e_n) (u_{\phi}^+(y,0)-u_{\phi}^-(y,0),w(y)\rangle dy\\
&=\int_{\gamma}\langle A_n \varphi(y,0)u(y,0),w(y,0)\rangle dy=\int_{\gamma}\varphi(y,0) \langle A_n u(y,0),A_n u(y,0)\rangle dy>0.
\end{align*}
Hence $u_\phi\notin \text{dom}_{p,\widetilde{\Om}}(\A_S)$ and so $\A_S\neq \A_w$. 
\end{ex}

One can ask now ask if in fact there are domains $\Om$ and operators $\A$ such that $\A_s\neq \A_w$. An example where $\A_w\neq \A_s$ was provided by S. Schwarz in \cite{Sch}. Here $\A=\dv_x\dv_y$ is a non-elliptic operator in the plane and $\Om=B_1(0)\setminus\{0\}$. 
\begin{Thm}[Thm. 1 in \cite{Sch}]
Let $\Om=B_1(0)\setminus\{0\}\subset \R^2$ and let $\A=\dv_1\dv_2$ be the wave operator. Then $\A_s\neq \A_w$.
\end{Thm}

In fact, a more general statement is shown in \cite{Sch}, see \cite[Thm. 4]{Sch}. Both \cite[Thm. 1]{Sch} and \cite[Thm. 4]{Sch} however require the order of the operator to be larger than one. In Theorem \ref{thm:sW} we show that for first order operators of the form \eqref{eq:OpA} we always have $\A_s=\A_w$. Before we show this we will first prove Friedrichs equivalence theorem showing that $\A_{ls}=\A_w$, a truly fundamental result in PDE theory. We first need to introduce Friedrichs mollifiers.

\begin{Def}[Friedrichs Mollifier]
\label{def:FMolli}
Let $\phi\in C^{\infty}(\R^n)$ satisfy the following properties:
\begin{itemize}
\item[]
\item[(i)] $\displaystyle \int_{\R^n}\phi(x)dx=1$
\item[]
\item[(ii)] $\phi(x)\geq 0$
\item[]
\item[(iii)] $\text{supp}(\phi)\subset B_1(0)$. 
\item[]
\end{itemize}
Set $\phi_\eps(x)=\eps^{-n}\phi(x/\eps)$ and define the operator 
\begin{align*}
(\Phi_\eps u)(x)=\int_{\Om}\phi_{\eps}(x-y)u(y)dy. 
\end{align*}
for $u\in L^1_{loc}(\Om,E)$. Then the operator $\Phi^\eps$ a \emph{Friedrich mollifier}. 
\end{Def}

\begin{ex}
A common choice of Friedrich mollifier is to chose
\begin{align}\label{eq:Molli}
\phi(x)=c\exp\bigg(\frac{1}{1-\vert x\vert^2}\bigg)\chi_{B_1(0)}(x),
\end{align}
with the constant chosen to satisfy the normalisation condition $(i)$. Note that 
\begin{align}\label{eq:Molli2}
\nabla \phi(x)=-\frac{2cx}{(1-\vert x\vert^2)^2}\exp\bigg(\frac{1}{1-\vert x\vert^2}\bigg)\chi_{B_1(0)}(x)=-\frac{2x}{(1-\vert x\vert^2)^2}\phi(x)
\end{align}
\end{ex}

In the original work of Friedrichs in \cite{KF} as well as that of Fuglede in \cite{F2} the mollifier used is of the form 
\begin{align*}
j_\eps(x)=\prod_{j=1}\phi_{\eps}(x_j)
\end{align*}
such that $j_{\eps}$ is supported in a cube of side length proportional to $\eps$. The different choices of mollifiers are however immaterial.  

The mollifier has a many useful properties, some of which are summarised in the theorem below.  
\begin{Thm}
Let $\Om_\eps=\{x\in \Om: \text{dist}(x,\dv \Om)>\eps\}$ and let $\Phi_\eps$ be a Friedrich mollifier. Let  $u\in L^1_{loc}(\Om,E)$ and set $u_\eps=\Phi_\eps u$. 
\begin{itemize}
\item[(i)]For each $\eps>0$, $u_\eps\in C^\infty(\Om_\eps,E)$.
\item[(ii)] If $u\in C(\Om,E)$ then $\lim_{\eps \to 0^+}u_\eps=u$ uniformly on relatively compact subsets of $\Om$.
\item[(ii)] If $u\in L^p_{loc}(\Om,E)$ for some $1\leq p<\infty$, then $\lim_{\eps \to 0^+}u_\eps=u$ in $L^p_{loc}$. 
\end{itemize}
\end{Thm}

For a proof of see for example \cite[Thm. 1, p. 123]{EG}. In particular,
\begin{align*}
\Phi_\eps^\ast u(x)=\int_{\R^n}\phi_{\eps}(y-x)u(y)dy. 
\end{align*}

Note however that $u_\eps$ is smooth only on $\Om_\eps \subsetneq \Om$. Thus $u_\eps\notin C^\infty(\Om,E)\cap L^p(\Om,E)$. The use of Friedrich mollifiers can therefore only directly show that $\A_{ls}=\A_w$ rather than $\A_s=\A_w$. On the other hand the use of Friedrich mollifiers do show that 
\begin{align*}
\lim_{\eps \to 0^+}(\A u_\eps ,w)=(u,\A^\ast w)
\end{align*}
for all $w\in C_0^1(\Om,F)$ and so if $\A$ has smooth coefficients $\A u_\eps \to v$ in the sense of distributions. 
Note that for $\phi$ as in \eqref{eq:Molli}
\begin{align*}
\Aa(x,\nabla \phi(x-y))=-\frac{2}{(1-\vert x\vert^2)^2}\Aa(x,x-y)\phi(x-y).
\end{align*}

\begin{Lem}[Friedrichs lemma]
\label{lem:F}
Let $\Phi_\eps$ be a Friedrichs mollifier according to Definition \ref{def:FMolli}. Let $\Om \subset \R^n$ be a domain and $\A$ be a first order operator of the form
\begin{align*}
\A=\sum_{j=1}A_j(x)\dv_j+B(x)
\end{align*}
such that $A_j\in C^1(\Om,\LL(E,F)$ for $j=1,2,...,n$ and $B\in C(\Om, \LL(E,F))$. Let $1\leq  p\leq \infty$ and $u\in \text{dom}_{p,\Om}(\A_w)$. Let $\A u=v$ and let $u_\eps=\Phi_\eps u$ and $v_\eps=\Phi_\eps v$. Then 
\begin{align*}
\A u_\eps-v_\eps=[\A,\Phi_\eps]u
\end{align*}
where
\begin{align*}
[\A,\Phi_\eps]u(x)=\int_{\R^n}\bigg(\sum_{j=1}^n\dv_{y_j}\big[(A_j(x)-A_j(y))\phi_\eps(x-y)\big]+(B(x)-B(y))\phi_\eps(x-y) \bigg)u(y)dy.
\end{align*}
\end{Lem}
We note in particular that the commutator $[\A,\Phi_\eps]$ is an \emph{integral operator} and that the partial derivatives $\dv_{y_j}$ act on both $A$ and $\phi$ but {\bf not} on $u$, so we have gained a derivative for $u$.

\begin{proof}
We follow the proof in \cite{KF}. If $u\in \text{dom}_{p,\Om}(\A_w)$ then by definition there exists a $v\in L^p(\Om,F)$ such that for all $w\in C^\infty_0(\Om,F)$
\begin{align*}
\int_{\Om}\langle v(x),w(x) \rangle dx=\int_{\Om}\langle u(x), \A^\ast w(x) \rangle dx,
\end{align*}
and where we have defined $\A u(x)=v(x)$. It is however very important to note that we can {\bf not} write
\begin{align*}
\A u(x)=\sum_{j=1}^nA_j(x)\dv_j u(x)+B(x)u(x)
\end{align*}
as the individual terms $A_j(x)\dv_j u(x)$ need not have any meaning. In particular, if $\A$ has smooth coefficients, then in general $\dv_j u(x)$ are distributions for $j=1,2,...,n,$ that cannot be identified with functions in $L^p(\Om,E)$. 
First using that $\nabla_x\phi_\eps(x-y)=-\nabla_y\phi_\eps(x-y)$ and that for a scalar function $\phi\in C^\infty(\Om)$ and a fixed vector $v\in E$ we have $\A \phi(x)v=\Aa(x,\nabla \phi(x))v$
we get 
\begin{align*}
(\A \circ \Phi_\eps) u(x)&=\int_{\R^n}\Aa(x,\nabla_x\phi_{\eps}(x-y))u(y)dy+\int_{\R^n}B(x)\phi_{\eps}(x-y)u(y)dy\\
&=-\int_{\R^n}\Aa(x,\nabla_y\phi_{\eps}(x-y))u(y)dy+\int_{\R^n}B(x)\phi_{\eps}(x-y)u(y)dy. 
\end{align*}

Furthermore,

\begin{align*}
\Aa(x,\nabla_y\phi_{\eps}(x-y))=\sum_{j=1}^n\dv_{y_j}(\Aa(x,e_j)\phi_{\eps}(x-y)).
\end{align*}

We now consider $(v_\eps,w)$. 
\begin{align*}
(v_\eps,w)=(\Phi_\eps v,w)=(v, \Phi_\eps^\ast w)=(u, \A^\ast \Phi_\eps^\ast w). 
\end{align*}

Noting that 
\begin{align*}
-\Aa(y,\nabla_y\phi_{\eps}(y-x))^\ast w(x)-\text{div}\,\Aa(y)^\ast\phi_{\eps}(y-x))w(x)=\sum_{j=1}^n\dv_{y_j}(\Aa(y,e_j)^\ast\phi(y-x)w(x))
\end{align*}
and using Fubini's theorem we get

\begin{align*}
&(u,\A^\ast \circ \Phi_\eps^\ast w)\\&=\int_{\R^n}\bigg\langle u(y),\int_{\R^n}\Aa^\ast (y,\nabla_x\phi_{\eps}(x-y))w(x)dx+\int_{\R^n}(B(y)^\ast -\text{div}\,\Aa(y)^\ast)\phi_{\eps}(x-y))w(x)dx\bigg\rangle dy\\
&=\int_{\R^n}\int_{\R^n}\langle u(y),\Aa^\ast (y,\nabla_x\phi_{\eps}(x-y))w(x)\rangle dx dy +\int_{\R^n}\int_{\R^n}\langle u(y),(B(y)^\ast -\text{div}\,\Aa(y)^\ast)\phi_{\eps}(x-y))w(x)\rangle dx dy\\
&=\int_{\R^n}\int_{\R^n}\langle u(y),-\Aa^\ast (y,\nabla_x\phi_{\eps}(x-y))w(x)\rangle dx dy +\int_{\R^n}\int_{\R^n}\langle u(y),(B(y)^\ast -\text{div}\,\Aa(y)^\ast)\phi_{\eps}(x-y))w(x)\rangle dx dy\\
&=\int_{\R^n}\int_{\R^n}\bigg\langle u(y),-\sum_{j=1}^n\dv_{y_j}(\Aa(y,e_j)^\ast\phi_\eps(y-x))w(x)+B(y)\phi_\eps(y-x)w(x)\bigg\rangle dx dy \\
&=\int_{\R^n}\int_{\R^n}\bigg\langle -\sum_{j=1}^n\dv_{y_j}(\Aa(y,e_j)\phi_\eps(y-x))u(y)+B(y)\phi_\eps(y-x)u(y),w(x)\bigg\rangle dx dy \\
&=\int_{\R^n}\bigg\langle \int_{\R^n}\bigg(-\sum_{j=1}^n\dv_{y_j}(\Aa(y,e_j)\phi_\eps(y-x))u(y)+B(y)\phi_\eps(y-x)u(y)\bigg)dy,w(x)\bigg\rangle dx 
\end{align*}

Thus,
\begin{align*}
(\A \circ \Phi_\eps) u(x)-v_\eps(x)&=\int_{\R^n}\bigg(\sum_{j=1}^n\dv_{y_j}\big[(\Aa(y,e_j)-\Aa(x,e_j))\phi_{\eps}(x-y)]+(B(x)-B(y))\phi_\eps(x-y)\bigg)u(y)dy.
\end{align*}

\end{proof}

Friedrichs' lemma has been extended also to pseduodifferential operators, see  \cite[Prop. 7.2, p.53]{Taylor1}.

\begin{Def}\label{def:FKernel}
For each mollifier $\Phi_\eps$ define the integral operator  
\begin{align*}
\mathcal{K}_\eps^{\Phi}u(x):=\int_{\R^n}K_\eps^\Phi(x,y)u(y)dy,
\end{align*}
where 
\begin{align}\label{eq:FKernel}
K_\eps^\Phi(x,y):=\sum_{j=1}^n\dv_{y_j}\big[(A_j(x)-A_j(y))\phi_\eps(x-y)\big]+(B(x)-B(y)\phi_\eps(x-y)
\end{align}
\end{Def}

Observe that the total symbol $\Aa\in C(\Om, \LL(\R^n,\LL(E,F)))$ equals 
\begin{align*}
\Aa(x)=\sum_{j=1}^nA_j(x)\otimes e_j.
\end{align*}
Therefore, using the divergence operator acting row-wise on the linear map $\LL(\R^n,\LL(E,F))$ we get 
\begin{align*}
K_\eps^\Phi(x,y)=\text{div}_y\,\big[ (\Aa(x)-\Aa(y))\phi_\eps(x-y)\big]+(B(x)-B(y))\phi_\eps(x-y).
\end{align*}

\begin{Def}
Let $\mathcal{K}$ be an integral operator with kernel $K\in C(\R^n \times \R^n, \LL(E,F))$ acting on functions $u\in C(\Om,E)$ through 
\begin{align*}
\mathcal{K}u(x)=\int_{\Om}K(x,y)u(y)dy. 
\end{align*}
Let 
\begin{align}
\Vert \mathcal{K}\Vert_\infty&=\esup_{y\in \R^n}\bigg\{\int_{\R^n}\Vert K(x,y)\Vert dx\bigg\}\\
\Vert \mathcal{K} \Vert_1&=\esup_{x\in \R^n}\bigg\{\int_{\R^n} \Vert K(x,y)\Vert dy\bigg\}
\end{align}
\end{Def}

\begin{Lem}[Schur estimates]
\label{lem:Schur}
Let $u\in L^p(\R^n,E)$. Then for $1\leq p\leq \infty$ and $1/p+1/q=1$
\begin{align*}
\Vert \mathcal{K}u\Vert_{L^p}\leq \Vert \mathcal{K}\Vert_1^p\Vert \mathcal{K}\Vert_\infty^q\Vert u\Vert_{L^p} 
\end{align*}
\end{Lem}

This lemma is a standard result, however for completeness, we give its proof. 
\begin{proof}
Firstly, 
\begin{align*}
\vert \mathcal{K}u(x)\vert &\leq \int_{\R^n}\Vert K(x,y)\Vert \vert u(y)\vert dy\leq \Vert u\Vert_\infty \int_{\R^n}\Vert K(x,y)\Vert dy
\end{align*}
which shows that $\Vert \mathcal{K}u\Vert_\infty \leq \Vert \mathcal{K}\Vert_\infty \Vert u\Vert_\infty$. Secondly, using Fubini's theorem 
\begin{align*}
\int_{\R^n}\vert \mathcal{K}u(x)\vert dx &\leq \int_{\R^n}\bigg(\int_{\R^n}\Vert K(x,y)\Vert dx\bigg) \vert u(y)\vert dy\leq \Vert \mathcal{K}\Vert_1 \Vert u\Vert_1.
\end{align*}
Finally, for $1<p<\infty$ by Hölder's inequality
\begin{align*}
\vert \mathcal{K}u(x)\vert^p &\leq \bigg(\int_{\R^n} \Vert K(x,y)\Vert^{1/q} \Vert K(x,y)\Vert^{1/p} \vert u(x)\vert dy\bigg)^p\\
&\leq \bigg(\int_{\R^n}\Vert K(x,y)\Vert dy \bigg)^{p/q}\int_{\R^n}\Vert K(x,y)\Vert \vert u(y)\vert^p dy,
\end{align*}
and so 
\begin{align*}
\int_{\R^n}\vert \mathcal{K}u(x)\vert^p dx&\leq \int_{\R^n}\bigg(\int_{\R^n}\Vert K(x,y)\Vert dy \bigg)^{p/q}\int_{\R^n}\Vert K(x,y)\Vert \vert u(y)\vert^p dydx\\
&\leq \Vert \mathcal{K}\Vert_\infty^{p/q}\Vert \mathcal{K}\Vert_1\Vert u\Vert_{p}^p.
\end{align*}
\end{proof}

Note that in Friedrichs' original work in \cite{KF} he uses 
\begin{align}
\Vert \mathcal{K}\Vert:=\max\{\Vert \mathcal{K}\Vert_1,\Vert \mathcal{K}\Vert_\infty\},
\end{align}
in which case one gets 
\begin{align*}
\Vert \mathcal{K}u\Vert_{L^p}\leq \Vert \mathcal{K}\Vert\Vert u\Vert_{L^p}
\end{align*}
for all $1\leq p\leq \infty$. 

\begin{Def}\label{def:K}
Let $\{\mathcal{K}_\eps\}$ be a one-parameter family of integral operators as in Definition \ref{def:FKernel} with kernels $K_\eps$. The family $\{\mathcal{K}_\eps\}$ satisfies the Friedrich properties $\text{(I)-(III}_0)$ if:\\
\begin{itemize}
\item The family $\{\mathcal{K}\}_\eps$ satisfy property (I) if $K_\eps(x,y)=0$ for all $\vert x-y\vert>\eps$. \\
\item The family $\{\mathcal{K}\}_\eps$ satisfy property (II) if there exist a $k>0$ such that $\Vert \mathcal{K}_\eps\Vert \leq k$ for all $\eps>0$. \\
\item Assume that $E=F$. The family $\{\mathcal{K}\}_\eps$ satisfy property (III) if there exists a number $\kappa\in \R$ such that 
\begin{align*}
\int_{\R^n}K_\eps(x,y)dy=\kappa I
\end{align*}
for all $x\in \R^n$, where $I$ is then identity transformation on $\LL(E)$.
\item Assume that $E\neq F$. The family $\{\mathcal{K}\}_\eps$ satisfy property $\text{III}_0$
\begin{align*}
\int_{\R^n}K_\eps(x,y)dy=0
\end{align*}
for all $x\in \R^n$.
\end{itemize}
\end{Def}

\begin{Lem}
\label{lem:F2}
Let $\{\mathcal{K}_\eps\}$ be a one-parameter family of integral operators as in Definition \ref{def:FKernel} with integral kernels $K_\eps$ satisfying $\text{(I),(II) and (III}_0)$ in Definition \ref{def:K}. Then for every $1\leq p\leq   \infty$,
\begin{align*}
\lim_{\eps\to 0^+}\Vert \mathcal{K}_\eps u\Vert_p=0. 
\end{align*}
\end{Lem}

\begin{proof}
First assume that $u\in C(\overline{U})$. By properties (I) and (III$)_0$ 
\begin{align*}
\bigg\vert \int_{\R^n}K_{\eps}(x,y)u(y)dy\bigg\vert&=\bigg\vert \int_{\R^n}K_{\eps}(x,y)u(y)dy-\int_{\R^n}K_{\eps}(x,y)u(x)dy\bigg\vert\\
& \bigg\vert \int_{B_\eps(x)}K_{\eps}(x,y)(u(y)-u(x))dy\bigg\vert.
\end{align*}
Thus, Lemma \ref{lem:Schur} and property (II) give
\begin{align*}
\Vert \mathcal{K}_\eps u\Vert_p^p\leq \Vert \mathcal{K}_\eps\Vert \Vert u-u(x)\Vert_{L^p(B_\eps(x))}\leq k\Vert u-u(x)\Vert_{L^p(B_\eps(x))}
\end{align*}

Assume that $u\in L^p(U, E)$ and $1\leq p<\infty$. By \cite[Corollary 1, p. 16]{EG} for any $\delta>0$ there exists an $u_\delta \in C(\overline{U},E)$ such that $\Vert u-u_\delta \Vert_p\leq \delta$. Thus
\begin{align*}
\lim_{\eps \to 0^+}\Vert \mathcal{K}_\eps u\Vert_p=\lim_{\eps \to 0^+}\Vert \mathcal{K}_\eps (u-u_\delta)\Vert_p+\Vert \mathcal{K}_\eps u_\delta\Vert_p\leq k\delta.
\end{align*}
Since $\delta>0$ was arbitrary the result follows. Now assume that $p=\infty$. Furthermore, we may assume $\vert U\vert<+\infty$. Then again by \cite[Corollary 1, p. 16]{EG} for any $\delta>0$ there exists a compact set $C$, such that $\vert A\setminus C\vert<\delta$ and a continuous function $u_\delta$ such that $u_\delta(x)=u(x)$ for $x\in C$. Thus 
\begin{align*}
\bigg\vert \int_{B_\eps(x)}K_{\eps}(x,y)u(y)dy\bigg\vert&\leq \bigg\vert \int_{B_\eps(x)\cap C}K_{\eps}(x,y)(u(y)-u(x))dy\bigg\vert+\bigg\vert \int_{B_\eps(x)\cap C}K_{\eps}(x,y)(u(y)-u(x))dy\bigg\vert\\
&\leq k\sup_{y\in B_\eps(x)}\{\vert u_\delta(x)-u_\delta(y)\vert\}+2k\Vert u\Vert_\infty \vert U\setminus C\vert.
\end{align*}
Hence
\begin{align*}
\lim_{\eps \to 0^+}\Vert \mathcal{K}_\eps u\Vert_\infty\leq 2\delta k\Vert u\Vert_\infty. 
\end{align*}
Since $\delta>0$ was arbitrary the result follows. 
\end{proof}

\begin{Lem}
\label{lem:F3}
Let $\mathcal{K}_\eps^{\Phi}$ be the Friedrichs commutator operator with $\Phi$ given by \eqref{eq:Molli}.  Then $\mathcal{K}_\eps^{\Phi}$ satisfies properties (I)-$\text{(III)}_0$ of Definition \ref{def:K}. 
\end{Lem}

\begin{proof}
Property $(I)$ is clear by the choice of $\phi$ and the formula \eqref{eq:FKernel}. We now investigate property (II). Using that 
\begin{align*}
\sum_{j=1}^n\dv_{y_j}\big[(A_j(x)-A_j(y))\phi_\eps(x-y)\big]=\Aa(x,\nabla \phi_\eps(x-y))-\Aa(y,\nabla \phi_\eps(x-y))-\text{div}\,\Aa(y)\phi_\eps(x-y), 
\end{align*}
we find that
\begin{align*}
\Vert K_{\eps}^{\Phi}(x,y)\Vert \leq \Vert \Aa(x,\nabla \phi_\eps(x-y))-\Aa(y,\nabla \phi_\eps(x-y))\Vert +\Vert \text{div}\,\Aa(y)\Vert \phi_\eps(x-y)
\end{align*}
Since $\Aa\in C^1(\Om\times \R^n,\LL(E,F))$ there exists a constant $C$ such that 
\begin{align*}
\Vert \Aa(x,\nabla \phi_\eps(x-y))-\Aa(y,\nabla \phi_\eps(x-y))\Vert \leq   C\vert x-y\vert \Vert \vert \nabla \phi_\eps(x-y)\vert
\end{align*}
Furthermore,
\begin{align*}
&\int_{\R^n}\vert x-y\vert\vert \nabla \phi_\eps(x-y)\vert dx=\int_{\R^n}\vert x\vert \vert \nabla \phi_\eps(x)\vert dx=\int_{\R^n}\frac{\vert x\vert}{\eps^{n+1}}\vert \nabla \phi(x/\eps)\vert dx\\
&\leq 2c\int_{B_\eps(0)}\frac{\eps}{\eps^{n+1}}\frac{\vert x/\eps\vert }{1-\vert x/\eps\vert^2}\phi(x/\eps)dx=2c\int_{B_1(0)}\frac{\vert x \vert }{1-\vert x\vert^2}\phi(x)dx\leq c'
\end{align*}
for some positive $c'>0$. In addition, since 
\begin{align*}
\Vert \text{div}\,\Aa(y)\Vert \phi_\eps(x-y)\leq c''\phi_\eps(x-y)
\end{align*}
we find that 
\begin{align*}
\max\bigg \{\int_{\R^n}\Vert K_{\eps}^{\Phi}(x,y)\Vert dx,\int_{\R^n}\Vert K_{\eps}^{\Phi}(x,y)\Vert dx\bigg\}\leq k
\end{align*}
independent of $\eps$. We now check property (III$)_0$. By Stokes' theorem 
\begin{align*}
\int_{\R^n}K_{\eps}^{\Phi}(x,y)dy&=\sum_{j=1}^n\int_{B_\eps(x)}\dv_{y_j}\big[(A_j(x)-A_j(y))\phi_\eps(x-y)\big]dy\\
&=\sum_{j=1}^n\int_{\dv B_\eps(x)}\nu(y_j)\big[(A_j(x)-A_j(y))\phi_\eps(x-y)\big]d\sigma(y)=0
\end{align*}
since $\phi_{\eps}\vert_{\dv B_\eps(x)}=0$.
\end{proof}

\begin{Thm}[Friedrichs equivalence theorem]
\label{thm:FriEquiv}
Assume that $A_j\in C^1(\Om,\LL(E,F))$, $B\in C(\Om,\LL(E,F))$. Then $\A_{ls}=\A_w$. 
\end{Thm}

\begin{proof}
The proof follows by combining Lemmata \ref{lem:F} -\ref{lem:F3}. 
\end{proof}

\begin{rem}
Theorem \ref{thm:FriEquiv} holds true if the assumptions on the coefficient matrices are relaxed to $A_j\in \text{Lip}(\Om,\LL(E,F))$ and $B\in L^\infty(\Om,\LL(E,F))$ since Lemma \ref{lem:F2} applies to this case as well. 
\end{rem}

It is interesting to note that Theorem \ref{thm:FriEquiv} does not hold for second order scalar operators with variable coefficients as was shown by Hörmander in \cite[Thm. 4.1]{H2}. There he proved the following result.
\begin{Thm}
Consider the variable coefficient hyperbolic operator 
\begin{align*}
\A u(x)=x_2\dv_1^2u(x)-\dv_1\dv_2 u(x). 
\end{align*}
Let $\psi \in C^\infty_0(\R^2)$ and consider the Friedrich mollifier $u_\eps =\Psi_\eps u=\psi_\eps\ast u$. Let $U,V$ be open sets in $\R^n$ such that $U\Subset V$. If for every $u\in L^2(V)$, satisfying the equation $\A u=0$ in the sense of distributions the norms $\Vert \A(\psi_\eps\ast u)\Vert_{L^2(U)}$ remains bounded as $\eps\to 0$, then $\psi=0$.  
\end{Thm}

What is surprising with the example is of course that we can rewrite the equation $\A u=0$ as a first order system using the factorisation 
\begin{align*}
\A u(x)=x_2\dv_1^2u(x)-\dv_1\dv_2 u(x)=x_2\dv_1(\dv_1-\dv_2)u(x)
\end{align*}
so that the equation becomes equivalent to the system 
\begin{align*}
\begin{bmatrix}
\dv_1-\dv_2 & 0\\
0 & x_2\dv_1
\end{bmatrix}
\begin{bmatrix}
u(x)\\
v(x)
\end{bmatrix}
+
\begin{bmatrix}
 0& -1\\
0 & 0
\end{bmatrix}
\begin{bmatrix}
u(x)\\
v(x)
\end{bmatrix}
=\begin{bmatrix}
0\\
0
\end{bmatrix}.
\end{align*}
However what Friedrich's equivalence theorem tells us is that 
\begin{align*}
\Vert ((\dv_1-\dv_2)u_\eps-v_\eps)^2+x_2^2(\dv_1 v_\eps)^2\Vert_{L^2(U,\R^2)}
\end{align*}
converges to $0$ as $\eps \to 0$, which is not the same as the expression $\A(\psi_\eps\ast u)$.

\begin{Def}[Operator of local type]
A PDE operator $\A$ is said to be of \emph{local type} on $\Om$ if for any $\zeta\in C^\infty_0(\Om)$ and any $u\in \text{dom}_{p,\Om}(\A_w)$, $\zeta u\in \text{dom}_{p,\Om}(\A_0)$. 
\end{Def}

One may also ask if $u\in \text{dom}_{p,\Om}(\A_w)$, is it true that for any $\varphi\in C^\infty(\overline{\Om})$, $\varphi u\in \text{dom}_{p,\Om}(\A_w)$? This property is very important for various types of localization arguments. It turns out that for general PDE operators this is false as the following two examples of Hörmander show.

\begin{ex}[Laplace operator]
Let $\A=\Delta$ be the Laplace operator in the plane and let $\Om=\Di$. By an example of Hadamard there exists a harmonic function $u\in C(\overline{\Di})$ such that 
\begin{align*}
\int_{\Di}\vert \nabla u(x)\vert^2 dx=+\infty, 
\end{align*} 
in particular we can chose $u$ so that $\dv_ru\notin L^2(\Di)$. Such a $u \in \text{dom}_{2,\Di}(\Delta)$. Now chose a $\varphi\in C^\infty(\overline{\Di})$ so that $\varphi(x)=\vert x\vert^2$ outside the origin. Then 
\begin{align*}
\Delta(\varphi u)(x)=u(x)\Delta \varphi(x)+2\langle \nabla u(x),\nabla \varphi(x)\rangle.
\end{align*}
The first term belongs to $L^2$ but the second term, being equal to $2\dv_r u$ outside the origin does not. Thus $\varphi u\notin \text{dom}_{2,\Di}(\Delta)$. 
\end{ex}

\begin{ex}[Wave operator]
Let $\A=\dv_1\dv_2$, the wave operator in dimension 2, and let $u=u(x_1)$ be an absolutely continuous function such that $\dv_1u\notin L^2$ in the neighbourhood of any point. Let $\Om=\{(x_1,x_2)\in \R^2: \vert x_1\vert<1, \,\,\, \vert x_2\vert<1\}$. Since for any $\varphi\in C^\infty(\overline{\Om})$
\begin{align*}
\dv_1\dv_2(\varphi u)(x)=\dv_1u\dv_2\varphi+u\dv_1\dv_2\varphi
\end{align*}
it follows that $\dv_1\dv_2(\varphi u)\notin L^2(\Om)$ unless $\varphi$ is a function of $x_1$ only. Thus $\varphi u\notin \text{dom}_{2,\Om}(\dv_1\dv_2)$ in general.
\end{ex}

However, for first order operators the situation changes for the better.

\begin{Lem}
Let $\A: C^1(\overline{\Om},E)\to C^0(\overline{\Om},F)$ be a first order PDE operator of the form \eqref{eq:OpA} Then $C^\infty(\overline{\Om}) \text{dom}_{p,\Om}(\A_w)\subset  \text{dom}_{p,\Om}(\A_w)$. 
\end{Lem}

\begin{proof}
For any $\varphi\in C^\infty(\overline{\Om})$ and $w\in C^{\infty}_0(\Om,F)$, using that $\varphi w\in C^{\infty}_0(\Om,F)$, we find
\begin{align*}
(\varphi u,\A^\ast w)&=(u, \varphi \A^\ast w)=(u, \A^\ast(\varphi w) )-\Aa(\nabla \varphi)^\ast w)=(\A u, \varphi w)-(\Aa(\nabla \varphi) u,w)\\
&=( \varphi \A u-\Aa(\nabla \varphi)u,w).
\end{align*}
Since $\varphi \A u-\Aa(\nabla \varphi)u\in L^p(\Om,F)$, it follows that $\varphi u\in \text{dom}_{p,\Om}(\A_w)$.  
\end{proof}

\begin{Thm}\label{thm:sW}
Let $\Om\subset \R^n$ be a domain. Then for any $\A$ of the form \eqref{eq:OpA} $\A_s=\A_{w}$. 
\end{Thm}

\begin{proof}
Here we follow the idea of the proof of \cite[Thm. 3.12]{H3} modified according to the proof of \cite[Thm. 2, p. 125]{EG} for Sobolev spaces. Fix $\eps>0$, set $U_0=\varnothing$ and 
\begin{align*}
U_k:=\bigg\{x\in \Om: \text{dist}(\dv \Om,U_k)>\frac{1}{k}\bigg\}, \quad k=1,2,3,...
\end{align*}
Set 
\begin{align*}
V_k:=U_{k+1}-\overline{U_{k-1}}
\end{align*}
and let $\{\zeta_k\}_k$ be a smooth partition of unity subordinate the covering $\cup_k V_k$ of $\Om$ so that $\zeta_k\in C^\infty_0(V_k)$, $k=1,2,3,..,$ $0\leq \zeta_k(x)\leq 1$ and $\sum_{k=1}^\infty\zeta_k(x)\equiv 1$. 
For $u\in L^p(\Om,E)$ and $v\in L^p(\Om,F)$ let $u_k=\zeta_k u$, $v_k=\zeta_k v$. Then $\text{supp}(u_k)\subset V_k$ and $\text{supp}(v_k)\subset V_k$. Let $\Phi_\eps$ be a Friedrich mollifier.  
Then for each $k=1,2,3,...$ there exists an $\eps_k>0$ such that $\text{supp}(\Phi_{\eps_k}u_{k})\subset V_k$ and $\text{supp}(\Phi_{\eps_k}v_{k})\subset V_k$ and 
 \begin{align*}
 \Vert \Phi_{\eps_k} (\zeta_k u)-\Phi_{\eps_k} u\Vert_{L^p(V_k,E)}<\frac{\eps}{2^k}, \quad \Vert \Phi_{\eps_k} (\zeta_k v)-\Phi_{\eps_k} v\Vert_{L^p(V_k,F)}<\frac{\eps}{2^k}.
 \end{align*}
Define 
\begin{align*}
u_\eps=\sum_{k=1}^\infty \Phi_{\eps_k} (\zeta_k u), \quad v_\eps=\sum_{k=1}^\infty \Phi_{\eps_k} (\zeta_k v).
\end{align*}
Since the sums above are finite for each $x\in \Om$, it follows that $u_\eps,v_\eps\in C^\infty$. Furthermore, by construction of $u_\eps$
\begin{align*}
\Vert u_\eps-u\Vert_{L^p(\Om,E)}\leq \sum_{k=1}^\infty \Vert u_\eps-\zeta_ku\Vert_{L^p(\Om,E)}\leq \eps,
\end{align*} 
and similarly for $v_\eps$, it follows that $\lim_{\eps\to 0^+}u_\eps=u$ and $\lim_{\eps\to 0^+}v_\eps=v$ in $L^p(\Om)$. As in Friedrich's lemma we compute $\A u_\eps -v_\eps$. 
For any $\phi\in C^\infty_0(\Om,F)$
\begin{align*}
(v_\eps,\phi)&=\sum_{k=1}^\infty(\Phi_{\eps_k} (\zeta_k v),\phi)=\sum_{k=1}^\infty( \zeta_k v,\Phi_{\eps_k}^\ast \phi)=(v,\sum_{k=1}^\infty\zeta_k\Phi_{\eps_k}^\ast \phi)\\
&=(u, \A ^\ast \sum_{k=1}^\infty\zeta_k\Phi_{\eps_k}^\ast \phi)=\sum_{k=1}^\infty(u, \Aa^\ast  (\nabla \zeta_k)^\ast \Phi_{\eps_k}^\ast \phi))+(u, \zeta_k\A^\ast \Phi_{\eps_k}^\ast \phi))\\
&=\sum_{k=1}^\infty(\Phi_{\eps_k}\Aa (\nabla \zeta_k)u,  \phi)+((\A^\ast \Phi_{\eps_k}^\ast )^\ast \zeta_k u, \phi)).
\end{align*}
Thus, 
\begin{align*}
\A u_\eps -v_\eps&=\A \sum_{k=1}^\infty \Phi_{\eps_k} (\zeta_k u)-v_\eps=\sum_{k=1}^\infty (\A \Phi_{\eps_k}-(\A^\ast \Phi_{\eps_k}^\ast )^\ast )\zeta_k u))-\sum_{k=1}^\infty\Phi_{\eps_k}\Aa (\nabla \zeta_k)u\\
&=\sum_{k=1}^\infty \mathcal{K}^{\Phi}_{\eps_k}\zeta_k u-\sum_{k=1}^\infty\Phi_{\eps_k}\Aa (\nabla \zeta_k)u. 
\end{align*}

Thus,
\begin{align*}
\A u_\eps -v_\eps=\sum_{k=1}^{\infty}\mathcal{K}_{\eps_k}\zeta_k u-\sum_{k=1}^{\infty}\Phi_{\eps_k}\Aa (\nabla \zeta_k)u
\end{align*}
Since
\begin{align*}
\lim_{\eps \to 0^+}\sum_{k=1}^{\infty}\Phi_{\eps_k}\Aa (\nabla \zeta_k)u&=\sum_{k=1}^{\infty}\Aa (\nabla \zeta_k)u=\Aa\bigg(\sum_{k=1}^{\infty} \nabla \zeta_k\bigg)u\\
&=\Aa\bigg(\nabla \sum_{k=1}^{\infty} \zeta_k\bigg)u=\Aa(\nabla 1)u=\Aa(0)u=0,
\end{align*}
and by Theorem \ref{thm:FriEquiv}, $\mathcal{K}_{\eps_k}\zeta_k u\to 0$ in $L^p$, $\A u_\eps -v_\eps \to 0$ as $\eps \to 0^+$ in $L^p$ and the proof is complete. 
\end{proof}

\begin{Def}[Contractible domains]
A bounded domain $\Om \subset \R^n$ is said to be \emph{contractible} in the sense of Friedrichs if for every $\eps>0$ there exists a homeomorphism $f_\eps: \Om \to \Om_\eps$ where $\Om_\eps \Subset \Om$ such that $\text{dist}(\Om_\eps,\dv \Om)\geq \eps$. 
\end{Def}

\begin{ex}
The following is an example of a non-contractible domain. Let $B_1(0)$ be the open unit ball in $\R^2$ and let $\Om =B_1(0)\setminus (\cup_{k=1}^\infty\{k^{-1}e_1\})$. 
\end{ex}

\begin{Def}[Smoothly contractible domains]
A bounded domain $\Om \subset \R^n$ is said to be \emph{smoothly contractible} if for every $\eps>0$ there exists a homeomorphism $f_\eps: \Om \to \Om_\eps$ where $\Om_\eps \Subset \Om$ such that $\text{dist}(\Om_\eps,\dv \Om)\geq \eps$ and in addition the following conditions hold on $f_\eps$:
\begin{itemize}
\item[]
\item[(i)] $f_\eps\in C^2(\Om,\Om_\eps)$.
\item[]
\item[(ii)] $\displaystyle \sup_{\eps>0}\sup_{x\in \Om}\frac{\vert f_\eps(x)-x\vert}{\eps}\leq \infty$.
\item[]
\item[(iii)] $\displaystyle \sup_{\eps>0}\sup_{x\in \Om}\frac{\vert Df_\eps(x)-I\vert}{\eps}\leq \infty$.
\item[]
\item[(iii)] $\displaystyle \sup_{\eps>0}\sup_{x\in \Om}\vert D^2f_\eps(x)\vert\leq \infty$ and $\displaystyle \lim_{\eps\to 0^+}\Vert D^2f_\eps\Vert_{L^\infty(U)}=0$ for all $U\Subset \Om$. 
\end{itemize}
\end{Def}

An example of a contractible domain which is not smoothly contractible is $\Om=B_1(0)\setminus \{0\}\subset \R^2$.

\begin{Thm}[Friedrichs equivalence theorem II]
\label{thm:FrThm2}
Assume that $A_j\in C^1(\overline{\Om})$, $B\in C(\overline{\Om})$. Assume that $\Om\subset \R^n$ is a smoothly contractible domain. Then $\A_{S}=\A_w$. 
\end{Thm}

The proof of \ref{thm:FrThm2} is based on considering mollification with $\phi_\eps(x-f_\eps(x))$ and using the definition of smoothly contractible domains. The proof otherwise proceeds along the same lines as the proof of Theorem \ref{thm:FriEquiv}. Instead of reproving Theorem \ref{thm:FrThm2} we will sketch a proof of a more general theorem for strong Lipschitz domains mimicking the proof for Sobolev functions given in \cite[Thm. 3, p. 127]{EG}. 

\begin{Def}\label{def:LipDom}
A domain $\Om \subset \R^n$ is called a strong Lipschitz domain if for each point $x\in \dv \Om$ there exists a $r>0$ and a Lipschitz map $\rho: \R^{n-1}\to \R$ such that after rotating and relabelling the coordinate axes 
\begin{align*}
\dv \Om \cap Q(x,r)=\{y: \rho(y_1,...,y_{n-1})<y_n\}\cap Q(x,r),
\end{align*}
where $Q(x,r)=\{y: \vert y_j-x_j\vert<r, \,\,\, j=1,2,...,n\}$. 
\end{Def}

\begin{Thm}\label{thm:BoundaryApprox}
Let $\Om\subset \R^n$ be a bounded strong Lipschitz domain and assume that $A_j\in \text{Lip}(\overline{\Om})$, $B\in L^\infty(\overline{\Om})$. Then for every $u\in \text{dom}_{p,\Om}(\A)$ there exists a sequence $\{u_j\}_j\subset C^\infty(\overline{\Om},E)$ such that 
\begin{align*}
\lim_{j\to \infty}\Vert u-u_j\Vert_p+\Vert \A u-\A u_j\Vert_p=0.
\end{align*}
\end{Thm}

\begin{proof}
Since $\A$ satisfies the principle of locality, if $u\in \text{dom}_{p,\Om}(\A)$, then $\varphi u\in \text{dom}_{p,\Om}(\A)$ for any $\varphi\in C^\infty(\overline{\Om})$. For $x_0\in \dv \Om$, take $r>0$ and $\rho:\R^{n-1}\to \R$ as in Definition \ref{def:LipDom}. Let $Q'=Q(x_0,r/2)$. After multiplying by some $\varphi\in C^\infty(\overline{\Om})$ we may assume that $u$ vanishes near $\dv Q'\cap \Om$. For $x\in \Om \cap Q'$, $\eps>0$, $\alpha>0$, define 
\begin{align*}
x^\eps=x+\eps \alpha e_n
\end{align*}
and note that $B(y^\eps,\eps)\subset \Om\cap Q(x_0,r)$ for all $\eps>0$ sufficiently small provided $\alpha$ is large enough, $\alpha=\text{Lip}(\rho)+2$ suffice. Define a mollifying operator through 
\begin{align*}
\Phi_\eps u(x)=\int_{\R^n}\phi_\eps(y)u(x^\eps-y)dy= \int_{B_{\eps}(x^\eps)}\phi_\eps(x-y+\eps\alpha e_n) u(y)dy
\end{align*} 
where $\phi$ is as in \eqref{eq:Molli}. Exactly as in the proof  \cite[Thm. 3, p. 127]{EG} one verifies that $\Phi_\eps u\in C^\infty(\overline{\Om \cap Q'})$. After this step one computes the Friedrich commutator with respect to $\Phi_\eps$ and $\A$. The rest of the proof follows along the same lines as the proof of Theorem \ref{thm:FriEquiv} and the details are left to the reader. 
\end{proof}

A proof of Theorem \ref{thm:BoundaryApprox} under the assumption that the domain is $C^1$ can also be found in \cite[Prop. 1]{R1}.

So far we have said nothing about the closure of differential operators with respect to boundary values. We will not go into this topic as it deserves a separate treatise. The interested reader is instead referred to \cite{Sara,Sara62,Sara73,LP60}.

Both the concept of weak and strong extensions have many virtues, however what is not immediately clear from their definitions is what type of more general functions that are allowed in $\text{dom}(\mathscr{A}_w)$ nor what type of regularity that is enforced by weak (and strong) solutions of $\mathscr{A}_wu=0$. It is for these matters we believe that a different concept of extension due to B. Fuglede  called \emph{flux extension} may sometimes be more helpful. In many ways it is a more geometric counterpart of the notion weak extensions and it is the topic of the next section.

%============NEW SUBSECTION=================================================================

\section{\sffamily Flux extensions}\label{sec:Flux}

\begin{Def}[Baire function]
A function $f$ on $\Om$ is called a \emph{Baire function} if $f$ is the pointwise limit of a sequence of continuous functions on $\Om$, i.e., there exists a sequence $\{f_j\}_j\subset C(\Om)$ such that 
\begin{align*}
f(x)=\lim_{j\to \infty}f_j(x). 
\end{align*}
\end{Def}
Baire functions have stronger continuity properties than general Lebsegue integrable functions. In particular the points of continuity of a Baire function are a comeagre $G_\delta$ set, see \cite[Thm. 24.14, p. 193]{Kech95}

\begin{Def}[Baire sets and Borel sets]
Let $X$ be a topological space. The Baire $\sigma$-algebra $\mathscr{B}^\ast(X)$ on $X$ is the smallest $\sigma$-algebra so that all bounded continuous functions on $X$ are measurable. The Borel $\sigma$-algebra $\mathscr{B}(X)$ on $X$ is the smallest $\sigma$-algebra of $X$ that contains all open sets of $X$.
\end{Def}

 In general $\mathscr{B}^\ast(X)\subset \mathscr{B}(X)$, however if $X$ is separable metric space then $\mathscr{B}^\ast(X)= \mathscr{B}(X)$. Furthermore, Baire functions are precisely the function that are Baire measurable. For more about the theory of Borel and Baire sets we refer to \cite{Bic,Kech95}.

If $u\in C^1(\Om,E)$, then for any smooth domain $U\Subset \Om$ Stokes theorem implies that 
\begin{align*}
\int_{U}\mathscr{A}u(x)dx=\int_{\dv U}\mathbb{A}(y,\nu(y))u(y)d\sigma(y)+\int_{U}(B(x)-\text{div}\, \Aa(x))u(x)dx
\end{align*}
If $u$ merely belongs to $L^p(\Om,E)$ then by Fubini's theorem the boundary integral 
\begin{align*}
\int_{\dv U}\mathbb{A}(y,\nu(y))u(y)d\sigma(y)
\end{align*}
is well-defined only for ``almost all'' hypersurfaces $\dv U$. Of course this statement does not make sense a priori since the family of all smooth domains compactly contained in $\Om$ is infinite dimensional. However should $U$ be a family of balls centered at $x\in \Om$ say, then the statement holds for all most every $B_\eps(x)$ in the usual sense. It was the great insight of B. Fuglede in \cite{F1} to give a precise meaning to the notion of all most all domains which we will now recall.

\begin{Def}[$k$-dimensional Lipschitz surfaces]
Let $1\leq k\leq n-1$. A non-empty subset $S\subset \R^n$ is called a \emph{$k$-dimensional Lipschitz surface} if there exists to each point $x\in S$ an open set $U\subset\R^n$ such that $x\in U$ and $S\cap U$ is the image of a Lipschitz map $f: V\subset \R^k\to \R^n$. 
\end{Def}

\begin{Def}[Moduli of family of $k$-dimensional Lipschitz surfaces]
Let $\mathscr{F}^k$ be a family of $k$-dimensional Lipschitz surfaces in $\R^n$. Let $\mathcal{B}(\mathscr{F}^k)$ denote the the set of all non-negative Baire functions $\rho: \R^n \to \R$ be a such that for every $\Sigma\in \mathscr{F}^k$
\begin{align*}
\int_{\Sigma}\rho(x)d\sigma(x)\geq 1, 
\end{align*}
where $\sigma=\mathscr{H}^k \lfloor \Sigma$, and $\mathscr{H}^k$ is the $k$-dimensional Hausdorff measure. The $p$-\emph{moduli} $\mathbf{\text{\sffamily M}}_p(\mathscr{F}^k)$ of the family $\mathscr{F}^k$ for $1\leq p<+\infty$ is defined according to 
\begin{align}\label{eq:Moduli}
\mathbf{\text{\sffamily M}}_p^k(\mathscr{F}^k):=\inf_{\rho\in \mathcal{B}(\mathscr{F}^k)}\int_{\R^n}\rho(x)^pdx.
\end{align}
A $k$-dimensional family of Lipschitz surface $\mathscr{F}'\in \mathscr{F}^k$ is said to be \emph{$p$-exceptional} if $\mathbf{\text{\sffamily M}}_p(\mathscr{F}')=0$. 
\end{Def}

\begin{rem}
Note that if $\mu$ is a Radon measure and $f$ a nonnegative Borel measurable function then 
\begin{align*}
\int f(x)d\mu(x)=&\inf\bigg\{\int g(x) d\mu(x): g\geq f \text{ and }g\text{ is lower semicontinuous} \bigg\}. 
\end{align*}
(see \cite[Prop. 7.14, p. 219]{Folland}). It therefore suffice to consider only lower-semicontinuous functions $\rho$ in \eqref{eq:Moduli}. 
\end{rem}

Let $\Om\subset \R^n$ be a domain. Denote by $\mathcal{M}(\Om)$ the space of Radon measures on $\Om$ and by $\mathcal{M}_0(\Om)$ the space of finite Radon measures. We may embed any Green domain $\Sigma\subset \Om$ into $\mathcal{M}(\Om)$ by identifying $\Sigma$ with the Radon measure $\sigma=\mathcal{H}^{n-1}\lfloor \Sigma$. It is therefore natural to extend the concept of moduli also to the space of Radon measures. This is indeed done in \cite{F1}, to which we refer the reader there for more details.

\begin{Thm}[Fuglede]
The $p$-moduli $\mathbf{\text{\sffamily M}}_p$ is monotone and countably subadditive. That is:
\begin{itemize}
\item[(i)]  $\mathbf{\text{\sffamily M}}_p(\mathscr{F}^k)\leq \mathbf{\text{\sffamily M}}_p(\mathscr{F'}^k)$ if $\mathscr{F}^k\subset \mathscr{F'}^k$
\item[(ii)] $\mathbf{\text{\sffamily M}}_p(\mathscr{F}^k)\leq \sum_i\mathbf{\text{\sffamily M}}_p(\mathscr{F}_i^k)$ if $\mathscr{F}^k=\cup_i\mathscr{F}^k_i$. 
\end{itemize} 
\end{Thm}

For a proof we refer the reader to the original paper \cite{F1} by Fuglede.

\begin{Def}[$p$-exceptional family of Green domains]
\label{Def:pExc}
Let $1\leq p<+\infty$. A family $\mathscr{E}$ of Green domains in $\R^n$ is called \emph{exceptional of order $p$}, abbreviated $p$-exc, if there exists a non-negative Borel function $\rho\in L^p(\R^n)$ such that \begin{align*}
\int_{\Sigma}\rho(x)d\sigma(x)=+\infty
\end{align*} 
for every $\Sigma\in \mathscr{E}$. If a proposition regarding a specified family of Green sets form a family $\mathscr{F}$ holds for all but a subfamily $\mathscr{F}'$ which is $p$-exc we say that it holds for $p$-a.e. Green set. 
\end{Def}

Note that by \cite[Thm. 2]{F1} Definition \ref{Def:pExc} is equivalent to $\mathbf{\text{\sffamily M}}_p(\mathscr{E})=0$.

It is instructive to see how the notion of $p$-exceptional family of Green domains coincides with the usual measure theoretic notion derived from Fubini's theorem in the case when then family $\mathscr{E}$ is finite dimensional. 

\begin{ex}\label{ex:NullFam}
Let $p=1$. Consider the finite dimensional family $\mathscr{F}$ of all spheres $\dv B_r(x)\subset \R^n$ for all $x\in \R^n$ and all $r>0$. Thus we can identify $\mathscr{F}\cong \R^n\times \R_+$ and equip it with the Lebesgue measure on $\R^{n+1}$ restricted to $\R^n\times \R_+$. Consider the subset $\mathscr{F}'$ consisting of all balls with center in $\mathbb{Q}^n$ and such that the radii are all rational numbers. Let $\{r_j\}_j$ be an enumeration of the positive rational numbers and let $\{p_j\}$ likewise be an enumeration of $\mathbb{Q}^n$. Let $1<\alpha<n$ and define
\begin{align*}
\rho_0(x)=\sum_{j=1}^\infty\frac{1}{2^j}\frac{1}{\vert \vert x\vert-r_j\vert^{\alpha}}.
\end{align*}
Then $\rho_0(r_j)=+\infty$ and 
\begin{align*}
\int_{\R^n}\rho_0(x)dx\leq \sum_{j=1}^\infty\frac{1}{2^j}\int_{\R^n}\frac{1}{\vert \vert x\vert-r_j\vert^{\alpha}}dx<+\infty
\end{align*}
Now let 
\begin{align*}
\rho(x)=\sum_{j=0}^\infty\frac{1}{2^j}\rho_0(x-p_j).
\end{align*}
Then clearly,
\begin{align*}
\int_{\dv B_{r_j}(p_j)}\rho(x)dx>\frac{1}{2^j}\int_{\dv B_{r_j}(p_j)}\rho_0(x-p_j)dx=+\infty. 
\end{align*}
Thus the family  $\mathscr{F}$ is $1$-exc. 
\end{ex}

Since families of balls have inconvenient covering properties, \cite{F2} consider the families of cubes  (where they are called intervals). An open cube $Q\subset \R^n$ is defined according to 
\begin{align*}
Q=\{x\in \R^n:\, a_j<x_j<b_j: \,\, j=1,2,...,n.\}
\end{align*}
We will consider the family of open sets $\mathscr{F}_Q$ such that each $U\in \mathscr{F}_Q$ is given by 
\begin{align*}
U=\bigg(\bigcup_{j=1}^N \overline{Q_j}\bigg)
\end{align*}
for some $N$ where $Q_i\cap Q_j=\varnothing$ if $i\leq j$. Furthermore open cubes generate the standard euclidean topology of $\R^n$ and also the Borel $\sigma$-algebra of $\R^n$. As in Example \ref{ex:NullFam}, it is shown in \cite[p. 25]{F2} that the concept of $p$-a.e. cube coincides with the measure theoretic notion of a.e. cube by identifying $\mathscr{F}_Q\cong \R^n \times \R$ and endowing $\mathscr{F}_Q$ with the Lebesgue measure on $\R^{n+1}$.

We have the following result from \cite[Thm. 2]{F2} about unions and subfamilies of $p$-exceptional $k$-dimensional Lipschitz surfaces.

\begin{Thm}
\label{thm:Exc}
Let $\Om\subset \R^n$ be a domain.
\begin{itemize}
\item[(a)] Any subfamily of a $p$-exc family is $p$-exc.
\item[(b)] The union of a countable family of $p$-exc family is $p$-exc.
\item[(c)] If $p>q$, then every $p$-exc system is likewise $q$-exc. 
\item[(d)] If $A \subset \Om$ and $\vert A\vert=0$, then 
\begin{align*}
\int_{A\cap \dv U}d\sigma=0
\end{align*}
for $p$-a.e. Green domain $U$. 
\item[(f)] If a sequence $\{u_j\}_j\subset L^p_{loc}(\Om)$ converges in $L^p_{loc}$ to some  $u\in L^p_{loc}(\Om)$, then there exists a subsequence $\{u_{j_k}\}_k$ which, for $p$-a.e. Green domain $U$ converges to $u$ in $L^1(\dv U,\sigma)$, i.e. 
\begin{align*}
\lim_{k\to \infty}\int_{\dv U}\vert u_{j_k}-u\vert d\sigma =0
\end{align*}
for $p$-a.e. Green domain $U$. 
\end{itemize}
\end{Thm}

Furthermore, Lipschitz maps preserves $p$-exceptional families of Lipschitz surfaces. More precisely we have
\begin{Thm}[Theorem 4 in \cite{F1}]
Let $U\subset \R^n$ be open and let $f: U\to \R^n$ be a Lipschitz map such that $V=f(U)$ is open. Then any family of $p$-exceptional $k$-dimensional Lipschitz surfaces contained in $U$ is mapped to a $p$-exceptional system $k$-dimensional Lipschitz surfaces contained in $V$. 
\end{Thm}

There is a close connection between capacity and moduli, where we recall that the $p$-capacity of a compact set $K\subset \R^n$ is given by 
\begin{align*}
\text{cap}_p(K):=\inf_{\substack{u\in W^{1,p}(\R^n)\\ u(x)=1, \,\,\, x\in K} }\bigg\{\int_{\R^n \setminus K}\vert \nabla u(x)\vert^pdx\bigg\}
\end{align*}

\begin{Thm}[Ziemer]
Let $K\subset \R^n$ be a compact set with Hausdorff dimension 1. Then 
\begin{align}
\text{cap}_p(K)=\mathbf{\text{\sffamily M}}_p^1(K). 
\end{align}
\end{Thm}

For a proof see \cite{Z2,Z1}.

Using the notion of moduli one can give a beautiful characterisation of Sobolev spaces due to J. Heinonen. We first recall the following definition. 
\begin{Def}[ACL on lines]
A function $u:\R^n \to \R$ is said to be \emph{absolutely continuous on lines} if $u$ is absolutely continuous on almost every line segment parallel to the coordinate axes. Furthermore, a function $u\in L^p(\R^n)$ is said to be of class $\text{AC}L_p(\R^n)$ if it is absolutely continuous on lines and and the distributional gradient satisfy $\nabla u\in L^p(\R^n,\R^n)$. A vector valued function $u: \R^n \to E$ is said to satisfy the same definition if each coordinate function is in $\text{ACL}_p(\R^n).$
\end{Def}

A classical equivalent characterisation of the Sobolev space $W^{1,p}(\R^n)$ is that for each $[u]\in W^{1,p}(\R^n)$ there exists a representative $u\in \text{ACL}_p(\R^n)$, where $[u]$ is to emphasise that $u\in W^{1,p}(\R^n)$ is really an equivalence class of functions. This characterization leads to the natural question of how many curves (not just lines) that a function $u\in W^{1,p}(\R^n)$, or rather a representative of its equivalence class, can be absolutely continuous on. Since we are now typically interested in an infinite dimensional family of curves, a statement along the lines of almost every curve would not make any sense in the measure theoretic sense. Here the notion of $p$-moduli comes to the rescue. We have (\cite[Thm. 7.4, Thm 7.6]{Hei1}):

\begin{Thm}[\cite{Hei1}]
Let $\Gamma$ be the family of all rectifiable curves in $\R^n$. Then every function in $W^{1,p}(\R^n)$ has a representative that is absolutely continuous on $p$-almost every curve in $\Gamma$ in the sense of moduli. Furthermore, a function $u\in L^p(\R^n)$ belongs to $W^{1,p}(\R^n)$ if and only if there exists a Borel function $\rho\in L^p(\R^n)$ such that the inequality 
\begin{align*}
\vert u(\gamma(a))-u(\gamma(b))\vert \leq \int_\gamma \rho(x)d\sigma(x)
\end{align*}
for $p$-a.e. rectifiable curve $\gamma:[0,1]\to \R^n$. 
\end{Thm}

The connection between moduli and capacity in fact extends also to compact $(n-1)$-dimensional Lipschitz surfaces in the following way. 
\begin{Thm}[Theorem 9 in \cite{F1}]
Let $K\subset \R^n$ be an arbitrary compact set, $n\geq 3$ and let $U^{ex}_K$ be the unbounded component of $\R^n-K$. Let $\mathscr{H}_{Lip}$ be the family of all compact $(n-1)$-dimensional Lipschitz surfaces such that every $\Sigma \in \mathscr{H}_{Lip}$, $\Sigma \subset U^{ex}_K$ separates $K$ from infinity. Then 
\begin{align*}
\mathbf{\text{\sffamily M}}_2^{n-1}(\mathscr{H}_{Lip})=\frac{1}{(n-2)\omega_{n}\text{cap}_2(K)}. 
\end{align*}
\end{Thm}

We are now ready to give the definition of flux extensions. 

\begin{Def}[Flux extensions]
\label{def:FluxE}
Let $1\leq p<+\infty$. Consider the set of all pairs $\text{gr}(\mathscr{A}_f)=(u,v)\in L^p(\Om,E)\times L^p(\Om,F)$ such that the Stokes' or Gauss-Green formula 
\begin{align}\label{eq:FluxE}
\int_{U}v(x)dx=\int_{\dv U}\mathbb{A}(y,\nu(y))u(y)d\sigma(y)+\int_{U}(B(x)-\di \Aa(x))u(x)dx
\end{align}
holds for the operator $\mathscr{A}$ for $p$-a.e. every Green set $U\Subset \Om$. The graph $\text{gr}(\mathscr{A}_f)$ now defines an extension $\mathscr{A}_f$ of $\mathscr{A}$ called the \emph{flux extension} and we define $\A_f u:=v$. 
\end{Def}
Note that the left hand side of \eqref{eq:FluxE} is well defined for $p$-a.e. Green set in view of Theorem \ref{thm:Exc} (b) and (e). In addition, flux extensions only depends on the equivalence classes of $u$ and $v$ in $L^p(\Om,E)$ and $L^p(\Om,F)$ respectively due to Theorem \ref{thm:Exc} (d).

We now come to the key point of using $p$-exceptional systems and insist that the relation \eqref{eq:FluxE} only holds for $p$-a.e. Green domain, namely that this will imply that the graph $\text{gr}(\mathscr{A}_f)$ in Definition \ref{def:FluxE} is closed.

\begin{Thm}[Fuglede's completeness theorem]
\label{thm:CompleteFuglede}
The graph $\text{gr}(\mathscr{A}_f)$ is closed in $L^p(\Om,E)\times L^p(\Om,F)$.
\end{Thm}

\begin{proof}
If $(0,v)\in \text{gr}(\mathscr{A}_f)$, then the relation $\int_Uv(x)dx=0$ holds for $p$-a.e. Green domain $U\Subset \Om$. By definition this implies that $\int_Uv(x)dx=0$ for every Borel set $U\Subset \Om$ and hence $v(x)=0$ for a.e. $x\in \Om$. Next let $\{(u_j,v_j)\}_j\subset \text{gr}(\mathscr{A}_f)$ such that $\lim_{j}u_j=u$ and  $\lim_{j}v_j=v$ in $L^p$. We want to show that $(u,v)\in \text{gr}(\mathscr{A}_f)$. According to Theorem \ref{thm:Exc} (b) and (f) there exists a subsequence $\{(u_{j_k},v_{j_k})\}_k\subset \text{gr}(\mathscr{A}_f)$ and a $p$-exceptional system $\mathscr{E}$ of Green domain in $\Om$ such that 
\begin{align*}
\lim_{k\to \infty}u_{j_k}=u
\end{align*}
in $L^1(\dv U,d\sigma; E)$ for every Green domain $U\Subset \Om$ such that $U\notin \mathscr{E}$. For every index $i$ there exists a $p$-exceptional system $\mathscr{E}_i$ such that 
\begin{align*}
\int_{U}\mathscr{A}u_i(x)dx=\int_{\dv U}\mathbb{A}(y,\nu(y))u_i(y)d\sigma(y)+\int_{U}(B(x)-\di \Aa(x))u_i(x)dx
\end{align*}
for every Green domain $U$ such that $u\notin \mathscr{E}_i$. We conclude that \eqref{eq:FluxE} holds for the subsequence $\{u_{j_k},v_{j_k}\}_k$ and every $U\notin \bigcup_{j_k}\mathscr{E}_{j_k}$. By Theorem \ref{thm:Exc} $ \bigcup_{j_k}\mathscr{E}_{j_k}$ is a again a $p$-exceptional system and so we conclude that 
\begin{align*}
\int_U v(x)dx&=\lim_{k\to \infty}\int_U v_{j_k}(x)dx=\lim_{k\to \infty}\int_{\dv U}\mathbb{A}(y,\nu(y))u_{i_k}(y)d\sigma(y)+\int_{U}(B(x)-\di \Aa(x))u_{i_k}(x)dx\\
&=\int_{\dv U}\mathbb{A}(y,\nu(y))u(y)d\sigma(y)+\int_{U}(B(x)-\di \Aa(x))u(x)dx
\end{align*}
for $p$-a.e. Green domain $U$. 
\end{proof}

\begin{Thm}[Fuglede's equivalence theorem]
\label{thm:EquivalenceFuglede}
\begin{align*}
\mathscr{A}_f=\mathscr{A}_{ls}
\end{align*}
\end{Thm}

Together with Friedrichs equivalence theorem, Theorem \ref{thm:FriEquiv}, we have the equality $\mathscr{A}_w=\mathscr{A}_f$. 

\begin{proof}
We want to show that $\mathscr{A}_fu=v$ implies $\mathscr{A}_{ls}u=v$. By the principle of locality it is enough to show this for all $\Om'\Subset \Om$. Chose a mollifier $\rho_\eps$ such that $\text{supp}(\rho_\eps)+\Om'\Subset \Om$ for all $\eps>0$ and $\text{diam}(\text{supp}(\rho_\eps))<\eps$. Let 
\begin{align*}
u_\eps(x)=\rho_\eps \ast u(x), \quad v_\eps(x)=\rho_\eps\ast v(x). 
\end{align*}
Then for $p$-a.e. Green domain $U\Subset \Om'$ and a.e. $\xi \in \R^n$ such that $\vert \xi\vert<\eps$ we have by assumption that $u\in \text{dom}\,\mathscr{A}_f$

\begin{align*}
\int_{U-\xi}v(y)dy&=\int_{\dv U-\xi}\mathbb{A}(y,\nu(y))u(y)d\sigma(y)+\int_{ U-\xi}V(y)u(y)dy
\end{align*}
where we defined $V(x):=B(x)-\di \Aa(x)$, or 
\begin{align*}
\int_{U}v(y-\xi)dy&=\int_{\dv U}\mathbb{A}(y-\xi,\nu(y))u(y)d\sigma(y)+\int_{U}V(y-\xi)u(y-\xi)dy\\
\end{align*}
where we note that we get $\nu(y)$ and not $\nu(y-\xi)$ in the integral. By Fubini's theorem 
\begin{align*}
\int_{U}v_\eps(x)dx&=\int_U\int_{\R^n}\rho_\eps(x-\xi)v(\xi)d\xi dx=\int_{\R^n}\rho_\eps(\xi)\bigg(\int_Uv(x-\xi)dx\bigg)d\xi\\
&=\int_{\R^n}\rho_\eps(\xi)\bigg(\int_{\dv U }\mathbb{A}(x-\xi,\nu(x))u(x-\xi)d\sigma(x)+\int_{U }V(x-\xi)u(x-\xi)dx\bigg)d\xi\\
&=\int_{\dv U }\bigg(\int_{\R^n}\mathbb{A}(x-\xi,\nu(x))u(x-\xi)\rho_\eps(\xi)d\xi \bigg)d\sigma(x)+\int_{U }\bigg(\int_{\R^n}V(x-\xi)u(x-\xi)\rho_\eps(\xi)d\xi\bigg)dx\\
&=\int_{\dv U }\bigg(\int_{\R^n}\mathbb{A}(y,\nu(x))u(y)\rho_\eps(x-y)dy\bigg)d\sigma(x)+\int_{U }\bigg(\int_{\R^n}V(y)u(y)\rho_\eps(x-y)dy\bigg)dx\\
\end{align*}
where we made the change of variables $y=x-\xi$. Thus,
\begin{align*}
\int_{U}v_\eps(x)dx&=\int_{\dv U }\bigg(\int_{\R^n}\mathbb{A}(y,\nu(x))u(y)\rho_\eps(x-y)dy\bigg)d\sigma(x)+\int_{U }\bigg(\int_{\R^n}V(y)u(y)\rho_\eps(x-y)dy\bigg)dx.
\end{align*}
We now observe that the function 
\begin{align*}
g_\eps(x):=\int_{\R^n}\mathbb{A}(y,\nu(x))u(y)\rho_\eps(x-y)dy=\sum_{j=1}^n\nu_j(x)\int_{\R^n}\mathbb{A}(y,e_j)u(y)\rho_\eps(x-y)dy
\end{align*}
is $C^1$ and hence we can apply Stokes theorem in the form 
\begin{align*}
\int_{U}\dv_{j}u(x)dx=\int_{\dv U}\nu_j(x)u(x)d\sigma(x)
\end{align*}
for each term in the right hand side and in addition differentiate under the integral sign giving 
\begin{align*}
&\int_{\dv U}\bigg(\int_{\R^n}\mathbb{A}(y,\nu(x))\rho_\eps(x-y)u(y)dy\bigg)d\sigma(x)=\sum_{j=1}^n\int_{\dv U}\nu_j(x)\bigg(\int_{\R^n}A_j(y,e_j)\rho_\eps(x-y)u(y)dy\bigg)d\sigma(x)\\
&=\sum_{j=1}^n\int_{ U}\dv_{x_j}\bigg(\int_{\R^n}A_j(y,e_j)\rho_\eps(x-y)u(y)dy\bigg)dx=\int_{ U}\bigg(\int_{\R^n}\sum_{j=1}^nA_j(y,e_j)\dv_{x_j}\rho_\eps(x-y)u(y)dy\bigg)dx\\
&=\int_{ U}\bigg(\int_{\R^n}\mathbb{A}(y,\nabla_x \rho_{\eps}(x-y))u(y)dy\bigg)dx\\
\end{align*}
for $p$-a.e. Green domain $U\Subset \Om'$ and since $v_\eps$ is $C^1$ on $\Om'$ we have the identity 
\begin{align*}
v_\eps(x)=\int_{\R^n}\mathbb{A}(y,\nu(x))u(y)\rho_\eps(x-y)dy+\int_{\R^n}V(y)u(y)\rho_\eps(x-y)dy
\end{align*}
for all $x\in \Om'$.

On the other hand 
\begin{align*}
\mathscr{A}_fu_\eps(x)&=\int_{\R^n}\mathscr{A}_x\rho_\eps(x-y)u(y)dy=\int_{\R^n}\Big[\mathbb{A}(x,\nabla_x \rho_\eps(x-y))+V(x)\rho_\eps(x-y)\Big]u(y)dy
\end{align*}

\begin{align*}
v_\eps(x)-\mathscr{A}_fu_\eps(x)&=\int_{\R^n}\bigg[\mathbb{A}(y,\nabla_x \rho_{\eps}(x-y))-\mathbb{A}(x,\nabla_x \rho_\eps(x-y))+(V(y)-V(x))\rho_\eps(x-y)\bigg]u(y)dy\\
&=\int_{\R^n}\mathcal{K}^\rho_\eps(x,y)u(y)dy
\end{align*}
where $\mathcal{K}^\rho_\eps$ is the Friedrichs mollifier kernel with respect to $\rho$ as in Lemma \ref{lem:F3}. The rest of the proof is identical to the proof of Theorem \ref{thm:FriEquiv}. 
\end{proof}

\begin{rem}
A proof of Theorem \ref{thm:EquivalenceFuglede} in the special case of the Cauchy-Riemann operator in the plane following the same ideas is given in \cite{Zalc}, apparently unaware of the more general case proven by Fuglede. 
\end{rem}

We will now give a different direct proof that $\text{dom}_{p,\Om}(\A_w)\subset \text{dom}_{p,\Om}(\A_f)$ that avoids going through the locally strong extension $\A_{ls}$ and Friedrichs lemma. 

\begin{Prop}
\begin{align*}
\text{gr}_{p,\Om}(\A_w)\subset \text{gr}_{p,\Om}(\A_f)
\end{align*}
\end{Prop}

\begin{proof}
Let $u\in \text{dom}_{p,\Om}(\A_w)$. Then for $p$-a.e. Green domain $U\Subset \Om$ ($C^1$), $u$ has an $L^p$-trace on $\dv U$. Let $\phi_\eps$ be a Friedrichs mollifier and let $\chi_\eps=\phi_\eps\ast \chi_U$ and chose a fixed vector $f$ and test function $\varphi_\eps=\chi_\eps f\in C^\infty_0(\Om,F)$. Since $u\in \text{dom}_{p,\Om}(\A_w)$, by definition there exists a $\A_wu\in L^p(\Om,F)$ such that 
\begin{align*}
\int_{\Om}\langle \A_w u(x),\varphi_\eps(x)\rangle dx=\int_{\Om}\langle u(x),\A^\ast \varphi_\eps(x)\rangle dx
\end{align*}
Furthermore, 
\begin{align*}
\lim_{\eps \to 0^+}\int_{\Om}\langle \A_w u(x),\varphi_\eps(x)\rangle dx=\int_{U}\langle \A_w u(x),f\rangle dx.
\end{align*}
In addition 
\begin{align*}
\A^\ast \varphi_\eps(x)=\Aa(\nabla \chi_\eps(x))f,
\end{align*}
and $\Aa(\nabla \chi_\eps)u\to \Aa(\nu)u\delta_{\dv U}$ in $\mathscr{E}'(\Om, \LL(E,F))$ since $U$ is a set of finite perimeter.  Thus 
\begin{align*}
\int_{\Om}\langle \A_w u(x),f\rangle dx=(\Aa(\nu)u\delta_{\dv U},f)=\int_{\dv_U}\langle \Aa(\nu(x))u(x),f\rangle d\sigma(x). 
\end{align*}
By letting $\{f_j\}_j$ be an ON-basis for $F$ and expanding $\A_w u$ in this basis the result follows. 
\end{proof}

The concept of flux extension also leads to the notion of \emph{flux solution} of $\mathscr{A}_fu=0$ which of course in view of Fuglede's theorem is equivalent to the notion of weak solution. 
\begin{Def}[Flux solutions]
$u\in \text{dom}(\mathscr{A}_f)\subset L^p(\Om,E)$ is a flux solution of the equation 
\begin{align*}
\A u(x)=f(x)
\end{align*}
in $\Om$ with $f\in L^p(\Om,F)$ if 
\begin{align*}
\int_{\dv U}\mathbb{A}(y,\nu(y))u(y)d\sigma(y)+\int_{U}(B(x)-\text{div}\, \Aa(x))u(x)dx=\int_Uf(x)dx
\end{align*}
for $p$-a.e. every Green domain $U\Subset \Om$. 
\end{Def}

%===================SUBSECTION========================================================
\subsection{\sffamily Characterization of the domain of flux extensions}

In this section we will give a characterization of the domain of flux extensions due to Fuglede. To that end we will recall an extension of a theorem of F. Riesz (\cite{Riesz}) also due to to Fuglede in \cite{F3}. 
We will first need some preliminary notions from measure theory.

\begin{Def}
Let $X\subset \R^n$ and let $\mathcal{A}\subset \mathcal{P}(X)$ be any family of subsets that contains the empty set. $\mathcal{A}$ is called a \emph{set algebra} or a \emph{field of sets} if 
\begin{enumerate}
\item $\varnothing \in \mathcal{A}$.
\item If $A\in \mathcal{A}$ then $X\setminus A\in \mathcal{A}$.
\item $\mathcal{A}$ is closed under finite unions and intersections.
\end{enumerate}
If in addition $\mathcal{A}$ is closed under countable unions and intersections we say that $\A$ is a \emph{$\sigma$-algebra}.
\end{Def}

\begin{Def}
Let $X\subset \R^n$ and let $\mathcal{A}\subset \mathcal{P}(X)$ be any family of subsets that contains the empty set. $\mathcal{A}$ is called a \emph{semiring} if: 
\begin{enumerate}
\item $\varnothing \in \mathcal{A}$.
\item For all $A,B\in \mathscr{A}$ then $A\cap B\in \mathcal{A}$:
\item For all $A,B\in \mathscr{A}$ there exists pairwise disjoint sets $I_k\in \mathscr{A}$ $J_k\in \mathscr{A}$ such that 
\begin{align*}
A\setminus B=\bigcup_{k=1}^NI_j, \quad A\cup B=\bigcup_{k=1}^MJ_k.
\end{align*}
\end{enumerate}
In particular any algebra of sets is a semiring. 
\end{Def}

For any family $\mathcal{A}\subset \mathcal{P}(X)$ of subsets that contains the empty set we denote by $\mathcal{A}_{\sigma}$ all countable unions from $\mathcal{A}$ and by $\mathcal{A}_{\delta}$ all countable intersections from $\mathcal{A}_{\sigma}$. We further set $\mathcal{A}_{\delta \sigma}=(\mathcal{A}_{\delta})_\sigma$.

\begin{Def}
Let $X\subset \R^n$ and let $\mathcal{A}\subset \mathcal{P}(X)$ be any family of subsets that contains the empty set. Let $\alpha: \mathcal{A}\to \overline{\R}_+=\{t\in \R: t\geq 0\}\cup \{\infty\}$ be a set function (where we always require that $\varnothing\in \mathcal{A}$). 
\begin{enumerate}
\item $\alpha$ is called monotone if for all $A,B\in \mathcal{A}$ with $A\subset B$ implies $\alpha(A)\leq \alpha(B)$. 
\item  $\alpha$ is called additive if for any finite family $\{A_k\}_k^N$ of pairwise disjoint sets in $\mathcal{A}$ with $\bigcup_{k=1}^NA_k\in \mathcal{A}$
\begin{align*}
\alpha\bigg(\bigcup_k^NA_k\bigg)=\sum_{k=1}^N\alpha(A_k).
\end{align*}
\item $\alpha$ is called countably additive if for any countable family $\{A_k\}_k$ of pairwise disjoint sets in $\mathcal{A}$ with $\bigcup_{k=1}^\infty A_k\in \mathcal{A}$
\begin{align*}
\alpha\bigg(\bigcup_k^\infty A_k\bigg)=\sum_{k=1}^\infty \alpha(A_k)
\end{align*}
\item $\alpha$ is called countably subadditive if for any countable family $\{A_k\}_k$ of sets in $\mathcal{A}$ with $\bigcup_{k=1}^\infty A_k\in \mathcal{A}$
\begin{align*}
\alpha\bigg(\bigcup_k^\infty A_k\bigg)\leq \sum_{k=1}^\infty \alpha(A_k)
\end{align*}
\end{enumerate}
\end{Def}

Any set function $\alpha$ defined on any family $\mathcal{A}\subset \mathcal{P}(X)$ of subsets can be extended to an outer measure on $\mathcal{P}(X)$ in the following way. 
\begin{align}\label{eq:ExtensionM}
\mu_\alpha^\ast(A):=\inf\bigg\{\sum_{j=1}^\infty \alpha(A_j): \bigcup_{j=1}^\infty A_j\supset A, \,\,\, A_j\in \mathcal{A}\bigg\}
\end{align}
with the understanding that $\mu_\alpha^\ast(A)=+\infty$ if no sequence $\{A_j\}\subset \mathcal{A}$ exists such that $\bigcup_{j=1}^\infty A_j\supset A$. In general however, $\mu_\alpha^\ast$ need not be a an extension of $\alpha$ and sets in $\mathcal{A}$ need not be $\mu_\alpha^\ast$ measurable. To avoid these pathologies further properties of $\mathcal{A}$ and $\alpha$ are needed. 

\begin{Thm}\label{thm:ExtensionM}
Let $X\subset \R^n$ and let $\mathcal{A}\subset \mathcal{P}(X)$ be a semiring. Let $\alpha: \mathcal{A}\to \overline{\R}_+$ be a countably additive set function. Then $\alpha$ is monotone and countably subadditive on $\mathcal{A}$. Furthermore, let $\mu_\alpha^\ast$ be the outer measure extension of $\alpha$ defined by \eqref{eq:ExtensionM},  $\Sigma_{\mu^\ast}$ the corresponding class of $\mu_\alpha^\ast$-measurable sets and $(\mu_{\alpha}^\ast,\Sigma_{\mu^\ast})$ the measure associated to $\alpha$.  Then the following holds:
\begin{enumerate}
\item $\mu_\alpha^\ast$ extends $\alpha$, i.e, $\mu_\alpha^\ast(A)=\alpha(A)$ for all $A\in \mathcal{A}$.
\item A set $A\in \mathcal{P}(X)$ is $\mu_\alpha^\ast$-measurable if and only if 
\begin{align*}
\mu_\alpha^\ast(A\cap B)+\mu_\alpha^\ast(A\cap B^c)\leq \mu_\alpha^\ast(B)
\end{align*}
for all $B\in \mathcal{A}$ with $\mu^\ast_\alpha(B)<+\infty$.
\item $\mathcal{A}\subset \Sigma_{\mu^\ast}$.
\item For all $A\in X$ with $\mu_\alpha^\ast(A)<+\infty$, there exists a decreaseing sequence of sets $\{D_k\}_k$ each being $\mu^\ast_\alpha$-measurable with $\mu^\ast_\alpha(D_k)<+\infty$ and such that 
\begin{align*}
D_k=\bigcup_{l=1}^\infty I_l^{(k)}
\end{align*}
where $I_j^{(k)}\in \mathcal{A}$ for all $j,k$ and such that $A\subset \cap_{k=1}^\infty D_k$, and $\mu^\ast_\alpha(A)=\mu^\ast_\alpha(\cap_kD_k)$. 
\end{enumerate}
\end{Thm}

For a proof of Theorem \ref{thm:ExtensionM} see the nice exposition in \cite[Thm. 5.29, p. 299]{GM12}. Also note that should $\mathcal{A}$ generate the Borel $\sigma$-algebra of $X$, then $(\mu^\ast_\alpha,\Sigma_{\mu^\ast_\alpha})$ is a Borel measure. Theorem \ref{thm:ExtensionM} can also be extended to $E$-valued countably additive set functions $\alpha$. Namely, after fixing a basis for $E$ we can identify $\alpha$ with $(\alpha_1,...,\alpha_N)$, where $N=\text{dim}(E)$ and such that $\alpha_j$ is an $\R$-valued set function. Furthermore, we have 
\begin{align*}
\alpha_j=\alpha_j^+-\alpha_j^-, \quad j=1,2,...,N,
\end{align*}
where $\alpha_j^+$ and  $\alpha_j^-$ is the positive and negative part of $\alpha_j$ respectively. We can then apply Theorem  \ref{thm:ExtensionM} to each $\alpha_j^\pm$. This gives an $E$-valued outer measure $\mu_\alpha^\ast$. In principle one should also check that different choices of bases for $E$ gives rise to the same extension $\mu_\alpha^\ast$, this step however will be omitted. 

Before proceeding to the theorem of Fuglede we will first recall the original theorem of F. Riesz.

\begin{Thm}[F. Riesz]
\label{thm:Riesz}
Let $1<p<+\infty$. A necessary and sufficient condition for a function $F$ to be the indefinite integral of a $f\in L^p((a,b))$ is that 
\begin{align}\label{eq:PVariation}
\sup_{\mathcal{P}}\sum_{j}\frac{\vert F(x_j)-F(x_{j+1})\vert^p}{\vert x_j-x_{j+1}\vert^p}<+\infty
\end{align}
where $a\leq x_0<x_1<....<x_N=b$ and the supremum ranges over all finite partitions of $(a,b)$. Furthermore, the supremum is given by 
\begin{align*}
\int_{a}^{b}\vert f(x)\vert^pdx.
\end{align*}
If $p=1$, then a necessary and sufficient condition for a function $F$ to be the indefinite integral of a $f\in L^1((a,b))$ is that to every $\eps>0$ there exists a $\delta=\delta(\eps)$ such that 
\begin{align}\label{eq:1Variation}
\sum_{j=1}^{n-1}\vert F(x_j)-F(x_{j+1})\vert< \eps
\end{align}
for every finite system of mutually disjoint intervals $(x_1,x_2),....(x_{n-1},x_n)$ of $(a,b)$ for which 
\begin{align*}
\sum_{j=1}^{n-1}\vert x_j-x_{j+1}\vert <\delta. 
\end{align*}
\end{Thm}

We note in particular that in the case $p=1$ condition \eqref{eq:PVariation} is the condition that $F$ is absolutely continuous and 
\begin{align}\label{eq:1Variation}
V_{[a,b]}(f)=\sup_{\mathcal{P}}\sum_{j=1}\vert F(x_j)-F(x_{j+1})\vert
\end{align}
is the total variation of $F$. We can therefore think of \eqref{eq:PVariation} as the total $p$-variation of $F$. Furthermore, F. Riesz Theorem \ref{thm:Riesz} can be seen as a sharpening of Theorem \ref{thm:ExtensionM} in the sense that the set countably additive set function defined through
\begin{align*}
\alpha_f(\cup_jI_j)=\sum_{j=1}F(x_j)-F(x_{j-1})
\end{align*} 
has an extension $\mu_\alpha^\ast$ that is absolutely continuous with repspect to the Lebesgue measure $\mathcal{L}^1$ on $[a,b]$, and such that the  Radon-Nikodym derivative of $\mu^\ast_\alpha$ with respect to $\mathcal{L}^1$ if a function $f\in L^p([a,b])$.

\begin{Thm}[Fuglede-Riesz representation theorem]
\label{thm:RieszFuglede}
Assume that $(X,\Sigma, \mu)$ is a measure space. Let $\mathcal{U}$ be a collection of $\mu$-measurable subsets of a set $X$ and assume that all finite unions of disjoint sets form $\mathcal{U}$ together with $\varnothing$ form an algebra of sets $\mathscr{F}$ over $X$.  Assume furthermore that $\Sigma$ is generated by $\mathscr{F}$.  Let $E$ be a euclidean vector space and let $\varphi: \mathscr{U}\to E$ be an additive set function. \newline
Assume that $1<p<+\infty$.
Then there exists a function $f\in L^p(X,\Sigma,\mu)$ such that 
\begin{align}\label{eq:RepFR}
\varphi(A)=\int_Af(x)d\mu(x)
\end{align}
for every $A\in \mathscr{U}$ if and only if there exists a constant $c_p>0$ such that 
\begin{align}
\sum_{j}^n\frac{\vert \varphi(A_j)\vert^p}{\mu(A_j)^{p-1}}\leq c_p
\end{align}
for every finite collection of disjoint sets $A_1,A_2,...,A_n$ from $\mathscr{U}$. The function $f$ is unique in its equivalence class and the smallest possible value for $c_p$ is $\int_X\vert f(x)\vert^pd\mu(x)$.  \newline
Assume that $p=\infty$. Then \eqref{eq:RepFR} holds for some $f\in L^\infty(X,\Sigma,\mu)$ if and only if there exists a constant $c_\infty>0$ such that
\begin{align}
\vert \varphi(A)\vert\leq c_\infty\mu(A)
\end{align}
for every set $A\in \mathscr{U}$. The equivalence class of $f$ is unique and the smallest possible value for $c$ is $\text{ess sup}_{x\in X}\vert f(x)\vert$. \newline
Assume that $p=1$. 
Then \eqref{eq:RepFR} holds for some $f\in L^1(X,\Sigma,\mu)$ if and only if the following two conditions hold:
\begin{enumerate}
\item For every $\eps>0$ there exists a $\delta=\delta(\eps)>0$ with the property that 
\begin{align}
\sum_{j}^n\vert \varphi(A_j)\vert\leq \eps
\end{align}
for every finite collection of disjoint sets $A_1,A_2,...,A_n$ from $\mathscr{U}$ such that 
\begin{align}
\sum_{j}^n\vert \mu(A_j)\leq \delta.
\end{align}
\item There exists a finite constant $c_1>0$ such that 
\begin{align}
\sum_{j}^n\vert \varphi(A_j)\vert\leq c_1
\end{align}
for every finite collection of disjoint sets $A_1,A_2,...,A_n$ from $\mathscr{U}$.
\end{enumerate}
The equivalence class of $f$ is unique and the smallest possible value for $c_1$ is $\int_X\vert f(x)\vert d\mu(x)$. 
\end{Thm}

We now come to Fuglede's characterization of the domain $\text{dom}_{p,\Om}(\A)$ of a first order operator satisfying the standard assumptions. The idea is the following. Define an $E$-valued set function on the collection of $p$-almost all Green domains of $\Om$ through 
\begin{align}
\alpha_u(U):=\int_{\dv U}\mathbb{A}(y,\nu(y))u(y)d\sigma(y)+\int_{U}(B(x)-\di \Aa(x))u(x)dx
\end{align}
We let $\mathcal{U}$ be any subcollection of Green domains such that all finite unions of disjoint sets of $\mathcal{U}$ together with $\varnothing$ form an algebra of sets $\mathscr{F}$ over $X$ that generate the Borel $\sigma$-algebra $\mathscr{B}(X)$. We then apply Theorem \ref{thm:RieszFuglede}. This yields the following conditions.

\begin{Def}[Fuglede's condition $\mathbf{A_p}$ and condition $\mathbf{B}$]
Let $u\in L^p(\Om,E)$ is said to satisfy condition $\mathbf{A}_p$ if there exists a constant $C_p>0$ independent of $u$ such that for $p$-almost every choice of finite collection of disjoint Green sets $U_j\Subset \Om$
\begin{align}
\sum_j \frac{1}{\vert U_j\vert^{p-1}}\bigg\vert \int_{\dv U_j}\mathbb{A}(y,\nu(y))u(y)d\sigma(y)+\int_{U_j}(B(x)-\di \Aa(x))u(x)dx\bigg\vert^p\leq C_p. 
\end{align}
Furthermore $u$ is said to satisfy condition $\mathbf{B}$ if for every $\eps>0$ there exists a $\delta=\delta(\eps)>0$ such that 
\begin{align}
\sum_j\vert U_j\vert<\delta 
\end{align}
implies that 
\begin{align}
\sum_j\bigg\vert \int_{\dv U_j}\mathbb{A}(y,\nu(y))u(y)d\sigma(y)+\int_{U_j}(B(x)-\di \Aa(x))u(x)dx)\bigg\vert<\eps 
\end{align}
for $p$-almost every choice of finite collection of disjoint Green sets $U_j\Subset \Om$.
\end{Def}

This finally gives us the following characterization of the domain of flux extensions proved by Fuglede in \cite{F2}. 

\begin{Thm}[Characterisation of domains of flux extensions]
\label{thm:FugledeChar}
Let $u\in L^p(\Om,E)$. A necessary condition for $u$ to belong to $\text{dom}(\mathscr{A}^f)$ is that condition $\mathbf{A}_p$ holds. If $1<p<+\infty$ this condition also sufficient. If $p=1$, condition $\mathbf{B}$ is necessary and together with condition $\mathbf{A}_1$ sufficient. If $\vert \Om\vert<+\infty$ then condition $\mathbf{B}$ alone is sufficient. 
\end{Thm}

\begin{proof}
We first show the necessity of condition $\mathbf{A}_p$. Assume that $\A_f u=v$. Let $\mathscr{B}(\Om)$ denote the Borel $\sigma$-algebra of $\Om$. Define the function $\varphi: \mathscr{B}(\Om)\to F$ through
\begin{align*}
\varphi(U):=\int_U v(x)dx
\end{align*}
for all $U\in \mathscr{B}(\Om)$ with $\vert U\vert<+\infty$. By Jensen's inequality 
\begin{align*}
\bigg\vert \frac{1}{\vert U\vert}\int_U v(x)dx\bigg\vert^p\leq \frac{1}{\vert U\vert}\int_U \vert v(x)\vert^p dx
\end{align*}
or 
\begin{align*}
\frac{1}{\vert U\vert^{p-1}}\bigg\vert \int_U v(x)dx\bigg\vert^p\leq \int_U \vert v(x)\vert^p dx.
\end{align*}
Thus in particular for any finite disjoint family of Green domains 
\begin{align*}
\sum_j\frac{1}{\vert U_j\vert^{p-1}}\bigg\vert \int_{U_j} v(x)dx\bigg\vert^p\leq \int_{\cup_{j} U_{j}} \vert v(x)\vert^p dx\leq \Vert v\Vert_{L^p(\Om,F)}<+\infty. 
\end{align*}
We now show the sufficiency of condition $\mathbf{A}_p$. Define the set function $\varphi: \mathscr{F}\to F$ on $p$-a.e Green domain through 
\begin{align*}
\varphi(U):=\int_{\dv U}\Aa(x,\nu(x))u(x)d\sigma(x)+\int_U(B(x)-\di \Aa(x))u(x)dx. 
\end{align*}
$\varphi$ is finitely additive and fulfills $\mathbf{A}_p$ for $p$-a.e. every disjoint family of Green domains. The family of $p$-a.e. Green domains $\mathscr{F}$ generates the Borel $\sigma$-algebra $\mathscr{B}(\Om)$. Thus by Theorem \ref{thm:ExtensionM}, $\varphi$ admits an extension $\mu_\varphi: \mathscr{B}(\Om)\to F$. By the Fuglede-Riesz representation Theorem \ref{thm:RieszFuglede} there exists a function $v\in L^p(\Om,F)$ such that 
\begin{align*}
\mu_{\varphi}(U)=\int_U v(x)dx
\end{align*}
for every Borel set $E\in \mathscr{B}(\Om)$ with $\vert U\vert<+\infty$. 
\end{proof}

%==============NEW SUBSECTIONS=================================

\subsection{\sffamily When Stokes theorem holds for every Lipschitz Green set}\label{subsec:StokesAll}

It is natural to ask for which first order operators 
\begin{align}\label{eq:StokesGreen}
\int_U\A u(x)dx=\int_{\dv U}\Aa(y,\nu(y))u(y)d\sigma(y)+\int_{U}(B(x)-\text{div}\,\Aa(x))u(x)dx
\end{align}
holds \emph{for every} Green domain $U\Subset \Om$ with some additional regularity say Lipschitz? To this end we recall the definition of elliptic operator. 
\begin{Def}\label{def:Elliptic}
A first order partial differential operator operator $\A$ on a domain $\Om$ is called \emph{elliptic} if the principal symbol $\Aa$ satisfies
\begin{align}\label{def:Elliptic}
\Aa(x,\xi): E\to F \quad \text{is injective for all $x\in \Om$ and $\xi\in \R^n\setminus\{0\}$}.
\end{align}
If $\text{dim}(F)>\text{dim}(E)$ we say that $\A$ is \emph{overdetermined elliptic} and if $\text{dim}(F)=\text{dim}(E)$ we say that $\A$ is \emph{determined elliptic}. Finally, $\A$ is called \emph{$\C$-elliptic} if 
\begin{align}\label{def:CElliptic}
\Aa(x,\zeta): E\to F \quad \text{is injective for all $x\in \Om$ and $\zeta\in \C^n\setminus\{0\}$}.
\end{align}
\end{Def}
Of course these definitions also extend to higher order differential operators, which however is not the focus of this treatise. For constant coefficient homogeneous elliptic first order operators we have the following classical theorem of Caldéron and Zygmund. 

\begin{Thm}\cite{CZ}
Let $1<p<\infty$ and let $\A$ be a homogeneous constant coefficients partial differential operator. Then the estimate
\begin{align}
\Vert D^k u\Vert_{L^p}\leq C\Vert \A u\Vert_{L^p}
\end{align}
holds for every $u\in C^\infty_0(\Om,E)$ and some positive constant $C$ independent of $u$ if and only if $\A$ is elliptic. 
\end{Thm}
This theorem has later been extended to include variable coefficient operators, see for example \cite{SW}. A consequence of these theorems is that for any elliptic first order operator (and more generally of arbitrary order) we have for every $1<p<\infty$
\begin{align*}
\text{dom}_{p,\Om}(\A_0)=W^{1,p}_0(\Om,E), \quad \text{dom}_{p,\Om}(\A_w)\subset W^{1,p}_{loc}(\Om,E).
\end{align*}

Furthermore, if $\A$ is a $\C$-elliptic and $\Om$ is Lipschitz domain even stronger results hold.

\begin{Thm}[Aronszajn-Ne{\u c}as-Smith]
\label{thm:ANS}
Let $\Om\subset \R^n$ be a bounded Lipschitz domain and let $\A$ be first order operator with $C^\infty$-coefficients. Assume that the principal symbol $\Aa$ is elliptic on $\Om$ and $\C$-elliptic on $\dv \Om$. Then there exists a constant $c>0$ such that for every $u\in \text{dom}_{p,\Om}(\A_w)$
\begin{align*}
c\Vert u\Vert_{W^{1,p}(\Om,E)}\leq \Vert \A u\Vert_{L^p(\Om,F)}+\Vert u\Vert_{L^p(\Om,E)}
\end{align*}
and $\text{ker}(\A_w)$ on $\text{dom}_{p,\Om}(\A_w)$ is finite dimensional.
\end{Thm}
For a proof of Theorem \ref{thm:ANS} see \cite{Ar54,Nec,Sm61,Sm70}, in particular \cite[Thm. 8.15]{Sm70}. Examples of $\C$-elliptic first order operators besides the total derivative $D$ when $E=\R^n$ are: 
\begin{itemize}
\item[]
\item[] The symmetric gradient $D_Su(x)=\frac{1}{2}(Du(x)+Du(x)^\ast)$  where \newline
$\text{ker}(\mathcal{E})=\{Ax+b: A\in \LL(\R^n), \,\,\, A^\ast =-A,\,\,\, b\in \R^n\}$.
\item[]
\item[] The deviatoric derivative $\mathcal{E}u(x)=\frac{1}{2}(Du(x)+Du(x)^\ast)-\frac{1}{n}\text{tr}(Du(x))$ where \newline
$\text{ker}(\mathcal{E})=\{\langle a,x\rangle x-\frac{a}{2}\vert x\vert^2+Ax+b+cx: A\in \LL(\R^n), \,\,\, A^\ast =-A,\,\,\, a,b\in \R^n, \,\,\, c\in \R\}$.

\item[]
\item[] The Ahlfors operator $\mathcal{A}u(x)=\frac{1}{2}(Du(x)-Du(x)^\ast)+\frac{1}{n}\text{tr}(Du(x))$ \newline
$\text{ker}(\mathcal{A})=\{a+2\langle c,x\rangle x-c\vert x\vert^2+Bx: B\in \LL(\R^n), \,\,\, B^\ast =-B,\,\,\, a,c\in \R^n\}$.
\item[]
\end{itemize}
For a reference on the Ahlfors operator see \cite{Ahl1}.

A consequence of Theorem \ref{thm:ANS} is that for any Lipschitz domain $\Om$ and any $\C$-elliptic operator $\A$ we have for any $1<p<\infty$,
\begin{align*}
\text{dom}_{p,\Om}(\A)=W^{1,p}(\Om,E). 
\end{align*}

The above considerations show that for any first order elliptic operator $\A$ and any $u\in \text{dom}_{p,\Om}(\A)$ we have for every Green domain $U\Subset \Om$
\begin{align*}
\text{Tr}_{\dv U}=u\vert_{\dv U}\in W^{1-1/p,p}(\dv U,\sigma). 
\end{align*}
Thus $u$ has a well-defined trace in the classical sense on the boundary of any Lipschitz domain and \eqref{eq:StokesGreen} holds for every Lipschitz domain $U$ and every $u\in \text{dom}_{p,\Om}(\A)$. On the other hand if $\A$ is not an elliptic operator then we do not have $\text{dom}_{p,\Om}(\A)\subsetneq W^{1,p}_{loc}(\Om,E)$.

If $U$ is a $C^1$-domain, and the coefficients of $\A$ satisfies the assumptions of Rauch's trace theorem, Theorem \ref{thm:Rauch1}, then $\Aa(\cdot,\nu(\cdot))u\in W^{-1/2,2}(\dv \Om,\sigma;F)$. This allows us to view 
\begin{align*}
\int_{\dv U}\Aa(x,\nu(x))u(x)d\sigma(x)=\sum_{k=1}^m\int_{\dv U}\langle \Aa(x,\nu(x))u(x),f_m\rangle f_md\sigma(x)
\end{align*}
for $u\in C^1(\Om,E)$ and $\{f_k\}_k$ any ON-basis for $F$ as an $F$-valued duality, so that 
\begin{align*}
\int_{\dv U}\Aa(x,\nu(x))u(x)d\sigma(x):=\sum_{k=1}^m\langle\Aa(\cdot,\nu(\cdot))u,f_k \rangle f_k
\end{align*}
is well-defined whenever $\Aa(\cdot,\nu(\cdot))u\in W^{-1/2,2}(\dv \Om,\sigma;F)$. Furthermore, by Rauch's continuity theorem, Theorem \ref{thm:RauchCont}, this duality is continuous in the $W^{-1/2,2}$-norm with respect to deformations of $\dv U$ by nearby foliations. Obviously, the continuity is false in general with respect to the $L^2$-norm on $\dv U$. In this sense Stokes theorem holds for all $C^1$ domains when $p=2$.

%==================NEW SECTION=================================================================

\subsection{\sffamily Applications of Fuglede's Flux Extensions}

Using Fugelede's characterization in Theorem \ref{thm:FugledeChar} in combination with Fugelede' equivalence theorem one can give the following generalizations of Cauchy's integral theorem and Morera's theorem.

\begin{Thm}[Generalized Cauchy integral theorem]
\label{thm:CauchyInt}
Let $1\leq p\leq \infty$ and let $\Om \subset \R^n$ be a domain. Let $u\in L^p(\Om,E))$ for $1\leq p\leq \infty$. Assume that $\A_w u=0$. Then for $p$-a.e.Green domain $U\Subset \Om$
\begin{align}\label{eq:CauchyInt}
\int_{\dv U}\Aa(x,\nu(x))u(x)d\sigma(x)+\int_{U}(B(x)-\text{div}\, \Aa(x))u(x)dx=0.
\end{align}
\end{Thm}

\begin{proof}
By Fuglede's equivalence theorem $\A_f=\A_w$. Since $\A_w u=0$, $u\in \text{dom}_{p,\Om}(\A)$. By the definition of flux extensions the theorem follows.  
\end{proof}

\begin{Thm}[Generalized Morera theorem]
\label{thm:Morera}
Let $\Om \subset \R^n$ be a domain and let $u\in L^p(\Om,E)$ for $1\leq p\leq \infty$. Assume that for $p$-a.e.Green domain $U\Subset \Om$
\begin{align}\label{eq:Morera}
\int_{\dv U}\Aa(x,\nu(x))u(x)d\sigma(x)+\int_{U}(B(x)-\text{div}\, \Aa(x))u(x)dx=0.
\end{align}
Then $\A u=0$ in the sense of distributions. 
\end{Thm}

\begin{proof}
If \eqref{eq:Morera} holds for $p$-a.e. Green domain $U\Subset \Om$ then both Fuglede's conditions $\mathbf{A}_p$ and $\mathbf{B}$ are satisfied. Hence $u\in \text{dom}_{p,\Om}(\A_f)$. By Fuglede's equivalence theorem $\A_f=\A_w$. By the definition of flux extensions $\A_w u=0$. 
\end{proof}

To put the above generalized Morera's theorem in perspective, we recall the classical Morera's theorem from complex analysis. Let $f\in C(\Om,\C)$ for some domain $\Om \subset \C$. If 
\begin{align}\label{eq:MoreraC}
\frac{1}{2\pi i}\int_\gamma f(z)dz=0
\end{align}
for every closed piece-wise $C^1$-curve $\gamma$, then $f$ is holomorphic in $\Om$. Let us rewrite the integral \eqref{eq:MoreraC} in a more transparent way from the point of view of Stokes theorem. Note that $dz=-i\nu(z)d\sigma(z)$. Thus 
\begin{align*}
\frac{1}{2\pi i}\int_\gamma f(z)dz&=\frac{1}{2\pi i}\int_\gamma f(z)(-i\nu(z))d\sigma(z)=-\frac{1}{\pi}\int_\gamma \frac{\nu(z)}{2}f(z)d\sigma(z)\\&=-\frac{1}{\pi}\int_\gamma \Aa_{\overline{\dv}}(\nu(z))f(z)d\sigma(z), 
\end{align*}
where $\Aa_{\overline{\dv}}(\nu)$ is the symbol of the Cauchy-Riemann operator $\overline{\dv}=\frac{1}{2}(\dv_x+i\dv_y)$. The classical Morera's theorem has various sharpenings. For example, we do not need to test condition \label{eq:MoreraC} for all closed piecewise $C^1$-curves. It is enough that \eqref{eq:MoreraC} holds for every triangle or every circle contained in $\Om$. Furthermore, the continuity condition on $f$ can be relaxed to $L^1_{loc}(\Om,\C)$ by a theorem of Royden \cite{Royden}, and that furthermore \eqref{eq:MoreraC} only needs to hold for a.e. rectangle contained in $\Om$. With this in mind it is natural to ask if something similar is true for arbitrary first order operators $\A$. For example, when is it enough that \eqref{eq:Morera} hold only for a finite dimensional subsystem of Green domains? Let $\mathscr{F}$ be any collection of Green domains of $\Om$ together with their complements such that they form a semiring of sets and in addition generate the Borel $\sigma$-algebra $\mathscr{B}(\Om)$. For example, all open and closed cubes or simplexes in $U$ and countable disjoint unions of such in $\Om$. By Theorem \ref{thm:ExtensionM}, the countably additive set function 
\begin{align*}
\alpha_u(U)=\int_{\dv U}\Aa(x,\nu(x))u(x)d\sigma(x)+\int_{U}(B(x)-\text{div}\, \Aa(x))u(x)dx
\end{align*} 
for any $U\in \mathscr{F}$ admits a unique extension to a regular Borel measure $\mu_\alpha$ on $\mathscr{B}(\Om)\otimes F$. Furthermore, by definition of flux extensions there exists a $v\in L^p(\Om,E)$ such that 
\begin{align}\label{eq:Borel}
\alpha_u(U)=\int_Uv(x)dx
\end{align}
for every $U\in \mathscr{F}$. Since $\mathscr{F}$ generates the Borel $\sigma$-algebra by assumption this implies that \eqref{eq:Borel} holds for all $U\in \mathscr{B}(\Om)$. This also implies that if $(0,v)\in \text{gr}(\mathscr{A}_f)$, then the relation $\int_Uv(x)dx=0$ holds for $p$-a.e. $U\in \mathscr{F}$. By definition this implies that $\int_Uv(x)dx=0$ for every Borel set $U\Subset \Om$ and hence $v(x)=0$ for a.e. $x\in \Om$. Hence Theorem \ref{thm:CompleteFuglede} holds in this case as well. Thus in the generalized Morera's theorem, Theorem \ref{thm:Morera} it is sufficient to assume that \eqref{eq:Morera} holds for $p$-a.e. $U\in \mathscr{F}$, where $U\in \mathscr{F}$ system of Green sets with the stated properties above. In particular all cubes or simplexes suffice. We can now combine Theorem \ref{thm:Morera} with Theorem \ref{thm:Div3} to give the following results about continuous extension of solutions across hypersurfaces that are Ahlfors-David regular. 

\begin{Thm}[Continuous extensions across hypersurfaces.]
\label{thm:HyperExt1}
Let $\Om\subset \R^n$ be a domain and let $u\in C(\Om,E)$. Let $\Sigma$ be an Ahfors-David regular hypersurface that such that $\Om \setminus \Sigma$ has two connected components $\Om_+$ and $\Om_-$. Let $\A$ be a first order operator satisfying the standard assumption and assume that $\A u =0$ in $\Om_+$ and $\Om_-$. Then $\A u=0$ in $\Om$.
\end{Thm}

\begin{proof}
By Theorem \ref{thm:Exc} $p$-a.e open cube $Q$ intersects $\Sigma$ such that $\mathcal{H}^{n-1}(\dv Q\cap \Sigma)=0$. Hence for $p$-a.e. $Q$, every connected component $U$ of $Q\setminus \Sigma$ is an Ahlfors-David regular domain such that $u$ satisfies the assumptions of Theorem \ref{thm:Div3} on $\dv U$ and we can apply Stokes Theorem. By assumption on $u$ and Theorem \ref{thm:Morera}, 
\begin{align*}
\int_{\dv U}\Aa(x,\nu(x))u(x)d\sigma(x)+\int_{U}(B(x)-\text{div}\, \Aa(x))u(x)dx=0
\end{align*} 
Let $V$ by any adjacent component of $U$. Then on $\Gamma=\dv U\cap \dv V$,
\begin{align*}
\int_{\Gamma\cap \dv U}\Aa(x,\nu(x))u(x)d\sigma(x)+\int_{\Gamma\cap \dv V}\Aa(x,-\nu(x))u(x)d\sigma(x)=0,
\end{align*} 
since the outward pointing unit normals have opposite directions. Hence if $\{U_j\}_j$ are the (countable) connected components of $Q\setminus \Sigma$, then 
\begin{align*}
\sum_{j} \int_{\dv U_j}\Aa(x,\nu(x))u(x)d\sigma(x)=\int_{\dv Q}\Aa(x,\nu(x))u(x)d\sigma(x).
\end{align*} 
Thus, for $p$-a.e. cube $Q$
\begin{align*}
\int_{\dv Q}\Aa(x,\nu(x))u(x)d\sigma(x)+\int_{Q}(B(x)-\text{div}\, \Aa(x))u(x)dx=0.
\end{align*} 
By Theorem \ref{thm:Morera}, $\A u=0$ in $\Om$. 
\end{proof}
Note that already for smooth hypersurfaces there are cubes such that $Q\cap \Sigma$ can be any closed subset of $\dv Q$, in particular $\mathcal{H}^{n-1}(\dv Q\cap \Sigma)>0$. For such cubes $Q$ the connected components need not be Ahlfors-David regular domain and we cannot apply Stokes theorem. We are saved by the fact that it is enough that this holds for $p$-a.e. cube.

In case $n=2$ and $\A=\overline{\dv}$, and $\Sigma=\gamma$ is simple rectifiable Jordan curve Theorem \ref{thm:HyperExt1} is classical result for holomorphic functions. Note that the assumption that $\gamma$ is rectifiable cannot be dispensed with, see the discussion in \cite[p. 122]{Zalc}.

\begin{Thm}[Extensions across hyperplane.]
Let $1< p< \infty$, $\Om\subset \R^n$ be a domain and let $u\in \text{dom}_{p,\Om}(\A)$. Let $\Sigma_\nu$ be a hyperplane with unit normal $\nu$ such that $\Om \setminus \Sigma_\nu$ has two connected components $\Om_+$ and $\Om_-$. After a rotation and translation we may assume $\nu=e_n$ and $0\in \Sigma_{e_n}$ in the new coordinate system. Assume that $u$ has traces $u_+$ and $u_-$ on $\Sigma_\nu$ in the sense that for $x=(t,y)$, $t\in \R$, $y\in \R^{n-1}$
\begin{align}\label{eq:Hyper1}
u_+(0,y)=\text{w}-\esslim_{t\to 0^+} u(t,y), \quad u_-(0,y)=\text{w}-\esslim_{t\to 0^-} u(t,y),
\end{align}
in $L^p(\Sigma_{e_n},\mathcal{L}^{n-1},E)$ where w$-\esslim$ denotes essential weak limit. If $\A u=0$ in $\Om_\pm$ and 
\begin{align*}
\Aa(x,\nu)(u_+(x)-u_-(x))=0, \quad \sigma-\text{a.e. $x\in \Sigma_\nu$},
\end{align*}
then $\A u=0$ in $\Om$. 
\end{Thm}

\begin{proof}
For $p$-a.e cube $Q$, $\sigma(\dv Q\cap \Sigma)=0$. Let $Q_+$ and $Q_-$ denote the connected components of $Q\setminus \Sigma$. Let
\begin{align*}
Q_+^\eps=\{x=(t,x)\in Q: t>\eps\}\subset Q^+, \quad Q_-^\eps=\{x=(t,x)\in Q: t<\eps\}\subset Q^-.
\end{align*}
Since $\A u=0$ in $\Om_+$ and $\Om_-$ for a.e. $\eps>0$
\begin{align*}
\int_{\dv Q_\pm^\eps}\Aa(x,\nu(x))u(x)d\sigma(x)+\int_{Q_\pm^\eps}(B(x)-\text{div}\, \Aa(x))u(x)dx=0.
\end{align*} 
Since $u\in \text{dom}_{p,\Om}(\A)$ we can apply Theorem \ref{thm:Div4}. By assumption on $u$
\begin{align*}
&\esslim_{\eps\to 0^\pm}\int_{\dv Q_\pm^\eps}\Aa(x,\nu(x))u(x)d\sigma(x)+\int_{Q_\pm^\eps}(B(x)-\text{div}\, \Aa(x))u(x)dx\\
&=\int_{\dv Q_\pm}\Aa(x,\pm \nu)u_\pm(x)d\sigma(x)+\int_{Q_\pm}(B(x)-\text{div}\, \Aa(x))u(x)dx=0
\end{align*} 
Hence 
\begin{align*}
&\int_{\dv Q}\Aa(x,\nu(x))u(x)d\sigma(x)+\int_{Q}(B(x)-\text{div}\, \Aa(x))u(x)dx=\int_{\dv Q}\Aa(x,\nu)u(x)d\sigma(x)\\
&+\int_{\Sigma\cap Q}\Aa(x, \nu)(u_+(x)-u_-(x))d\sigma(x)+\int_{Q}(B(x)-\text{div}\, \Aa(x))u(x)dx\\
&=\int_{\dv Q_+}\Aa(x,\nu)u(x)d\sigma(x)+\int_{Q_+}(B(x)-\text{div}\, \Aa(x))u(x)dx\\
&+\int_{\dv Q_+}\Aa(x,\nu)u(x)d\sigma(x)+\int_{Q_+}(B(x)-\text{div}\, \Aa(x))u(x)dx=0.
\end{align*} 
By Theorem \ref{thm:Morera}, $\A u=0$ in $\Om$. 
\end{proof}

We now want to discuss removable singularities. To that end we need to introduce yet another Theorem of Fuglede from his paper \cite{F2}. 

\begin{Thm}\cite[Thm. 6]{F2}
\label{thm:Intersect}
Let $p\geq 1$ and $kp\leq n$. Let $S_{\text{Lip}}^k(B)$ denote the family of all $k$-dimensional Lipschitz hypersurfaces that intersect a given non-empty set $B\subset\R^n$. A neccessary, and for $p>1$, sufficient condition for $S_{\text{Lip}}^k(B)$ to be a $p$-exceptional is that the there exists a non-negative function $f\in L^p(\R^n)$ such that the $k$-th Riesz potential 
\begin{align*}
\mathcal{I}_k(f)(x)=\frac{1}{c_k}\int_{\R^n}\frac{f(y)}{\vert x-y\vert^{n-k}}dy, \quad c_k=\pi^{n/2}2^k\frac{\Gamma(k/2)}{\Gamma((n-k)/2)}
\end{align*}
satisfies $\mathcal{I}_k(f)(x)=+\infty$ for all $x\in B$ while $\mathcal{I}_k(f)\not\equiv +\infty$. If $p=1$, then the condition is also sufficient provided 
\begin{align*}
\int_{\R^n}f(x)\log_+f(x)dx<\infty. 
\end{align*}
\end{Thm}

The proof of Theorem \ref{thm:Intersect} in \cite{F2} is based on the following fundamental integral identity

\begin{align}\label{eq:Integral identity}
\int_{\mathbf{Gr}_k(\R^n)}\bigg(\int_{L}f(x)d\sigma(x)\bigg)d\mu_{k,n}(L)=\frac{\omega_k}{\omega_n}\int_{\R^n}\frac{f(x)}{\vert x\vert^{n-k}}dx
\end{align}
in combination with Calderón-Zygmund theory. Here $L\subset \R^n$ denotes a $k$-dimensional linear subspace of $\R^n$ and $\mathbf{Gr}_k(\R^n)$ denotes the Grassmann manifold of $k$-dimensional linear subspaces of $\R^n$. The invariant measure $\mu$ on $\mathbf{Gr}_k(\R^n)$ is defined from the Haar measure $\theta_n$ of the orthogonal group $\text{O}(\R^n)$ by fixing $L\in \mathbf{Gr}_k(\R^n)$ and defining 
\begin{align*}
\mu_{k,n}(A):=\theta_n(\{g\in \text{O}(\R^n): gw\in A\})
\end{align*}
for a Borel subset $A\subset \mathbf{Gr}_k(\R^n)$. $\mu_{k,n}$ is an invariant Radon measure on $\mathbf{Gr}_k(\R^n)$, that is 
\begin{align*}
\mu_{k,n}(gA)=\mu_{k,n}(A)
\end{align*}
for all $g\in \text{O}(\R^n)$. For a proof of \eqref{eq:Integral identity} see \cite{F4}.

\begin{Thm}[Removable singularities]
\label{thm:RemovSing}
Let $p\geq 1$ and $kp\leq n$. Let $\Om \subset \R^n$ be a bounded domain and let $B\subset \R^n$ be a closed subset such that $S_{\text{Lip}}^k(B)$ is $p$-exceptional. Let $u\in \text{dom}_{p,\Om}(\A)$ such that $\A u=0$ in $\Om \setminus B$ in the sense of distribution. If there exists a decreasing family of Green domains $U_\eps$, such that for $\eps>0$, $B\Subset U_\eps$ and 
\begin{align}\label{eq:RemovableC}
\esslim_{\eps \to 0^+}\bigg\vert \int_{\dv U_\eps} \Aa(x,\nu(x))u(x)d\sigma(x)\bigg\vert=0,
\end{align}
then $\A u=0$ in $\Om$. 
\end{Thm}

\begin{proof}
By the generalized Morera theorem \ref{thm:Morera} we need to show that $\int_{\dv U} \Aa(x,\nu(x))u(x)d\sigma(x)=0$ for $p$-a.e. Green domain $U$. Since the set $S_{\text{Lip}}^k(B)$ is $p$-exceptional we may assume that $\dv U$ does not intersect $B$. By taking $\eps$ sufficiently small we may assume that $(\dv U_\eps\cap U)\Subset U$. Then $V=U\setminus (U_\eps \cap U)$ is a Green domain not containing $B$. By assumption 
\begin{align*}
0=\int_{\dv V} \Aa(x,\nu(x))u(x)d\sigma(x)=\int_{\dv U} \Aa(x,\nu(x))u(x)d\sigma(x)+\int_{\dv U_\eps\cap U} \Aa(x,\nu(x))u(x)d\sigma(x).
\end{align*}
Letting $\eps\to 0^+$ and using our assumption it follows that $\int_{\dv U} \Aa(x,\nu(x))u(x)d\sigma(x)=0$. 
\end{proof}
Of course condition \eqref{eq:RemovableC} is rather implicit and is not really a measure theoretic condition purely regarding the size of $B$. On the other hand if $\A$ is not an elliptic operator this is not expected to be the case as removability should be connected to properties of $\A$. More precisely, by Hölder's inequality with $1/p+1/q=1$
\begin{align*}
&\bigg\vert \int_{\dv U_\eps} \Aa(x,\nu(x))u(x)d\sigma(x)\bigg\vert\leq \int_{\dv U_\eps} \Vert \Aa(x,\nu(x))\Vert \vert \Aa(x,\nu(x))^+\Aa(x,\nu(x))u(x)\vert d\sigma(x)\\
\leq &\bigg(\int_{\dv U_\eps} \Vert \Aa(x,\nu(x))\Vert^q d\sigma(x)\bigg)^{1/q} \bigg(\int_{\dv U_\eps}\vert \Aa(x,\nu(x))^+\Aa(x,\nu(x))u(x)\vert^p d\sigma(x)  \bigg)^{1/p}\\
\leq &\sup_{(x,\xi)\in\Om\times S^{n-1})}\Vert \Aa(x,\nu(x))\Vert(\mathcal{H}^{n-1}(\dv U_\eps))^{1/q} \esup_{\eps>0}\bigg(\int_{\dv U_\eps}\vert \Aa(x,\nu(x))^+\Aa(x,\nu(x))u(x)\vert^p d\sigma(x)  \bigg)^{1/p},
\end{align*}
where  $\Aa(x,\nu(x))^+$ denotes the Moore-Penrose inverse and $\Aa(x,\nu(x))^+\Aa(x,\nu(x))$ is the orthogonal projection onto $\big(\text{ker}(\Aa(x,\nu(x)))\big)^\perp$, see also section \ref{sec:NonEllip}. Thus if 
\begin{align*}
\esup_{\eps>0}\int_{\dv U_\eps}\vert \Aa(x,\nu(x))^+\Aa(x,\nu(x))u(x)\vert^p d\sigma(x) <+\infty
\end{align*}
and $\mathcal{H}^{n-1}(\dv U_\eps)\to 0$ as $\eps \to 0^+$, then condition \eqref{eq:RemovableC} is satisfied.

It is instructive to compare the removable singularity Theorem \ref{thm:RemovSing} with those of Harvey and Polking in \cite{HP} relying on the weak extension $\A_w$ instead and judicious choices of test functions.

Closely connected to Morera's theorem and removable singularities in complex analysis is the concept of \emph{analytic capacity} of a compact set $K\subset \C$ which is defined as follows. For any domain $\Om \subset \C$ we let $H^\infty(\Om)$ denote the space of bounded holomorphic functions on $\Om$. The analytic capacity of a compact set $K\subset$ is 
\begin{align}\label{eq:AnalyticCap}
\gamma(K):=\sup\bigg\{\bigg\vert \frac{1}{2 \pi i}\int_\Gamma f(z)dz\bigg\vert: f\in H^\infty(\C\setminus K), \,\,\, \Vert f\Vert_{L^\infty}\leq 1\bigg\}
\end{align}
where $\Gamma$ is any rectifiable curve surrounding $K$. It was proven by Ahlfors in \cite{Ah47} that a compact $K\subset \C$ is a removable singularity for all $f\in H^\infty(U\setminus K)$ where $U$ is any open set containing $K$ if and only if $\gamma(K)=0$. Furthermore, in the case when $\gamma(K)>0$, there exists a unique holomorphic function, called the \emph{Ahlfors map} $\alpha_K$ that maximizes \eqref{eq:AnalyticCap}. If $K=\dv \Om$ of a finitely connected domain $\Om$ all whose connected boundary components are Jordan curves (say $n$ of them), then the Ahlfors map is an $n-1$-branched covering of $\Om$ onto the unit disc. It is natural to extend the notion of analytic capacity to general first order operators $\A$ by defining 
\begin{align}\label{eq:AnalyticCap2}
\gamma_\A^p(K):=&\sup\bigg\{\bigg\vert \int_{\dv U} \Aa(x,\nu(x))u(x)d\sigma(x)+\int_{U}(B(x)-\text{div}\,\Aa(x))u(x)dx\bigg\vert\\&: \A u=0\,\,\, x\in \R^n\setminus K, \,\,\, \Vert u\Vert_{L^p(\R^n\setminus K,E)}\leq 1\bigg\}\nonumber
\end{align}
for any smooth bounded open domain $U$ containing $K$. For the analytic capacity to be well-defined, i.e., being in dependent of the choice of $U$, some regularity of solutions of $\A u$ in $\R^n \setminus K$ is required so that the trace of $u$ on $\dv U$ is always defined. It is natural to ask for which operators $\A$ other than the Cauchy-Riemann operator, if any, $\gamma_\A^p(K)$ is well-defined and that this still characterizes removable singularities, and if there exists a generalization of the Ahlfors map. Let us give a partial answer to the last question first. 
\begin{Def}\label{def:WeakHypo}
An operator $\A: C^\infty(\R^n,E)\to C^\infty(\Om,F)$ with smooth coefficients is called \emph{weakly hypoelliptic} if for any $U\subset \Om$ and any $u\in \mathcal{D}'(U,E)$ $\A u=0$ implies that $u\in C^\infty(U,E)$. 
\end{Def}
The concept of weak hypoellipticity was introduced in \cite{Bär13} and is strictly weaker than hypoellipticity, i.e. that for any $u\in \mathcal{D}'(U,E)$ such $\A u\in C^\infty(U,F)$ implies that $u\in C^\infty(U,E)$. An explicit example of a weakly hypoelliptic but not hypoelliptic operator given in \cite{Bär13} is the operator $\A: C^\infty(\R^2,\R)\to C^\infty(\R^2,\R^2)$ given by 
\begin{align*}
\A u(x)=\begin{bmatrix} x\\
-y
\end{bmatrix}\dv_xu+
\begin{bmatrix} y\\
x
\end{bmatrix}\dv_yu-\begin{bmatrix} 2\\
0
\end{bmatrix}u. 
\end{align*}

Weakly hypoelliptic operators satisfy a generalized version of Montel's theorem from complex analysis. 
\begin{Thm}[Generalized Montel Theorem]
\label{thm:Montel}
Let $\A: C^\infty(\Om,E)\to C^\infty(\Om,E)$ be weakly hypoelliptic. Then any bounded sequence $\{u_j\}_j\in L^1_{loc}(\Om,E)$ such that $\A u_j=0$ in $\mathcal{D}'(\Om,E)$ has a subsequence (which we do not relabel) $\{u_j\}_j$ which converges to some $u\in C^\infty(\Om,E)$ such that $\A u=0$ in $\Om$. 
\end{Thm}

Theorem \ref{thm:Montel} is proven in \cite[Thm. 4]{Bär13}.  We see that the analytic capacity $\gamma_\A^p(K)$ is well-defined for any weakly hypoelliptic operator $\A$. Using the generalized Montel theorem \ref{thm:Montel} one can prove that the there exists at least one function $u$ the maximizes \eqref{eq:AnalyticCap2}. Indeed, let $\{u_j\}_j\subset L^p(\R^n\setminus K,E)$, with $\Vert u_j\Vert_{L^p(\R^n\setminus K,E)}\leq 1$ and $\A u_j=0$ such that 
\begin{align*}
\lim_{j\to \infty}\bigg\vert \int_{\dv U} \Aa(x,\nu(x))u_j(x)d\sigma(x)+\int_{U}(B(x)-\text{div}\,\Aa(x))u_j(x)dx\bigg\vert=\gamma_\A^p(K). 
\end{align*}
By the generalized Montel theorem the maximizing sequence contains a convergent subsequence convergning to some $u\in C^\infty(\R^n\setminus K,E)$ such that $\A u=0$. Consider now the generalized analytic capacity itself. We assume that $K$ is a null set so that $L^p(\R^n\setminus K)=L^p(\R^n)$.  Hence $\A u\in \mathcal{E}'(\Om,F)$ with $\text{supp}(\A u)\subset K$. Since $u\in C^\infty(\R^n,E)$, $u$ satisfies the assumptions for Stokes theorem in the form \eqref{eq:StokesA}. Hence 
\begin{align*}
\bigg\vert \int_{\dv U} \Aa(x,\nu(x))u(x)d\sigma(x)+\int_{U}(B(x)-\text{div}\,\Aa(x))u(x)dx\bigg\vert=\vert (\A u,1)_F\vert
\end{align*}
and 
\begin{align}\label{eq:AnalyticCap3}
\gamma_\A^p(K):=&\sup\{\vert (\A u,1)_F\vert: \A u=0\,\,\, \text{ in } \R^n\setminus K, \,\,\, \Vert u\Vert_{L^p(\R^n,E)}\leq 1\}. 
\end{align}
It is clear that if $K$ is removable in $L^p(U,E)$ then $\A u=0$ in $\mathcal{E}'(U,F)$ for all $u$ and thus that $\gamma_\A^p(K)=0$. Let us now assume that $\gamma_\A^p(K)=0$. This implies that $\vert (\A u,1)_F\vert=0$ for all $u\in L^p(\R^n,E)$ such that $\A u=0$ in $\R^n\setminus K$. The question now is if this actually implies $\A u=0$?  All proofs of this fact that the author is aware of in the case of holomorphic functions rely on the function algebra structure of holomorphic functions and does not generalize to this case. Trying to resolve these issues in this survey paper will lead us too far astray, and we leave this question as an interesting future research direction. Finally, to conclude the discussion on generalized analytic capacity, in the case of the classical analytic capacity in the plane, one wants a purely geometric characterization of analytic capacity. This is essentially Vitushkin's conjecture. It says that 
\begin{align*}
\gamma(K)=0 \,\,\, \Longleftrightarrow \,\,\, \int_{0}^\pi \mathcal{H}^1(\text{proj}_\theta(K))d\theta=0.
\end{align*}
If $\text{dim}_{\mathcal{H}}(K)\neq 1$ the conjecture is true. If $\text{dim}_{\mathcal{H}}(K)=1$ the conjecture is true if $K$ is $\mathcal{H}^1-\sigma$-finite by the work of G. David in \cite{David} and X. Tolsa in \cite{Tolsa} and false otherwise by a counterexample by P. Mattila in \cite{Matti}. In particular, the analytic capacity is in some sense independent of the Cauchy-Riemann operator $\overline{\dv}$. One can therefore ask how the generalized analytic capacity depends on the choice of operator $\A$? Again we leave this as an open problem.

%==================NEW SECTION=================================================================

\section{\sffamily On the pointwise value of flux operators}\label{sec:Point}

%============NEW SUBSECTION=================================================================
\subsection{\sffamily  Pointwise value of flux operators and Lesbegue differentiation theorem}

We begin by recalling the Lesbegue differentation theorem.

\begin{Thm}[Lesbegue]
Let $f\in L^1_{loc}(\Om,E)$. Then 
\begin{align*}
\lim_{r\to 0^+}\frac{1}{\vert B_r(x)\vert}\int_{B_r(x)}f(y)dy=f(x) \text{ in $L^1_{loc}$ and a.e. $x$.}
\end{align*}
\end{Thm}

The representative $f^\star$ in each equivalence class in $L^1_{loc}(\Om,E)$ defined pointwise according to 
\begin{align}\label{eq:PRep}
f^\star(x)=\lim_{r\to 0^+}\frac{1}{\vert B_r(x)\vert}\int_{B_r(x)}f(y)dy
\end{align}
is called the \emph{precise representative} of $f$ (\cite[Def. 1.26]{EG}).

\begin{Def}
Let $1\leq p<+\infty$. A point $x\in \Om$ is a $p$-Lesbegue point of $f\in L^p_{loc}(\Om,E)$,  if 
\begin{align*}
\lim_{r\to 0^+}\frac{1}{\vert B_r(x)\vert}\int_{B_r(x)}\vert f(x)-f^\star(y)\vert^pdy=0.  
\end{align*}
The set of $p$-Lebesgue points of $u$ is denoted by $\mathscr{L}_f^p$. 
\end{Def}
A well-known result in integration theory shows that the set of $p$-Lebesgue points has full-measure. Furthermore, the limit \eqref{eq:PRep} exists for all $x\in \mathscr{L}_u^p$. It is natural to ask if the same result holds if we replace the family of balls $\{B_r(x)\}_r$ with some other family of subsets of $\R^n$? 
\begin{Def}\label{def:RegularF}
A family of measurable subsets $\mathscr{F}$ of $\R^n$ is \emph{regular} at $x\in \R^n$ if
\begin{itemize}
\item[(i)] The sets are bounded and have positive measure.
\item[(ii)] There is a sequence $\{U_i\}_i\subset \mathscr{F}$ such that $\vert U_i\vert\to 0$ as $i\to \infty$. 
\item[(iii)] There is a constant $c>0$ such that for all $U\in \mathscr{F}$, $\vert U\vert\geq c\vert B_U\vert$, where $B_U$ is the smallest closed ball contained in $U$ centered at $x$. 
\end{itemize}
\end{Def}

Using the notion of regular family of measurable subsets we have the following extension of Lesbegue differentation theorem.

\begin{Thm}
Let $f\in L^p_{loc}(\Om,E)$ and let $\mathscr{F}$ be a family of measurable subsets \emph{regular} at $x\in \R^n$. Then for all $x\in \mathscr{L}^p_f$,
\begin{align}
\lim_{\substack{\vert U\vert \to 0\\U\in \mathscr{F}}}\frac{1}{\vert U\vert}\int_{U}f(y)dy=f(x)
\end{align}
and 
\begin{align}\label{eq:UniLSet}
\lim_{\substack{\vert U\vert \to 0\\U\in \mathscr{F}}}\frac{1}{\vert U\vert}\int_{U}\vert f(y)-f(x)\vert^p dy=0
\end{align}
for a.e. $x$ and in $L^p_{loc}$. 
\end{Thm}

For a proof of these facts we refer the reader to \cite[Sec. 1.8, p. 11]{SteinBook} and \cite[Thm. 2.28, p. 11]{Ha}. We see in particular that the notion of Lebesgue point is \emph{independent} of the choice of regular family $\mathscr{F}$ and choice of representative in the equivalence class in $L^p(\Om,E)$. Indeed, assume that $U\in \mathscr{F}$. By scaling we may assume that $U\subset B_r(x)$. Then, Jensen's inequality implies
\begin{align*}
\bigg\vert \frac{1}{\vert U\vert}\int_{U} f(y)dy-f(x)\bigg\vert^p\leq  \frac{1}{\vert U\vert}\int_{U} \vert f(y)-f(x)\vert^p dy\leq \frac{c}{\vert B_r(x)\vert}\int_{B_r(x)} \vert f(y)-f(x)\vert^p dy
\end{align*}
from which the statement follows. Combining the Lesbegue differentation theorem with Fuglede's flux extension, the pointwise value of $\A_f u(x)=\A_wu(x)$ for $u\in \text{dom}(\A_w)$ is then give as a limit. 

\begin{Thm}\label{thm:LimitFuglede}
Let $\mathscr{F}$ be a regular finite dimensional family of measurable sets of $\R^n$ at $x$, such that each $U\in \mathscr{F}$ is a Green domain. Let $\Om \subset \R^n$ be a domain, $1\leq p\leq \infty$ and $u\in \text{dom}_{p,\Om}(\A_w)$. Then for any family of regular domains  
\begin{align}\label{eq:FLimit}
\esslim_{\substack{\vert U\vert \to 0\\U\in \mathscr{F}}}\frac{1}{\vert U\vert}\int_{\dv U}\Aa(y,\nu(y))u(y)d\sigma(y)+\frac{1}{\vert U\vert}\int_{U}(B(x)-\text{div}\, \Aa(x))u(x)dx=\A_w u(x)
\end{align}
for all $x\in \mathscr{L}_{\A u}^p$ and in $L^p_{loc}(\Om,F)$. In particular
\begin{align}\label{eq:FLimit2}
\esslim_{\eps \to 0^+}\frac{1}{\vert B_\eps(x)\vert}\int_{\dv B_\eps(x)}\Aa(y,\nu(y))u(y)d\sigma(y)+\frac{1}{\vert B_\eps(x)\vert}\int_{B_\eps(x)}(B(y)-\text{div}\, \Aa(y))u(y)dy=\A _w u(x).
\end{align}
\end{Thm}

\begin{proof}
This follows directly from the Lebesgue differentiation theorem, the definition of flux extensions and Fuglede's equivalence theorem, Theorem \ref{thm:EquivalenceFuglede}. 
\end{proof}

\begin{rem}
Note that if in the limit we were to consider the two integrals in the left-hand side of \eqref{eq:FLimit2} separately, we would get that the limit exists for all $x\in \mathscr{L}_u^p\cap \mathscr{L}^p_{\A u}$. This is still a set of full measure, but $\mathscr{L}_u^p$ could still be much smaller that $\mathscr{L}_{\A u}^p$. This suggests that studying the limit of the integrals can be quite subtle and rely on cancellation effects between the two terms. 
\end{rem}

\begin{Def}[Pure flux operators]
A first order operator $\A$ satisfying the standard assumptions is called a \emph{pure flux} operator if 
\begin{align*}
B(x)-\text{div}\, \Aa(x)=0. 
\end{align*}
\end{Def}
In particular, any homogeneous constant coefficient first order operator is a pure flux operator.  For a pure flux operator, $\A u(x)$ becomes a limit of pure boundary integral.

Of central importance for establishing these facts is the \emph{Hardy-Littlewood maximal function},
\begin{align*}
\mathcal{M}f(x)=\sup_{r>0}\frac{1}{\vert B_r(x)\vert}\int_{B_r(x)}\vert f(y)\vert dy
\end{align*}
for $f\in L^1_{loc}(\R^n,E)$. Also for $f\in L^1_{loc}(\Om,E)$ and for a regular family $\mathscr{F}$ we set 
\begin{align*}
\mathcal{M}_{\Om}f(x)&=\sup_{\substack{r>0\\ B_r(x)\in \Om}}\frac{1}{\vert B_r(x)\vert}\int_{B_r(x)}\vert f(y)\vert dy\\
\mathcal{M}_{\mathscr{F},\Om}f(x)&=\sup_{U\in \mathscr{F}}\frac{1}{\vert U\vert}\int_{U}\vert f(y)\vert dy. 
\end{align*}
In particular we note that $\mathcal{M}_{\mathscr{F},\Om}f(x)\leq c^{-1} \mathcal{M}f(x)$ for some $c>0$ and all $f$.
In the next section we will consider a more refined notion of maximal operator associated to $\A_wu$.

%============NEW SUBSECTION=================================================================
\subsection{\sffamily Maximal operator for $\A$}\label{sec:MaxOp}

Let $u\in L^p(\Om,E)$ with $p\in [1,+\infty]$. By Fubini's theorem $u\vert_{\dv B_\eps(x)}\in L^p(\dv B_\eps(x),d\sigma)$ for a.e. $\eps>0$ such that $B_\eps(x)\Subset \Om$. Thus for every $x\in \Om$ and a.e. $0<\eps<\text{dist}(x,\dv \Om)$ the \emph{truncated operator} 
\begin{align}
\mathscr{A}_\eps u(x):=\frac{1}{\vert B_\eps(x)\vert}\int_{\dv B_\eps(x)}\mathbb{A}(y,\nu(y))u(y)d\sigma(y)+\frac{1}{\vert B_\eps(x)\vert}\int_{B_\eps(x)}(B(y)-\text{div}\, \Aa(y))u(y)dy,
\end{align}
 is well-defined. The truncated operator $\mathscr{A}_\eps$ has an appealing analogy with singular integrals. Assume we have a Calderón-Zygmund operator given by the principal valued integral
\begin{align*}
Tu(x)=\text{p.v.}\int_{\R^n}K(x,y)u(y)dy
\end{align*}
with kernel $K$. One defines the truncated operator according to 
\begin{align*}
T^\eps u(x)=\int_{\R^n\setminus B_\eps(x)}K(x,y)u(y)dy. 
\end{align*}
Since we are only interested in the behaviour of $\A_\eps u$ in the interior of $\Om$, we will assume henceforth in this section that $\Om= \R^n$. This we may, since after multiplying any $u\in \text{dom}_{p,\Om}(\A_w)$ with some $\phi\in C^\infty_0(\Om)$, we have $\phi u\in \text{dom}_{p,\R^n}(\A_0)$ since first order operators are of local type. 

Associated to the truncated operator we define a \emph{maximal operator} $\mathscr{A}^\star u(x)$ through 
\begin{align*}
\mathscr{A}^\star u(x):&=\esup_{\eps>0} \frac{1}{\vert B_\eps(x)\vert}\bigg\vert\int_{\dv B_\eps(x)}\mathbb{A}(y,\nu(y))u(y)d\sigma(y)+\frac{1}{\vert B_\eps(x)\vert}\int_{B_\eps(x)}(B(y)-\text{div}\, \Aa(y))u(y)dy\bigg\vert
\end{align*}
where the supremum here is the essential supremum. This again has an analogy with singular integrals through the maximal transform 
\begin{align*}
T^\star u(x):=\sup_{\eps>0}\bigg\vert\int_{\R^n\setminus B_\eps(x)}K(x,y)u(y)dy\bigg\vert
\end{align*}
 Here we remark that it is not immediate that $\A^\star$ is Lebesgue measurable since the supremum is over an uncountable set. This however follow from the fact that for a.e. fixed $\eps>0$ and every $x\in \R^n$ the integral average
\begin{align*}
\frac{1}{\vert B_\eps(x)\vert}\int_{B_\eps(x)}\A u(y)dy
\end{align*}
is continuous in $x$. Since the supremum of a family of continuous function is lower-semicontinuous it follows that $\A^\star  u$ is Lebesgue measurable. Furthermore, we obviously have the inequality 
 \begin{align*}
 \A^\star u(x)\leq \mathcal{M}( \A u)(x)
 \end{align*}
 for all $x$. 
 For the operator $T^\star$ we have the classical Cotlar's inequality
 \begin{align*}
 T^\star u(x)\leq C(\mathcal{M}(T f)(x)+\mathcal{M}(f)(x)), \quad x\in \R^n.
 \end{align*}

In particular we note that by definition, if $u\in \text{dom}_{p,\R^n}(\A_w)$ we have for a.e. $\eps>0$
\begin{align*}
\vert \mathscr{A}_\eps u(x)\vert &\leq\frac{1}{\vert B_\eps(x)\vert}\bigg\vert \int_{\dv B_\eps(x)}\mathbb{A}(y,\nu(y))u(y)d\sigma(y)\bigg\vert+\frac{1}{\vert B_\eps(x)\vert}\bigg\vert \int_{B_\eps(x)}(B(y)-\text{div}\, \Aa(y))u(y)dy\bigg\vert\\ 
&\leq \frac{1}{\vert B_\eps(x)\vert}\int_{ B_\eps(x)}\vert \A u(y)\vert dy +\Vert \Vert B-\text{div}\,\Aa\Vert \Vert_{L^\infty(\Om)}\frac{1}{\vert B_\eps(x)\vert}  \int_{B_\eps(x)}\vert u(y)\vert dy
\end{align*}
and so 
\begin{align*}
\A^\star u(x)\leq C'(\mathcal{M}(\A u)(x)+\mathcal{M}( u)(x))
\end{align*}
for some constant $C'$ independent of $u$ in complete analogy with Cotlar's inequality for singular integrals.

In fact the similarities with singular integral operators goes further. 
\begin{Def}[Directional Hilbert transform]
Let $\xi\in S^{n-1}$. Then for any $f\in L^p(\R^n,E)$ we define the \emph{directional Hilbert transform}
\begin{align*}
\mathcal{H}_{\xi}(f)(x)=\frac{1}{\pi}\pv\int_{-\infty}^{\infty}\frac{f(x-t\xi)}{t}dt. 
\end{align*}
\end{Def}
Recall the following result for homogeneous singular integral operators.

\begin{Thm}[Method of rotations]
Assume that $K\in L^1(S^{n-1})$ is odd and consider the associated singular integral operator 
\begin{align*}
T_Ku(x)=\pv\int_{\R^n}\frac{K(x-y)}{\vert x-y\vert^n}u(y)dy.
\end{align*}
Then for any $u\in L^p(\R^n,E)$ and $1<p<\infty$ 
\begin{align}\label{eq:SingF}
T_Ku(x)=\frac{\pi}{2}\int_{S^{n-1}}K(\xi)\mathcal{H}_{\xi}(u)(x)d\sigma(\xi),
\end{align}
where $\mathcal{H}_{\xi}(u)(x)$ is well-defined for a.e. $\xi\in S^{n-1}$ and belongs to $L^p(\R^n,E)$. 
\end{Thm}

For a proof see \cite[Thm. 5.2.7, p. 339]{Gar}.

We now compare formula \eqref{eq:SingF} with the corresponding formula for a constant coefficient first order PDE operator 
\begin{align*}
\A u(x)=\frac{1}{\omega_n}\int_{S^{n-1}}\Aa(\xi)u'(x;\xi)d\sigma(\xi),
\end{align*}
and we see that the directional derivative $u'(x;\xi)$ takes the role of the directional Hilbert transform in \eqref{eq:SingF}. In both cases we have the cancellation conditions 
\begin{align}\label{eq:Cancel1}
\int_{S^{n-1}}\Aa(\xi)d\sigma(\xi)=0, \quad \int_{S^{n-1}}K(\xi)d\sigma(\xi)=0. 
\end{align}

Finally, as an alternative to the study of the pointwise convergence of $\A^\eps u$ in Theorem \ref{thm:LimitFuglede}, we could instead use $\A^\star u$ by considering the oscillations 
\begin{align*}
\omega_\A u(x):=\vert \esslimsup_{\eps \to 0^+} \vert \A^\eps u(x)\vert-\essliminf_{\eps \to 0^+}\vert \A^\eps u(x)\vert \vert
\end{align*}
and use that $\omega_\A u(x)\leq 2\A^\star u(x)$. The details of such an argument can be found in \cite[Thm. 13.34, p. 528]{Tah}.

%============NEW SECTION=================================================================

\subsection{\sffamily Convergence in the sense of distributions}

In the previous sections we studied the convergence of $\A^\eps u$ to $\A u$ in the case when $u\in \text{dom}_{p,\Om}(\A)$. It is natural to ask if for any $u\in L^p(\Om,E)$ we have 
$\A^\eps u\to \A u\in \mathcal{D}'(\Om,E)$ in the sense of distributions. For this to make sense we will need that the coefficients of $\A$ are smooth since otherwise $\A u$ does not make sense as a distribution. To further simplify the question we will in this section assume that $\A$ has constant coefficients. 

\begin{Lem}\label{lem:TestFunc}
Let $\phi\in C^\infty_0(\Om,F)$. Then for all $0<2\eps <\text{dist}(\text{supp}(\phi), \dv \Om)$, $\A_\eps \phi\in C_0^\infty(\Om,F)$ and $\displaystyle \lim_{\eps\to 0^+}\A_\eps \phi=\A \phi$ in $C^\infty_0(\Om,F)$. 
\end{Lem}

\begin{proof}
For every $x\in \Om$
\begin{align*}
\A_\eps \phi(x)=\frac{1}{\omega_n}\int_{S^{n-1}}\Aa(\xi)\frac{(\phi(x+\eps \xi)-\phi(x))}{\eps}d\sigma(\xi)=0
\end{align*}
whenever $x\in \{y\in \Om: \text{dist}(y, \text{supp}(\phi)\geq \eps\}$, and so $\A_\eps u$ is compactly supported in $\Om$. In addition, a Taylor expansion together with the cancellation condition \eqref{eq:Cancel1} yield
\begin{align*}
& \A_\eps \phi(x)-\A_\eps \phi(y)=\frac{1}{\omega_n}\int_{S^{n-1}}\Aa(\xi)\frac{(\phi(x+\eps \xi)-\phi(y+\eps \xi))}{\eps}d\sigma(\xi)\\
&=\frac{1}{\eps\omega_n}\int_{S^{n-1}}\Aa(\xi)D\phi(x)(y-x)d\sigma(\xi)+O(\vert y-x\vert)=\frac{1}{\eps\omega_n}\int_{S^{n-1}}\Aa(\xi)d\sigma(\xi)D\phi(x)(y-x)+O(\vert y-x\vert)\\
&=O(\vert y-x\vert)
\end{align*}
independent of $\eps$, which implies that $\A_\eps \in C_0(\Om,E)$. Repeating the same argument with $\dv^{\alpha}\A_\eps \phi=\A_\eps \dv_{\alpha}\phi$ for any $\dv^{\alpha}=\dv^{\alpha_1}_{x_1}....\dv^{\alpha_n}_{x_n}$ shows that $\A^\eps \phi \in C^k_0(\Om,E)$ for any $k\geq 0$, which completes the proof. 
 \end{proof}

\begin{Lem}\label{lem:TruncatedINTPART}
If $u\in L^p(\Om,E)$, $1\leq p\leq \infty$ with $\text{supp}(u)\Subset \Om$ and $\phi \in C^\infty_0(\Om,F)$, then for a.e. $0<\eps <\text{dist}(\text{supp}(u),\Om)$ 
\begin{align*}
\langle \A_\eps u,\phi\rangle=-\langle u,\A_\eps^\ast \phi\rangle. 
\end{align*}
\end{Lem}

\begin{proof}
By Fubini's theorem
\begin{align*}
\langle \A_\eps u,\phi\rangle&=\frac{1}{\omega_n}\int_{S^{n-1}}\int_{\R^n}\eps^{-1}\langle \Aa(\xi)u(x+\eps \xi),\phi(x)\rangle dxd\sigma(\xi)\\
&=\frac{1}{\omega_n}\int_{S^{n-1}}\int_{\R^n}\eps^{-1}\langle \Aa(\xi)u(y),\phi(y-\eps \xi)\rangle dyd\sigma(\xi)\\
&=\int_{\R^n}\frac{1}{\omega_n}\int_{S^{n-1}}\eps^{-1}\langle u(y),\Aa(\xi)^\ast\phi(y-\eps \xi)\rangle d\sigma(\xi)dy\\
&=\int_{\R^n}\frac{1}{\omega_n}\int_{S^{n-1}}\eps^{-1}\langle  u(y),\Aa(-\xi')^\ast\phi(y+\eps \xi')\rangle d\sigma(\xi')dy\\
&=-\langle u,\A_\eps^\ast \phi\rangle.
\end{align*}
\end{proof}

\begin{Thm}
If $u\in L^p(\Om,E)$, $1\leq p\leq \infty$ with $\text{supp}(u)\Subset \Om$, then 
\begin{align*}
\esslim_{\eps \to 0^+}\A_\eps u=\A u
\end{align*}
in $\mathcal{D}'(\Om,F)$ 
\end{Thm}

\begin{proof}
First note that $\A_\eps u\in L^1_{loc}(\Om,F)$ for a.e. $\eps>0$. By Lemma \ref{lem:TruncatedINTPART},
\begin{align*}
\langle \A_\eps u,\phi\rangle=-\langle u,\A_\eps^\ast \phi\rangle. 
\end{align*}
By Lemma \ref{lem:TestFunc}, $\A_\eps^\ast \phi \in C^\infty_0(\Om,E)$, and by Lebesgue dominated convergence theorem 
\begin{align*}
\lim_{\eps\to 0^+}\langle \A^\eps u(x),\phi(x)\rangle=\lim_{\eps\to 0^+}-\langle u(x),\A_\eps^\ast \phi(x)\rangle =-\langle u,\A^\ast \phi \rangle =\langle \A u,\phi\rangle. 
\end{align*}
\end{proof}

%==================NEW SECTION=================================================================

\section{\sffamily Limit formula for flux operators}\label{sec:Limit}

%==================NEW SUBSECTION=================================================================

\subsection{\sffamily Notions of directional differentiability}

While Theorem \ref{thm:LimitFuglede} is highly satisfactory from the point of view of the existence of the limit \eqref{eq:FLimit2}, it does not ``explain'' in a more hands on way why the limit 
\begin{align*}
\lim_{\eps\to 0^+}\frac{1}{\vert B_\eps(x)\vert}\int_{B_\eps(x)}\Aa(x,\nu(x))u(x)d\sigma(x)
\end{align*}
exists, nor what information about $u$ is actually encoded in the limit. In particular, after a change of variables we find that 
\begin{align*}
\lim_{\eps\to 0^+}\frac{1}{\vert B_\eps(x)\vert}\int_{\dv U}\Aa(y,\nu(y))u(y)d\sigma(y)=\lim_{\eps\to 0^+}\frac{1}{\omega_n}\int_{\dv U}\Aa(x+\eps \xi,\xi)\frac{u(x+\eps \xi)}{\eps}d\sigma(\xi).
\end{align*}
Since $\eps^{-1}u(x+\eps \xi)$ is divergent as $\eps\to 0^+$ the existence of the limit depends on a suitable cancellation property of the integral. The understanding of this cancellation property will be the main topic of later subsections.  To answer these questions we will first introduce some notions of directionally differentiability that will be useful in later discussions. 

\begin{Def}[One sided Gâteaux directional derivatives]
\label{def:DerG}
$u\in L^p(\Om,E)$, $1\leq p\leq \infty$, is said to have a one sided directional derivative at $x\in \Om$ if for every $\xi\in S^{n-1}$ the limit 
\begin{align}
u'(x;\xi):=\lim_{\eps\to 0^+}\frac{u(x+\eps \xi)-u(x)}{\eps}
\end{align}
exists. 
\end{Def}

\begin{ex}
Let $P$ be an $E$-valued $m$-homogeneous polynomial and $q$ a $\R$-valued $m-1$-homogeneous polynomial. Then the field $u(x)=q(x)^{-1}P(x)$ is continuous at $x=0$ (where $u(0)=0$) but not differentiable at $x=0$, while the one-sided Gâteaux directional derivatives exist and equal 
\begin{align*}
u'(0;\xi)=\lim_{\eps\to 0^+}\frac{u(x+\eps \xi)-u(x)}{\eps}=\lim_{\eps\to 0^+}\frac{P(\eps \xi)}{\eps q(\eps \xi)}=\frac{P(\xi)}{q(\xi)}. 
\end{align*}
\end{ex}

The definition \ref{def:DerG} requires the limit to exist pointwise for all $\xi\in S^{n-1}$. This is typically much too strong for the purposes we are going to consider. Furthermore, it will turn out not to be the relevant notion of convergence. 

\begin{Def}[Weak one sided Gâteaux derivatives]
\label{def:DerGweak}
$u\in L^p(\Om,E)$, $1< p< \infty$, is said to have weak one-sided directional derivatives at $x\in \Om$ in $L^p(S^{n-1})$ if 
\begin{align}
\frac{u(x+\eps \xi)-u(x)}{\eps}\rightharpoonup u'(x;\xi)
\end{align}
weakly in $L^p(S^{n-1},\R^m)$ as $\eps \to 0^+$, i.e. for all $\phi\in L^q(S^{n-1},\R^m)$ with $\displaystyle \frac{1}{p}+\frac{1}{q}=1$
\begin{align*}
\lim_{\eps\to 0^+}\int_{S^{n-1}}\bigg\langle \frac{u(x+\eps \xi)-u(x)}{\eps},\phi(x)\bigg\rangle d\sigma(\xi)=\int_{S^{n-1}}\langle u'(x;\xi),\phi(x)\rangle d\sigma(\xi).
\end{align*}
We call $u'(x;\xi)$ \emph{the weak $L^p$ one sided directional derivative} of $u$ at $x$ and say that $u$ is \emph{$L^p$ weakly directionally differentiable} at $x$. 
\end{Def}

Many times it may be very difficult to prove that a function $u$ is weakly $L^p$-directionally differentiable for a.e. $x\in \Om$. However, by weak compactness it may be much easier to prove that that there exists a sequence $\{\eps_j\}_j$ of positive $\eps_j$ such that $\lim_{j\to \infty}\eps_j=0$ and 
\begin{align}
\wlim_{j\to \infty }\frac{u(x+\eps_j \xi)-u(x)}{\eps_j}=v(\xi)
\end{align}
for some $v\in L^p(S^{n-1},\R^m)$. In this case we say that $v$ is a \emph{$L^p$-weak directional derivative} of $u$ at $x$. Here the situation is very much analogous to the notion of \emph{tangent measure} of a Radon measure which we now recall. Consider the affine change of variables 
\begin{align*}
T_{a,r}(x)=\frac{x-a}{r},
\end{align*}
for $a>0$, which amount to ``zooming in'' around the point $a$. Given a vector valued Radon measure $\mu\in \mathcal{M}(\Om,\R^m)$ we consider the \emph{push-forward measure}
\begin{align*}
(T_{a,r})_\ast \mu(A):=\mu(a+rA)
\end{align*}
for all Borel sets $A\subset \R^n$. We say that a Randon measure $\nu\in  \mathcal{M}(\Om,\R^m)$ is a \emph{tangent measure} of $\mu$ (see \cite[Def. 3.3]{DeLe}), if there exists a positive sequence of radii $\{r_j\}_j$ converging to $0$ and positive sequence of real numbers $\{c_j\}_j$ such that 
\begin{align*}
\wslim_{j\to \infty}c_j(T_{a,r})_\ast \mu=\nu,
\end{align*}
where the convergence is the weak star convergence, i.e. for all $\phi\in C_0(\Om,\R^m)$,
\begin{align*}
\lim_{j\to \infty}c_j\int_{\Om}\langle \phi(x),d(T_{a,r})_\ast \mu(x)\rangle =\int_{\Om}\langle \phi(x),d\nu(x)\rangle. 
\end{align*}
Here typically the sequence $\{c_j\}_j$ is of the form $c_j=r_j^\alpha$ for some $\alpha>0$. 
The set of all tangent measures of $\mu$ at a point $x\in \Om$ is denoted by $\text{Tan}(\mu,x)$. By analogy we define:
\begin{Def}[Set of weak directional derivatives]
Let $u\in L^p(\Om,E)$, $1< p< \infty$. $\text{Der}_p(u,x)$ is defined be the set of all $L^p$-weak directional derivatives of $u$ at $x$, the set of all weak limits in $L^p(S^{n-1},\sigma;E)$ of sequences $\eps_j^{-1}(u(x+\eps_j\xi)-u(x))$ for some sequence $\{\eps_j\}_j$.
In particular if $u$ is $L^p$-weakly directionally differentiable at $x$, then $\text{Der}_p(u,x)=\{u'(x,\xi)\}$. 
\end{Def}

%==================NEW SUBSECTION=================================================================

\subsection{\sffamily Weak directional differentiability of Sobolev functions}

Given the notion of weak $L^p$ directional differentiability, it is natural to ask if functions $u\in W^{1,p}(\Om,E)$ are weakly $L^p$ directional differentiability for a.e. $x\in \Om$ or if $\text{Der}_p(u,x)\neq \varnothing$ for a.e. $x\in \Om$? 

In view of the ACL-characterization of Sobolev functions, see \cite{Hei1}, it is certainly true that for almost every line in a \emph{fixed direction} $u$ is absolutely continuous and thus have partial derivative in that direction for $\mathcal{H}^1$-a.e. every point on a.e. every such line. However, for the weak $L^p$ directional derivative we want to look at directional derivatives at every direction for a fixed point $x$, and it is not clear if Sobolev functions are  \emph{weakly $L^p$ directionally differentiable} at almost every point. In particular we recall that Sobolev functions are not differentiable a.e. when $p\leq n$, see \cite{EG}. 

In analysing this situation it is convenient to use a norm for Sobolev functions which does not directly employ weak derivatives. One such norm is provided by the Brezis-Bourgain-Mironesu formula, the BBM-formula for short, introduced in \cite{BBM}, which provides an equivalent norm for Sobolev functions. For that end we let $\{\rho_j\}_j$ be a sequence of mollifier satisfying 
\begin{itemize}
\item[(i)] $\rho_j\geq 0$ for all $j$.
\item[(ii)] $\rho_j(x)=\tilde{\rho}_j(\vert x\vert)$ for some $\tilde{\rho}_j:[0,+\infty)\to [0,+\infty)$ and all $j$.
\item[(iii)] $\int_{\R^n}\rho_j(x)dx=1$.
\item[(iv)] $\displaystyle \lim_{j\to \infty}\int_{\delta}^\infty\tilde{\rho}_j(r)r^{n-1}dr=0$ for every $\delta>0$. 
\item[(v)] $\tilde{\rho}_j=0$ for $r>1$. 
\end{itemize}

Define for $1\leq p<+\infty$ and $u\in L^1_{loc}(\R^n,\mathbb{E})$

\begin{align*}
D_{j,p}(u)(x)=\int_{\Om}\frac{\vert u(x)-u(y)\vert^p}{\vert x-y\vert^p}\rho_j(\vert x-y\vert)dy\quad \text{for a.e. $x\in \Om$}. 
\end{align*}
 
 It is shown in \cite{BBM,BN} that there exists constants depending only on $n$ and $p$ such that

\begin{align}
\int_{\R^n}D_{j,p}(u)(x)dx&\leq C_{n,p}\int_{\Om}\Vert Du(x)\Vert^pdx\label{eq:BBM1}\\
\lim_{j\to \infty}\int_{\Om} D_{j,p}(u)(x)dx&=\gamma_{n,p}\Vert Du(x)\Vert^p\quad \text{for a.e. $x\in \R^n$}\label{eq:BBM2}.
\end{align}
and in addition that 
\begin{align*}
\lim_{j\to \infty}\int_{\Om}\frac{\vert u(x+h)-u(x)-Du(x)h\vert^p}{\vert h\vert^p}\rho_j(\vert h\vert)dh=0\quad \text{for a.e. $x\in \R^n$}.
\end{align*}

provided that $u\in W^{1,p}(\R^n,E)$ and $p>1$.

\begin{Lem}\label{lem:WeakDiv}
Let $u\in W^{1,p}(\R^n,E)$ and $p>1$. Then for a.e. $x\in \R^n$ there exists a sequence $\{r_j(x)\}_j$ such that $\lim_{j\to \infty}r_j(x)=0$ and such that 
\begin{align*}
\sup_{j}\Vert r_j(x)^{-1}(u(x-r_j(x)\cdot)-u(x))\Vert_{L^2(S^{n-1}, E)}<+\infty. 
\end{align*}
\end{Lem}

\begin{proof}

Choose $\rho_\eps(\vert x\vert)=\eps^{-n}\chi_{B_\eps(0)}$. After spherical change of coordinates

\begin{align*}
D_{\eps,p}(u)(x)&=\int_{\R^n}\frac{\vert u(x+h)-u(h)\vert^p}{\vert h\vert^p}\rho_\eps(\vert h\vert)dh\\
&=\int_{\R^n}\bigg\vert\frac{ u(x+r\xi)-u(x)}{r}\bigg\vert^p\rho_\eps(r)r^{n-1}drd\sigma(\xi)\\
&=\int_{0}^\eps\frac{1}{\eps}\int_{S^{n-1}}\bigg\vert\frac{ u(x+r\xi)-u(x)}{r}\bigg\vert^p\bigg(\frac{r}{\eps}\bigg)^{n-1}drd\sigma(\xi).
\end{align*}

Define

\begin{align*}
g(x,r)&:=\int_{S^{n-1}}\bigg\vert\frac{u(x+r\xi)-u(x)}{r}\bigg\vert^pd\sigma(\xi)=\Vert r^{-1}(u(x+r\cdot)-u(x))\Vert^p_{L^p(S^{n-1},E)}\\
\end{align*}
 so that 
 
\begin{align*}
D_{\eps,p}(u)(x)&=\int_{0}^\eps \eps^{-n}g(x,r)r^{n-1}dr. 
\end{align*}
By \eqref{eq:BBM1}, $D_{\eps,p}(u)(x)\leq C_{n,p}\Vert Du(x)\Vert^p$ for a.e. $x$. On the other hand 

\begin{align*}
\int_{0}^{\eps} \eps^{-n}g(x,r)r^{n-1}dr\geq \frac{1}{2}\int_{\eps/2}^{\eps} \eps^{-1}g(x,r)dr
\end{align*}
and so $\eps \mapsto 2^{-k}\eps$ shows that for all $k\geq 0$ we have 

\begin{align*}
\int_{2^{-(k+1)}\eps}^{2^{-k}\eps} g(x,r)dr\leq 2^{-(k+1)}\eps C_{n,p}\Vert Du(x)\Vert^p.
\end{align*}

Thus, 

\begin{align*}
\int_{0}^\eps g(x,r)dr\leq \sum_{k=0}^\infty 2^{-(k+1)}\eps C_{n,p}\Vert Du(x)\Vert^p=\eps C_{n,p}\Vert Du(x)\Vert^p. 
\end{align*}

Let $I_N^\eps=\{x\in [0,\eps]: g(x,r)>N\}$. Since 

\begin{align*}
\int_{0}^\eps g(x,r)dr\geq \int_{0}^\eps \chi_{I_N^\eps}(x)g(x,r)dr\geq N\vert I_N^\eps\vert,
\end{align*}
and thus 
\begin{align*}
\vert I_N^\eps\vert\leq \frac{\eps}{N} C_{n,p}\Vert Du(x)\Vert^p. 
\end{align*}
By choosing $N=N(x)$ sufficently large we can ensure that $\eps^{-1}\vert I_N^\eps\vert<\delta<1$ for some $\delta\in (0,1)$. Thus we can always find a sequence $\{r_j\}_j$ tending to $0$ such that 
\begin{align*}
\sup_jg(x,r_j)<+\infty\quad \Longleftrightarrow \quad \sup_j\Vert r_j^{-1}(u(x+r_j\cdot)-u(x))\Vert^p_{L^p(S^{n-1},E)}<+\infty.
\end{align*}
\end{proof}

\begin{Prop}\label{prop:Der}
Let $u\in W^{1,p}_0(\Om,E)$ for $1<p<\infty$. Then $\text{Der}(u,x)\neq \varnothing$ for a.e. $x\in \Om$.
\end{Prop}

\begin{proof}
Since $u\in W^{1,p}_0(\Om,E)$, $u$ admits an extension to $W^{1,p}(\R^n,E)$. By Lemma \ref{lem:WeakDiv} the sequence $\{ r_j^{-1}(u(x+r_j\cdot)-u(x)\}_j$ is uniformly bounded in $L^p(S^{n-1},E)$ for a.e. $x\in \Om$. Hence it contains a weakly convergent subsequence which we do not relabel that converge weakly to some $v\in L^p(S^{n-1},E)$. 
\end{proof}

Proposition \ref{prop:Der} still leaves open the question whether or not a function $u\in W^{1,p}(\R^n,E)$ is still weakly $L^p$-directionally differentiable a.e. in $\R^n$. If we fix a $\nu\in S^{n-1}$ (or rather since lines in direction $\nu$ and $-\nu$ are the same we should really think of fixing a point in the real projective space $\mathbf{P}(\R^n)$) and let $L_\nu$ be the family of all lines in direction $\nu$, then by \cite[Thm. 7.4 ]{Hei1} $u\vert_{\ell}$ is $ACL_p$ for a.e. line. Hence $u'(x,\nu)$ exists pointwise for a.e. $x\in \R^n$. In particular, if 
\begin{align*}
U_\nu:=\{x\in \R^n: u'(x;\nu)\text{ exists}\},
\end{align*}
then for any bounded open set $\Om \subset \R^n$, $\vert \Om \setminus U_\nu\vert=0$. Moreover for any countable collection $\{\nu_j\}_j\subset S^{n-1}$, we have 
\begin{align*}
\vert \Om \setminus \bigcap_{j=1}^\infty U_{\nu_j}\vert=0
\end{align*}
since the countable collection of null sets is a null set. Thus for any $u\in W^{1,p}(\R^n,E)$, $u'(x,\nu)$ exists for a.e. $x\in \R^n$ and a dense set of directions in $S^{n-1}$. On the other hand the existence of the limit is in a stronger topology than needed. We leave it as an open problem to determine if in general a Sobolev function in $W^{1,p}(\R^n,E)$ is weakly $L^p$-directionally differentiable a.e. in $\R^n$.

%============NEW SUBSECTION=================================================================

\subsection{\sffamily Elliptic operators and the harmonic differential}

In this section we will only consider elliptic operators $\A$ as given in Definition \ref{def:Elliptic}. In addition, as has already been discussed in Section \ref{subsec:StokesAll} this implies that $\text{dom}_{p,\Om}(\A)\subset W^{1,p}_{loc}(\Om,E)$, and so by the observation in Section \ref{subsec:StokesAll} Stokes theorem holds for \emph{every} Green domain $U\Subset \Om$. Thus the essential limit can be replaced by the ordinary limit in \eqref{eq:FLimit2}. This will also motivate our notion of harmonic differential.  In particular we will see that the limit \eqref{eq:FLimit2} is independent of the symbol $\Aa$ when $\A$ is elliptic. To understand why this is the case, we first give an informal argument. If $\A$ is elliptic then that the symbol $\Aa(x,\xi)$ is injective for all $x\in \Om$ and all $\xi\in \R^n\setminus\{0\}$. This implies that $\Aa(x,\xi)$ has a left-inverse which can be taken to be the Moore-Penrose pseudoinverse $\Aa(x,\xi)^+$, so that $\Aa(x,\xi)^+\Aa(x,\xi)=I$ and $\Aa(x,\xi)^+$ is $C^1$ in $x$ and smooth for all $\xi \neq 0$. Assume that the pointwise limit
\begin{align*}
\lim_{\eps \to 0^+}\Aa(x,\xi)\frac{u(x+\eps \xi)-u(x)}{\eps}=v(x,\xi)
\end{align*}
exists. Then 
\begin{align*}
&\bigg\vert \frac{u(x+\eps \xi)-u(x)}{\eps}-\Aa^+(x,\xi)v(x,\xi)\bigg\vert=\bigg\vert \Aa^+(x,\xi)\Aa(x,\xi)\frac{u(x+\eps \xi)-u(x)}{\eps}-\Aa^+(x,\xi)v(x,\xi)\bigg\vert\\&\leq 
\Vert \Aa(x,\xi)^+\Vert\bigg\vert \Aa(x,\xi)\frac{u(x+\eps \xi)-u(x)}{\eps}-v(x,\xi)\bigg\vert
\end{align*}
and so limit 
\begin{align*}
\lim_{\eps \to 0^+}\frac{u(x+\eps \xi)-u(x)}{\eps}=\Aa^+(x,\xi)v(x,\xi)
\end{align*}
exists. The converse direction is clear. We now however see the fundamental difference between elliptic and non-elliptic operators. Assume for simplicity that $\A$ is a non-elliptic constant coefficient operator. Then there exists $\xi\neq 0$ such that $\text{ker}\,(\Aa(\xi))\neq \{0\}$. If we have that 
\begin{align*}
\frac{u(x+\eps \xi)-u(x)}{\eps}\in \text{ker}\,(\Aa(\xi))
\end{align*}
for all $\eps\in (0,\delta)$, then we have 
\begin{align*}
\lim_{\eps \to 0^+}\Aa(\xi)\frac{u(x+\eps \xi)-u(x)}{\eps}=0
\end{align*}
while $\lim_{\eps \to 0^+}\frac{u(x+\eps \xi)-u(x)}{\eps}$ need not exist.

We now recall the following results on spherical means due to E. Stein in \cite{St} and J. Bourgain in \cite{Bour1}.  

\begin{Def}
Let $u\in L^p(\R^n,E)$. Define the \emph{spherical maximal function} according to 
\begin{align}
\mathcal{M}_Su(x)=\sup_{t>0}\frac{1}{\sigma_{n-1}}\bigg\vert \int_{S^{n-1}}u(x+t\xi)d\sigma(\xi)\bigg\vert
\end{align}
\end{Def}

\begin{Thm}[Spherical means]
\label{thm:SM}
Let $u\in L^p(\R^n,E)$, for $\frac{n}{n-1}<p\leq \infty$ and $n\geq 3$ or $n=2$ and $p>2$. Then there exists a constant $C=C(p,n)>0$ such that
\begin{align}\label{eq:SM}
\Vert \mathcal{M}_Su\Vert_{L^p}\leq C \Vert u\Vert_{L^p}.
\end{align}
Consequently,
\begin{align*}
\lim_{\eps \to 0^+}\frac{1}{\sigma_{n-1}}\int_{S^{n-1}}u(x+\eps \xi)d\sigma(\xi)=u(x)
\end{align*}
for a.e. $x$ and in $L^p$. Furthermore, the inequality \eqref{eq:SM} fails for $p\leq \frac{n}{n-1}$ for $n\geq 3$ or $p\leq 2$ and $n=2$. 
\end{Thm}

For a proof of Theorem \ref{thm:SM} see for example \cite[Thm. 6.5.1, p. 476]{Gar}.

\begin{Lem}\label{lem:SLimit}
Let $1<p<\infty$, $u\in W^{1,p}(\Om,E)$ and assume that $u$ is weakly $L^p$-directionally differentiable at $x\in \Om$. Then 
\begin{align}\label{eq:SLimit}
\lim_{\eps\to 0^+}\frac{1}{\sigma_{n-1}}\int_{S^{n-1}}u(x+\eps\xi)d\sigma(\xi)=u(x).
\end{align}
\end{Lem}

\begin{proof}
Let $\{e_j\}_j$ be an ON-basis for $E$. Then $e_j\in L^q(S^{n-1},\sigma;E)$ with $\frac{1}{p}+\frac{1}{q}=1$. We have
\begin{align*}
&\frac{1}{\sigma_{n-1}}\int_{S^{n-1}} u(x+\eps\xi)d\sigma(\xi)=\frac{\eps}{\sigma_{n-1}}\int_{S^{n-1}} \frac{u(x+\eps\xi)-u(x)}{\eps}d\sigma(\xi)+u(x)\\
&=\sum_{j}\frac{\eps}{\sigma_{n-1}}\int_{S^{n-1}} \bigg\langle \frac{u(x+\eps\xi)-u(x)}{\eps},e_j\bigg\rangle d\sigma(\xi)e_j+u(x)\to u(x)
\end{align*}
as $\eps\to 0^+$ by assumption and the conclusion follows.
\end{proof}

By the Sobolev embedding theorem $W^{1,p}(\Om,E)\subset L^{p^\star}(\Om,E)$, with $\frac{1}{p^{\star}}=\frac{1}{p}-\frac{1}{n}$, and so $p^\star>p$, and if $p^\star>n/(n-1)$ then obviously \eqref{eq:SLimit} is true for a.e. $x\in \Om$ by the results of Stein and Bourgain, but the point of the lemma is that the point-wise limit exists independently of the assumption $p^\star>n/(n-1)$ if we instead assume that $u$ is weakly $L^p$-directionally differentiable at $x$.

\begin{Lem}\label{lem:SphereInt}
\begin{align*}
\frac{1}{\omega_n}\int_{S^{n-1}}\xi_j\xi_kd\sigma(\xi)=\delta_{jk}
\end{align*}
\end{Lem}

For a proof see \cite{Folland2}. 

\begin{Prop}\label{prop:Limit1}
Let $\Om \subset \R^n$ be a bounded domain and let $1< p<+\infty$ and assume that $\A$ satisfies the standard assumptions. Assume that $u\in W^{1,p}(\Om,E)$ is weakly $L^p$-directionally differentiable at $x$ so that 
\begin{align}
\frac{u(x+\eps\xi)-u(x)}{\eps} \rightharpoonup u'(x;\xi) \text{ in $L^p(S^{n-1},E)$}.
\end{align}
Then
\begin{align*}
&\lim_{\eps\to 0^+}\frac{1}{\vert B_\eps(x)\vert}\int_{\dv B_\eps(x)}\Aa(y,\nu(y))u(y)d\sigma(y)=\frac{1}{\omega_n}\int_{S^{n-1}}\Aa(x,\xi)u'(x;\xi)d\sigma(\xi)+\text{div}\,\Aa(x)u(x).
\end{align*}
\end{Prop}

\begin{proof}
We use the fact that $\Aa(x,\xi)$ is \emph{linear} in $\xi$ which implies the cancellation property 
\begin{align*}
\int_{S^{n-1}}\Aa(x,\xi)d\sigma(\xi)=0
\end{align*}
and so in particular for any $x$ in the Lebesgue set of $u$
\begin{align*}
\int_{S^{n-1}}\Aa(x,\xi)u(x)d\sigma(\xi)=0.
\end{align*}
Using this cancellation we find 
\begin{align*}
&\frac{1}{\vert B_\eps(x)\vert}\int_{\dv B_\eps(x)}\Aa(y,\nu(y))u(y)d\sigma(y)=\frac{1}{\eps^n \omega_n}\int_{S^{n-1}}\Aa(x+\eps \xi,\nu(x+\eps \xi))u(x+\eps \xi)\eps^{n-1}d\sigma(\xi)\\
&=\frac{1}{\eps \omega_n}\int_{S^{n-1}}\Aa(x+\eps \xi,\xi)u(x+\eps \xi)d\sigma(\xi)-\frac{1}{\eps \omega_n}\int_{S^{n-1}}\Aa(x,\xi)u(x)d\sigma(\xi)\\
&=\frac{1}{\omega_n}\int_{S^{n-1}}\frac{1}{\eps}(\Aa(x+\eps \xi,\xi)-\Aa(x,\xi))u(x+\eps \xi)d\sigma(\xi)+\frac{1}{ \omega_n}\int_{S^{n-1}}\Aa(x,\xi)\frac{u(x+\eps \xi)-u(x)}{\eps}d\sigma(\xi).
\end{align*}
Furthermore,
\begin{align*}
&\frac{1}{\omega_n}\int_{S^{n-1}}\frac{1}{\eps}(\Aa(x+\eps \xi,\xi)-\Aa(x,\xi))u(x+\eps \xi)d\sigma(\xi)\\&= \frac{1}{\omega_n}\int_{S^{n-1}}(\Aa(x+\eps \xi,\xi)-\Aa(x,\xi))\frac{u(x+\eps \xi)-u(x)}{\eps}d\sigma(\xi)
+ \frac{1}{\omega_n}\int_{S^{n-1}}\frac{1}{\eps}(\Aa(x+\eps \xi,\xi)-\Aa(x,\xi))u(x)d\sigma(\xi)\\
&=:I_1+I_2.
\end{align*}
First assume that $u\in C^1$. The mean value theorem implies that for some $\theta_\eps\in [0,\eps]$
\begin{align*}
\vert I_1\vert&\leq  \frac{1}{\omega_n}\int_{S^{n-1}}\Vert D_x\Aa(x+\theta_\eps \xi,\xi)\Vert\vert u(x+\eps \xi)-u(x)\vert d\sigma(\xi)\\
&\leq C\int_{S^{n-1}}\vert u(x+\eps \xi)-u(x)\vert d\sigma(\xi)=C\int_{S^{n-1}}\bigg\vert \int_0^\eps Du(x+t\xi)\xi dt\bigg\vert d\sigma(\xi)\\
&\leq C\int_{B_\eps(x)}\frac{\vert Du(y)\vert }{\vert x-y\vert^{n-1}}dy=\mathcal{I}_1(\vert Du\vert\chi_{B_\eps(x)}),
\end{align*}
where $\mathcal{I}_1$ denotes the $1$-Riesz potential. By approximation the same inequality holds for a general $u\in W^{1,p}(\Om,E)$. By \cite[Lem. 5.11]{Kinnunen}
\begin{align*}
\Vert \mathcal{I}_1(\vert Du\vert\chi_{B_\eps(x)})\Vert_{L^1(B_\eps(x))}\leq c(n)\vert B_\eps(x)\vert^{1/n}\Vert \vert Du\vert\Vert_{L^1(B_\eps(x))}\leq c'(n)\omega_n^{1/n}\eps^{1+n-n/p} \Vert \vert Du\vert\Vert_{L^p(B_\eps(x))}
\end{align*}
Thus 
\begin{align*}
\vert I_1\vert&\leq  d(n)\eps^{1+n-n/p} \Vert \vert Du\vert\Vert_{L^p(\Om)}
\end{align*}
for some constant $d=d(n)$. We now compute $I_2$. By our standard assumption $\Aa\in C^1(\Om, \LL(\R^n,\LL(E,F)))$ and Lemma \ref{lem:SphereInt}
\begin{align*}
I_2&= \frac{1}{\omega_n}\int_{S^{n-1}}\frac{1}{\eps}(\Aa(x+\eps \xi,\xi)-\Aa(x,\xi))u(x)d\sigma(\xi)=\frac{1}{\omega_n}\int_{S^{n-1}}(D_x\Aa(x,\xi)(\xi))u(x)d\sigma(\xi)+o(\eps)\\
&=\sum_{j,k=1}^n\frac{1}{\omega_n}(D_x\Aa(x,e_j)(e_k))\int_{S^{n-1}}\xi_j\xi_kd\sigma(\xi)u(x)+o(\eps)\\
&=\sum_{j}^n\frac{1}{\omega_n}(D_x\Aa(x,e_j)(e_j))\delta_{jk}u(x)+o(\eps)=\text{div}\,\Aa(x)u(x)+o(\eps). 
\end{align*}
Finally, using the assumption that is weakly $L^p$-directionally differentiable at $x$ we find that 
\begin{align*}
\lim_{\eps\to 0^+}\frac{1}{ \omega_n}\int_{S^{n-1}}\Aa(x,\xi)\frac{u(x+\eps \xi)-u(x)}{\eps}d\sigma(\xi)=\frac{1}{\omega_n}\int_{S^{n-1}}\Aa(x,\xi)u'(x;\xi)d\sigma(\xi),
\end{align*}
which completes the proof of the proposition.
\end{proof}

Proposition \ref{prop:Limit1} makes it obvious why the limit \eqref{eq:FLimit} can exist even though $u$ need not be differentiable at $x$. However, even the condition of weak $L^p$-directionally differentiability is on closer inspection unnecessary and much too strong. This can be seen by computing the integral 
\begin{align}
\frac{1}{\omega_n}\int_{S^{n-1}}\mathbb{A}(x,\xi)u'(x;\xi)d\sigma(\xi).
\end{align}
 To that end we recall the spherical harmonics of degree $k$ which are $k$-homogenous harmonic polynomials in $\R^n$ restricted to $S^{n-1}$. Let 
\begin{align*}
\mathcal{H}^k(S^{n-1}):=\{\text{$h$ is a $k$-homogenous harmonic polynomial in $\R^n$}\}. 
\end{align*}
Furthermore we let 
\begin{align*}
\mathcal{H}(S^{n-1})=\bigoplus_{k=0}^\infty\mathcal{H}^k(S^{n-1}).
\end{align*}
The spherical harmonics form a complete ON-basis for $L^2(S^{n-1},d\sigma)$ under the inner product
\begin{align*}
(f,g)=\frac{1}{\sigma_{n-1}}\int f(x)g(x)d\sigma(x).
\end{align*}
We now define $E$-valued spherical harmonics to be 
\begin{align*}
\mathcal{H}(S^{n-1},E)=\bigoplus_{k=0}^\infty\mathcal{H}^k(S^{n-1})\otimes E
\end{align*}

\noindent where an element $h\in \mathcal{H}^k(S^{n-1})\otimes E$ is simply an $E$-valued vector field on $\R^n$ whose coordinate functions all are spherical harmonics in $\mathcal{H}^k(S^{n-1})$. This space again form a complete ON-basis for $L^2(S^{n-1},\sigma;E)$ under the inner product 
\begin{align*}
(u,v):=\frac{1}{\sigma_{n-1}}\int \langle u(x),v(x)\rangle_Ed\sigma(x),
\end{align*}
where $\langle\cdot,\cdot\rangle_E$ denotes the inner product in $E$. In particular, we note that since every linear polynomial is harmonic we have the identification
\begin{align*}
\mathcal{H}^1(S^{n-1})\otimes E\cong \LL(\R^n, E). 
\end{align*}
In the case when $E=\R$, an ON-basis is given by the \emph{spherical harmonics}. By tensoring with $E$, this also gives an ON-basis for $L^2(S^{n-1},\sigma;E)$. 
\begin{Def}[Zonal spherical harmonics]
The \emph{zonal spherical harmonics} $Z_d^n(x,y)$ are defined through
\begin{align*}
Z_d^n(x,y)=\frac{1}{\sigma_{n-1}}(C^{n/2}_d(\langle x,y\rangle)-C^{n/2}_{d-2}(\langle x,y\rangle))
\end{align*}
for $x,y\in S^{n-1}$, and where $C^\lambda_d(t)$ are the Gegenbauer polynomials given by the generating function expansion 
\begin{align*}
\frac{1}{(1-2tr+r^2)^{\lambda}}=\sum_{d=0}^\infty C^\lambda_d(t)r^d, 
\end{align*}
for $\lambda\in \R$ and $d\in \N$. 
In particular,
\begin{align*}
Z_1^n(x,y)=\frac{1}{\sigma_{n-1}}C_1^{n/2}(\langle x,y\rangle)=\frac{n}{\sigma_{n-1}}\langle x,y\rangle=\frac{1}{\omega_{n}}\langle x,y\rangle
\end{align*}
\end{Def}

If $\{Y_l^d\}_{l}$ are the spherical harmonics homogeneous of degree $d$ then 
\begin{align*}
Z_d^n(x,y)=\sum_{l=1}^{n(d)}Y_l^d(x)Y_l^d(y). 
\end{align*}

The zonal spherical harmonics define orthogonal projections $\mathbf{Z}_d^n: L^2(S^{n-1},\sigma)\to \mathcal{H}^d(S^{n-1})$ through 
\begin{align*}
\mathbf{Z}_d^n(u)(x)=\frac{1}{\sigma_{n-1}}\int_{S^{n-1}}Z_d^n(x,y)u(y)d\sigma(y). 
\end{align*}
These projections extend directly to orthogonal projection $\mathbf{Z}_d^n: L^2(S^{n-1},\sigma;E)\to \mathcal{H}^d(S^{n-1},E)$ with the same kernel. In particular we have for any $u,v\in L^2(S^{n-1},\sigma;E)$
\begin{align*}
( \mathbf{Z}_d^n(u)(x),\mathbf{Z}_l^n(v)(x))=0 
\end{align*}
whenever $d\neq l$.

Let $u'(x;\xi)\in L^2(S^{n-1},\sigma;E)$ be a the weak $L^p$ directional derivative of $u$ at $x$. After expansion in $\mathcal{H}^d(S^{n-1},E)$ using the zonal harmonics we have 
\begin{align*}
u'(x;\xi)=\sum_{d=0}^\infty\mathbf{Z}_d^n(u(x;\cdot))(\xi)=:\sum_{d=0}^\infty u_d'(x;\xi)
\end{align*}
where in particular, 
\begin{align*}
u_0(x;\xi)&=\frac{1}{\sigma_{n-1}}\int_{S^{n-1}}u'(x;\xi)d\sigma(\xi),\\
u_1(x;\xi)&=\frac{1}{\sigma_{n-1}}\int_{S^{n-1}}Z_1(\xi,\eta)u'(x;\eta)d\sigma(\eta)=\frac{1}{\sigma_{n-1}}\int_{S^{n-1}}\langle \xi,\eta \rangle u'(x;\eta)d\sigma(\eta).
\end{align*}

\begin{Lem}\label{lem:LimitInt}
Let $v\in L^2(S^{n-1},E)$ and let 
\begin{align*}
v_d(\xi):=\frac{1}{\sigma_{n-1}}\int_{S^{n-1}}Z_d^n(\xi,\eta)v(\eta)d\sigma(\xi).
\end{align*}
Then 
\begin{align*}
\int_{S^{n-1}}\Aa(\xi)v(\xi)d\sigma(\xi)=\int_{S^{n-1}}\Aa(\xi)v_1(\xi)d\sigma(\xi).
\end{align*}
\end{Lem}

\begin{proof}
Since $\Aa(x,\xi)$ is linear in $\xi$, $\Aa(x,\cdot)\in \LL(\R^n,\LL(E,F))=\mathcal{H}^1(S^{n-1},\LL(E,F))$ and so 
\begin{align*}
\Aa(x,\xi)=\frac{1}{\sigma_{n-1}}\int_{S^{n-1}}Z^n_1(\xi,\eta)\Aa(x,\eta)d\sigma(\eta).
\end{align*}
Thus, by Fubini's theorem
\begin{align*}
\int_{S^{n-1}}\Aa(x,\xi)v(\xi)d\sigma(\xi)&=\int_{S^{n-1}}\int_{S^{n-1}}Z^n_1(\xi,\eta)\Aa(x,\eta)d\sigma(\eta)v(\xi)d\sigma(\xi)\\
&=\int_{S^{n-1}}\Aa(x,\eta)\frac{1}{\sigma_{n-1}}\int_{S^{n-1}}Z^n_1(\eta,\xi)v(\xi)d\sigma(\xi)d\sigma(\eta)\\
&=\int_{S^{n-1}}\Aa(x,\eta)v_1(\eta)d\sigma(\eta).
\end{align*}
\end{proof}

Using Lemma \ref{lem:LimitInt} we see that 
\begin{align*}
\mathscr{A}u(x)=\frac{1}{\omega_n}\int_{S^{n-1}}\mathbb{A}(x,\xi)u'(x;\xi)d\sigma(\xi)=\frac{1}{\omega_n}\int_{S^{n-1}}\mathbb{A}(x,\xi) u_1'(x;\xi)d\sigma(\xi),
\end{align*}

Thus $\A u$ {\bf only depends on $u'(x,\xi)$ through $u_1'(x;\xi)$ i.e. the linear part of $u'(x;\xi)$} in $\xi$. All other frequency information of $u'(x;\xi)$ is lost.

\begin{Def}[Harmonic differential]
\label{def:HarmonicDiff}
$u\in L^p(\Om,E)$ is said to have a \emph{harmonic differential} at $x$ if there exists a linear map $T_x\in \LL(\R^n, E)$ such that 
\begin{align*}
&\esslim_{\eps \to 0^+}\frac{1}{\omega_{n}}\int_{S^{n-1}}\langle \xi,\eta\rangle \frac{u(x+\eps \eta)-u(x)}{\eps}d\sigma(\eta)\\
&=\esslim_{\eps \to 0^+}\frac{1}{\omega_{n}}\int_{S^{n-1}}\frac{u(x+\eps \eta)-u(x)}{\eps}\otimes \eta d\sigma(\eta)(\xi)=T_x(\xi)
\end{align*}
We write $T_x=D_{\Delta}u(x)$, and call $D_{\Delta}u(x)$ the harmonic differential of $u$ at $x$. 
\end{Def} 

We can think of the harmonic differential in the following way. Consider the difference quotient $u(y)-u(x)$ on the dilated and translated sphere $x+\eps S^{n-1}$. After a linear change of variables consider the blow up $\eps^{-1}(u(x+\eps \xi)-u(x))$ on the unit sphere which we view as boundary values of an harmonic function $U_\eps$ in the unit ball $B_1(0)$. A Taylor expansion of $U_\eps$ at $0$ yields 
\begin{align*}
U_\eps(x)=U_\eps (0)+DU_\eps(0)x+ \sum_{\vert \alpha \vert=2}\frac{1}{\alpha !}\dv^\alpha U_\eps(0)x^\alpha.
\end{align*}
Then 
\begin{align*}
\lim_{\eps \to 0^+}DU_\eps(0)=D_{\Delta}u(x). 
\end{align*}

\begin{Thm}\label{thm:LimitSobolev}
If $u\in W^{1,p}(\Om, E)$, $1\leq p<+\infty$, then 
\begin{align*}
D_wu(x)=D_{\Delta}u(x), 
\end{align*}
and 
\begin{align}\label{eq:StrongDiff}
D_wu(x)=\lim_{\eps\to 0^+}\frac{1}{\omega_n}\int_{S^{n-1}}\frac{u(x+\eps \xi)-u(x)}{\eps}\otimes \xi d\sigma(\xi)
\end{align}
for a.e. $x\in \Om$ and in $L^p(\Om,E)$. Furthermore, for any elliptic operator $\A$
\begin{align*}
\A u(x)=\sum_{j=1}^n\mathbb{A}(x,e_j)D_wu(x)(e_j)+B(x)u(x)
\end{align*}
for a.e. $x\in \Om$.
\end{Thm}

\begin{proof}
By Theorem \ref{thm:EquivalenceFuglede}, $D_w=D_f$. By Theorem \ref{thm:LimitFuglede}, $D_f=D_{\Delta}$. Finally, 
\begin{align*}
&\frac{1}{\omega_n}\int_{S^{n-1}}\mathbb{A}(x,\xi)\frac{u(x+\eps \xi)-u(x)}{\eps}d\sigma(\xi)=\sum_{j=1}^n\mathbb{A}(x,e_j)\frac{1}{\omega_n}\int_{S^{n-1}}\xi_j\frac{u(x+\eps \xi)-u(x)}{\eps}d\sigma(\xi)\\
&=\sum_{j=1}^n\mathbb{A}(x,e_j)\frac{1}{\omega_n}\int_{S^{n-1}}\frac{u(x+\eps \xi)-u(x)}{\eps}\otimes \xi d\sigma(\xi)(e_j)
\end{align*}
and so 
\begin{align*}
\lim_{\eps \to 0^+}\frac{1}{\omega_n}\int_{S^{n-1}}\mathbb{A}(x,\xi)\frac{u(x+\eps \xi)-u(x)}{\eps}d\sigma(\xi)&=\sum_{j=1}^n\mathbb{A}(x,e_j)D_wu(x)(e_j)
\end{align*}
for a.e. $x\in \Om$. 
\end{proof}

Theorem \ref{thm:LimitSobolev} has the deficiency that we need to let $x\in \mathscr{L}_u^p\cap  \mathscr{L}_{\A u}^p$ while by Theorem \ref{thm:LimitFuglede} we know that the limit exists for all $x\in \mathscr{L}_{\A u}^p$, i.e. $x$ need not belong to the $p$-Lebesgue set of $u$. This leads us naturally to the next section when we study the limit \eqref{eq:FLimit2} also for non-elliptic operators.

%============NEW SUBSECTION=================================================================

\subsection{\sffamily Non-elliptic homogeneous constant coefficient operators}\label{sec:NonEllip}

In this section we study the limit \eqref{eq:FLimit2} also for non-elliptic operators homogeneous constant coefficient operators, i.e. operators for which the principal symbol $\Aa(\xi)$ is not injective for all $\xi\in S^{n-1}$. The reason for restricting to this class of constant coefficient operators is that certain cancellation phenomena becomes much more subtle in the variable coefficient case and a careful analysis of this case has to be deferred to future work.  Furthermore, we want to make the ideas as transparent as possible, which is easiest in the case of homogenous constant coefficient operators. For non-elliptic operators it is no longer true that $\text{dom}_{p,\Om}(\A)\subset W^{1,p}_{loc}(\Om,E)$, and in particular the trace of an element $u\in \text{dom}_{p,\Om}(\A)$ need not exists on $\dv U$ as an $L^p$-function for every Green domain $U\Subset \Om$. It is therefore necessary to work with the essential limit instead. Moreover, in the previous section we only considered points $x\in \Om$ which belonged to the Lebesgue set of  $u\in \mathscr{L}_u^p$. This assumption will be relaxed.  However before we do that we will first start with an informal discussion. Assume for simplicity that $u$ is continuous on $\R^n$ and consider the pointwise limit 
\begin{align}\label{eq:limNonE}
\lim_{\eps \to 0^+}\Aa(\xi)\frac{u(x+\eps \xi)-u(x)}{\eps}.
\end{align}
for a fixed $\xi\in S^{n-1}$. If the kernel of $\Aa(\xi)$ is non-trivial, then it is clear that the limit \eqref{eq:limNonE} can exists without the limit $\lim_{\eps \to 0^+}(u(x+\eps \xi)-u(x))/\eps$ existing. In particular, should $u(x+\eps \xi)-u(x)\in \text{ker}(\Aa(\xi))$ for all $0<\eps <\delta$ then the limit is 0 independent of whether the limit $\lim_{\eps\to 0^+}(u(x+\eps \xi)-u(x))/\eps$ exists or not. What matters here is of course whether the projected part of $(u(x+\eps \xi)-u(x))/\eps$ onto the orthogonal complement of the kernel of $\Aa(\xi)$ has a limit or not. A convenient way to express this is via the Moore-Penrose pseudoinverse of $\Aa(\xi)$, a concept which we now recall.

 Let $A\in \LL(E,F)$. The Moore-Penrose pseudoinverse of $A$ is the unique linear map $A^+\in \LL(F,E)$ that satisfies the four Moore-Penrose equations 
\begin{enumerate}
\item $AA^+A=A$
\item $A^+AA^+=A^+$
\item $(AA^+)^\ast=AA^+$
\item $(A^+A)^\ast=A^+A$
\end{enumerate}

In general any linear map $A^g\in \LL(F,E)$ that satisfies the first equation $AA^gA=A$ is called a pseudoinverse of $A$, which is in general not unique. Moreover, 
\begin{align*}
AA^+=\text{Proj}_{\text{im}(A)},\quad A^+A=\text{Proj}_{\text{im}(A^\ast)}=\text{Proj}_{(\text{ker}(A))^\perp}.
\end{align*}

Furthermore, in \eqref{eq:limNonE} we are implicitly assuming that $u$ is continuous at $x$ or at least that $x$ belongs to the Lebesgue set of $u$. This assumption is however unnecessary as we will now see.

\begin{Def}[Weak spherical blowup limit]
Let $u\in L^p(\Om,E)$, $1< p <+\infty$. We say that $u$ has a weak $L^p$ spherical blow up limit at $x$ if there exists a function $v\in L^p(S^{n-1},\sigma;E)$
\begin{align*}
u(x+\eps \xi)\rightharpoonup v_x(\xi) \quad \text{as $\eps \to 0^+$}
\end{align*}
with the limit being the essential weak limit (outside a null set of $\eps$). If $u$ admits a spherical blow up limit $v_x$ at a.e $x\in \Om$, we write $v_x(\xi)=u_S(x,\xi)$.
\end{Def}

\begin{ex}[Differential form with tangential discontinuity at a hyperplane]
\label{ex:DisHyper}
Let $\omega,\eta: \R^{n}\to \Lambda \R^n$ be $C^1$-differential forms such that $\omega(x)\neq \eta(x)$ for all $x$ such that $x_{n}=0$. Define the differential form
\begin{align*}
u(x)=\omega(x)\chi_{x_{n+1}>0}(x)+\eta(x)\chi_{x_{n}<0}(x).
\end{align*}
Then $u: \R^{n}\to \Lambda \R^{n}$ is differential form with a discontinuity across the hyperplane $x_{n}=0$ and the hyperplane does not belong to the Lebesgue set of $u$. Moreover, for any open set $U\subset \R^{n+1}$ which intersects the hyperplane $x_{n}=0$,
$u\notin W^{1,1}(U,\Lambda \R^{n})$ (though $u\in \text{SBV}(U,\Lambda \R^{n})$, the space of differential forms whose coefficients belong to the space of functions of special bounded variation, see \cite{AFP}). 
One easily sees that the spherical blow up limit for any point $x$ on the hyperplane $x_{n}=0$ equals $u_S(x,\xi)=\omega(x)\chi_{\langle \xi,e_{n+1}\rangle>0}+\eta(x)\chi_{\langle \xi,e_{n+1}\rangle<0}$.
\end{ex}

\begin{Def}[Adapted weak spherical blowup limit]
Let $u\in L^p(\Om,E)$, $1< p < +\infty$. We say that $u$ has a weak $L^p$ spherical blow up limit adapted to the PDE operator $\A$ at $x$ if there exists a function $v\in L^p(S^{n-1},\sigma;E)$
\begin{align*}
\Aa(\xi)^+\Aa(\xi)u(x+\eps \xi)\rightharpoonup v(x,\xi) \quad \text{as $\eps \to 0^+$}
\end{align*}
with the limit being the essential weak limit (outside a null set of $\eps$). If $u$ admits a spherical blow up limit $v_x$ adapted to $\A$ at a.e $x\in \Om$, we write $v_x(\xi)=u_{\A,S}(x,\xi)$.
\end{Def}

\begin{ex}\label{ex:BlowUp}
Consider the vector field $u(x)=x/\vert x\vert^n$. Then $u\in L^{p}(B_1(0),\R^n)$ for $1\leq p< n/(n-1)$ and $du=0$ in the sense of distributions. Furthermore, $0$ does not belong to the Lebesgue set of $u$ and $u$ does not have a spherical blow up limit in $L^p$ for any $1\leq p< n/(n-1)$. On the other hand $u$ does have a weak $L^p$-blow up limit adapted to the exterior derivative $d$ at $0$. Indeed,  
\begin{align*}
\Aa_d(\xi)w=\xi \wedge w
\end{align*}
and 
\begin{align*}
\Aa_d(\xi)^+w=\frac{\xi}{\vert \xi \vert^2} \ri w, 
\end{align*} 
 where $\xi\ri $ denotes interior multiplication by $\xi$ for any multivector $w\in \Lambda \R^n$. Thus 
 \begin{align*}
 \Aa(\xi)^+\Aa(\xi)u(0+\eps \xi)=\frac{\xi}{\vert \xi \vert^2}\ri\bigg(\xi \wedge \frac{\eps\xi}{\eps^n\vert \xi\vert^n}\bigg)=0. 
 \end{align*}
 Thus $u$ has a $L^p$-blow up limit $0$ at $0$ adapted to $d$ at $x$. 
 \end{ex}

\begin{Lem}\label{lem:Cancel}
Let $1<  p <+\infty$ and let $\A$ be a  homogeneous constant coefficient operators.  Assume that $u\in \text{dom}_{p,\Om}(\A)$, and that $u$ has a weak $L^p$ spherical blow up limit $u_{\A,S}(x,\cdot)$ at $x$ and that $x$ is a Lebesgue point of $\A u$. 
Then $u_{\A,S}(x,\cdot)$ satisfies the cancellation condition 
\begin{align}\label{eq:Cancel}
\int_{S^{n-1}}\Aa(\xi)u_S(x,\xi)d\sigma(\xi)=0. 
\end{align}
\end{Lem}

\begin{proof}
Since $\A$ is a pure flux operator we have for a.e. $\eps>0$
\begin{align*}
\frac{1}{\eps \omega_n}\int_{S^{n-1}}\Aa(\xi)u(x+\eps \xi)d\sigma(\xi)=\frac{1}{\vert B_\eps(x)\vert}\int_{B_\eps(x)}\A u(y)dy.
\end{align*}
By the assumption of the existence of $u_{\A,S}(x,\xi)$ and the first Moore-Penrose equation
\begin{align*}
&\frac{1}{ \omega_n}\int_{S^{n-1}}\Aa(\xi)u_S(x,\xi)d\sigma(\xi)=\esslim_{\eps \to 0^+}\frac{1}{ \omega_n}\int_{S^{n-1}}\Aa(\xi)\Aa(\xi)^+\Aa(\xi)u(x+\eps \xi)d\sigma(\xi)\\
&=\esslim_{\eps \to 0^+}\frac{1}{ \omega_n}\int_{S^{n-1}}\Aa(\xi)u(x+\eps \xi)d\sigma(\xi)\\
&=\lim_{\eps \to 0^+}\frac{\eps}{\vert B_\eps(x)\vert}\int_{B_\eps(x)}\A u(y)dy=0.
\end{align*}
\end{proof}

\begin{Lem}[Fundamental cancellation lemma]
\label{lem:CancelFund}
Let $u\in \text{dom}_{p,\Om}(\A_w)$. Assume the assumptions of Lemma \ref{lem:Cancel}.
\begin{align*}
&\lim_{\eps \to 0^+}\int_{S^{n-1}}\Aa(\xi)u(x+\eps \xi)d\sigma(\xi)\\
&=\lim_{\eps \to 0^+}\int_{S^{n-1}}\Aa(\eta)\frac{1}{\sigma_{n-1}}\int_{S^{n-1}}\Aa(\xi)^+\Aa(\xi)\frac{u(x+\eps \xi)-u_{\A,S}(x,\xi)}{\eps}\otimes \xi d\sigma(\xi)(\eta)d\sigma(\eta)
\end{align*}
\end{Lem}

\begin{proof}
By Lemma \ref{lem:Cancel} and the identity $\Aa(\xi)=\Aa(\xi)\Aa(\xi)^+\Aa(\xi)$
\begin{align*}
&\int_{S^{n-1}}\Aa(\xi)u(x+\eps \xi)d\sigma(\xi)=\int_{S^{n-1}}\Aa(\xi)(u(x+\eps \xi)-u_{\A,S}(x,\xi))d\sigma(\xi)\\
&=\int_{S^{n-1}}\Aa(\xi)\Aa(\xi)^+\Aa(\xi)(u(x+\eps \xi)-u_{\A,S}(x,\xi))d\sigma(\xi)\\
&=\int_{S^{n-1}}\frac{1}{\sigma_{n-1}}\int_{S^{n-1}}\langle \xi,\eta\rangle \Aa(\eta)d\sigma(\eta)\Aa(\xi)^+\Aa(\xi)(u(x+\eps \xi)-u_{\A,S}(x,\xi))d\sigma(\xi)\\
&=\int_{S^{n-1}}\int_{S^{n-1}}\Aa(\eta)\frac{1}{\sigma_{n-1}}\langle \xi,\eta\rangle \Aa(\xi)^+\Aa(\xi)(u(x+\eps \xi)-u_{\A,S}(x,\xi))d\sigma(\xi)d\sigma(\eta)\\
&=\int_{S^{n-1}}\Aa(\eta)\frac{1}{\sigma_{n-1}}\int_{S^{n-1}}\Aa(\xi)^+\Aa(\xi)(u(x+\eps \xi)-u_{\A,S}(x,\xi))\otimes \xi d\sigma(\xi)(\eta)d\sigma(\eta).
\end{align*}
\end{proof}

Lemma \ref{lem:CancelFund} suggests that one should consider the quantity 
\begin{align}
D_{\eps,\A}u(x):=\frac{1}{\sigma_{n-1}}\int_{S^{n-1}}\Aa(\xi)^+\Aa(\xi)\frac{u(x+\eps \xi)-u_{\A,S}(x,\xi)}{\eps}\otimes \xi d\sigma(\xi)
\end{align}
where $D_{\eps,\A}u(x)\in \LL(\R^n, E)$. If $\esslim_{\eps \to 0+}D_{\eps,\A}u(x)$ exists we denote it by $D_\A u$. Then by Lemma \ref{lem:CancelFund}

\begin{align}
\esslim_{\eps\to 0^+}\frac{1}{\eps \omega_n}\int_{S^{n-1}}\Aa(\xi)u(x+\eps \xi)d\sigma(\xi)=\frac{1}{\omega_n}\int_{S^{n-1}}\Aa(\eta)D_{\A}u(x)(\eta)d\sigma(\eta)
\end{align}
If $u\in \text{dom}_{p,\Om}(\A)$ for a non-elliptic operator, then in general $Du\in \mathcal{D}'(\Om, \LL(\R^n,E))$. We split $Du=D_au+D_su$, such that $D_au\in L^1_{loc}(\Om,\LL(\R^n,E))$ and $D_su\in \mathcal{D}'(\Om, \LL(\R^n,E))$ is not equal to an $L^1_{loc}$ function anywhere on its support. A natural question is if $D_\A u=D_au$? We leave this question for the future. 

\begin{ex}[Example with singularity at a point]
Let $u(x)=x/\vert x\vert^{n}$ and let $\A=d$, the exterior derivative on differential forms. Then $u\in L^p(B_1(0),\R^n)$ for $1\leq p<n/(n-1)$. Furthermore, $x=0$ does not belong to the Lebesgue set of $u$. The classical total derivative equals  
\begin{align*}
Du(x)=\sum_{j=1}^n\dv_ju(x) \otimes e_j=\sum_{j=1}^n\bigg(\frac{e_j}{\vert x\vert^n}-n\frac{x_jx}{\vert x\vert^{n+2}}\bigg)\otimes e_j=\frac{1}{\vert x\vert^n}I-n\frac{x\otimes x}{\vert x\vert^{n+2}}
\end{align*}
for $x\neq 0$. In addition 
\begin{align*}
\vert Du(x)\vert^2=\frac{n(n-1)}{\vert x\vert^{n-2}},
\end{align*}
where $\vert Du(x)\vert=\text{tr}(Du(x)^\ast Du(x))$ is the Hilbert-Schmidt norm.  We now compute its weak derivative. Let $\Phi\in C^\infty_0(B_1(0),\LL(\R^n))$. By definition
\begin{align*}
&( D_wu(x),\Phi(x))=-\int_{B_1(0)}\langle u(x),\text{div}\,\Phi(x)\rangle dx=\lim_{\eps\to 0^+}-\int_{B_1(0)\setminus B_\eps(0)}\langle u(x),\text{div}\,\Phi(x)\rangle dx\\
&=\lim_{\eps\to 0^+}\int_{B_1(0)\setminus B_\eps(0)}\langle Du(x),\Phi(x)\rangle dx-\int_{\dv B_\eps(0)}\langle u(x)\otimes \nu(x),\Phi(x)\rangle d\sigma(x)\\
&=\int_{B_1(0)}\langle Du(x),\Phi(x)\rangle dx-\lim_{\eps\to 0^+}\int_{S^{n-1}}\langle \xi\otimes \xi,\Phi(\eps \xi)\rangle d\sigma(\xi)\\
&=\int_{B_1(0)}\langle Du(x),\Phi(x)\rangle dx-\int_{S^{n-1}}\langle \xi\otimes \xi,\Phi(0)\rangle d\sigma(\xi)\\
&=\int_{B_1(0)}\langle Du(x),\Phi(x)\rangle dx-\sum_{j,k=1}^n\langle e_j\otimes e_k,\Phi(0)\rangle\int_{S^{n-1}} \xi_j\xi_kd\sigma(\xi)\\
&=\int_{B_1(0)}\langle Du(x),\Phi(x)\rangle dx-\sum_{j,k=1}^n\langle e_j\otimes e_k,\Phi(0)\rangle\omega_n\delta_{jk}\\
&=\int_{B_1(0)}\langle Du(x),\Phi(x)\rangle dx-\langle \omega_n I,\Phi(0)\rangle.
\end{align*}
Thus
\begin{align*}
D_w u(x)=Du(x)-\omega_n I\delta_0. 
\end{align*}
Moreover, the limit
\begin{align*}
\lim_{\eps\to 0^+}\frac{1}{\omega_n}\int_{S^{n-1}}\frac{u(0+\eps \xi)}{\eps}\otimes \xi d\sigma(\xi)=\lim_{\eps\to 0^+}\frac{1}{\eps^{n}\omega_n}\int_{S^{n-1}}\xi \otimes \xi d\sigma(\xi)
\end{align*}
does not exist. On the other hand, similar computations show that 
\begin{align*}
d_wu(x)=0, \quad \text{div}_w u(x)=\sigma_{n-1}\delta_0. 
\end{align*}
We also see that $u$ does not have a weak $L^p$-blow up limit on $S^{n-1}$ at $x=0$. However, by example \ref{ex:BlowUp}, $u$ has the $L^{p}$-weak blow up limit adapted to $d$ $u_{d,S}(0,\xi)=0$ on $S^{n-1}$ for any $1< p<\infty$ at $x=0$. Furthermore, $x=0$ belongs to the Lebesgue set of $d_wu$. Finally,
\begin{align*}
\lim_{\eps\to 0^+}\frac{1}{\omega_n}\int_{S^{n-1}}\xi \wedge \frac{u(0+\eps \xi)}{\eps}d\sigma(\xi)=\lim_{\eps\to 0^+}\frac{1}{\eps^{n}\omega_n}\int_{S^{n-1}}\xi \wedge \xi d\sigma(\xi)=0. 
\end{align*}

\end{ex}

\begin{ex}[Example with discontinuity across a hyperplane]
This is a variation of example \ref{ex:DisHyper}. Let $\omega,\eta: \R^{n}\to \Lambda \R^{n-1}$ be $C^1$-differential forms such that $\omega(x)\neq \eta(x)$ for all $x=(x_1,x_2,...,x_n)$ such that $x_{n}=0$ and satisfying $d\omega(x)=d\eta(x)=0$ for all $x$. Define the differential form
\begin{align*}
u(x)=e_n\wedge \omega(x)\chi_{x_{n}>0}(x)+e_n\wedge \eta(x)\chi_{x_{n}<0}(x).
\end{align*}
Then $u: \R^{n}\to \Lambda \R^{n}$ is differential form with a discontinuity across the hyperplane $x_{n}=0$ and the hyperplane does not belong to the Lebesgue set of $u$. One easily checks that the spherical blow up limit for any point $x$ on the hyperplane $x_{n}=0$ equals $u_S(x,\xi)=e_n\wedge \omega(x)\chi_{\langle \xi,e_{n}\rangle>0}(x)+e_n\wedge \eta(x)\chi_{\langle \xi,e_{n}\rangle<0}(x)$. We now compute $d_fu(x)$ when $x$ belongs to the hyperplane from the limit formula. Let $S^{n-1}_+:=\{x\in S^{n-1}: \langle x,e_n\rangle>0\}$, $S^{n-1}_-:=\{x\in S^{n-1}: \langle x,e_n\rangle<0\}$ and similarly for $B_\eps(0)^\pm$. 
\begin{align*}
&\frac{1}{\omega_n}\int_{S^{n-1}}\xi \wedge (u(x+\eps \xi)-u_S(x))d\sigma(\xi)\\&=\frac{1}{\omega_n}\int_{S^{n-1}_+}\xi \wedge (e_n\wedge \omega(x+\eps \xi)-e_n\wedge \omega(x))d\sigma(\xi)+\frac{1}{\omega_n}\int_{S^{n-1}_-}\xi \wedge (e_n\wedge \eta(x+\eps \xi)-e_n\wedge \eta(x))d\sigma(\xi)\\
&=\frac{1}{\vert B_\eps^+(x)\vert}\int_{ B_\eps^+(x)}d\omega(y)dy+\frac{1}{\vert B_\eps^-(x)\vert}\int_{ B_\eps^+(x)}d\eta(y)dy-\frac{1}{\omega_n}\int_{B_1(0)\cap \{x_n=0\}}e_n \wedge (e_n\wedge \omega(x+\eps \xi)-e_n\wedge \omega(x))d\sigma(\xi)\\&+\frac{1}{\omega_n}\int_{B_1(0)\cap \{x_n=0\}}e_n \wedge (e_n\wedge \omega(x+\eps \xi)-e_n\wedge \omega(x))d\sigma(\xi)=0.
\end{align*}
Thus, $d_fu(x)=0$.

We now compute the weak exterior derivative of $u$. Let $\phi\in C^\infty_0(Q,\Lambda \R^n)$, where $Q=(-1,1)\times (-1,1)$. Then 
\begin{align*}
( d_wu, \phi)=\int_{Q}\langle u(x), -d^\ast\phi(x)\rangle dx=\lim_{\eps\to 0^+}\int_{Q\setminus R_\eps}\langle u(x), -d^\ast\phi(x)\rangle dx,
\end{align*}
where $R_\eps=[-1,1]\times [-\eps,\eps]$. Set $Q_\eps^+=\{x\in Q: x_n>\eps\}$ and  $Q_\eps^-=\{x\in Q: x_n<-\eps\}$. Integration by parts give 
\begin{align*}
( d_wu, \phi)&=\lim_{\eps\to 0^+}\int_{Q\setminus R_\eps}\langle du(x), \phi(x)\rangle dx+\int_{\{x_n=\eps\}}\langle \nu(x)\wedge u(x), \phi(x)\rangle d\sigma(x)+\int_{\{x_n=-\eps\}}\langle \nu(x)\wedge u(x), \phi(x)\rangle d\sigma(x)\\
&=\lim_{\eps\to 0^+}\int_{Q_\eps^+}\langle -e_n\wedge d\omega(x), \phi(x)\rangle dx+\int_{Q_\eps^-}\langle -e_n\wedge d\eta(x), \phi(x)\rangle dx\\
&+\int_{\{x_n=\eps\}}\langle -e_n\wedge (e_n\wedge \omega(x)), \phi(x)\rangle d\sigma(x)+\int_{\{x_n=-\eps\}}\langle e_n\wedge (e_n\wedge \eta(x)), \phi(x)\rangle d\sigma(x)=0.
\end{align*}
This shows that 
\begin{align*}
d_wu(x)=0,
\end{align*}
and so $d_wu(x)=d_fu(x)$ as expected. 
\end{ex}

%============NEW SUBSECTION=================================================================

\subsection{\sffamily Constant rank operators}

\begin{Def}[Constant rank operators]
\label{def:ConstantR}
A first order operator 
\begin{align*}
\A=\sum_{j=1}^nA_j(x)\dv_j+A_0(x)
\end{align*}
on a domain $\Om\subset \R^n$ such that $A_j\in C^1(\Om,\LL(E,F)$ for $j=1,...,n$ and $A_0\in C(\Om,\LL(E,F))$ is called const rank operator if the principal symbol satisfy 
\begin{align*}
\text{rank}(\Aa(x,\xi))=k \text{ for all $x\in \Om$ and $\xi\in \R^n\setminus \{0\}$}
\end{align*}
for some non-negative integer $k$. 
\end{Def}

\begin{ex}
Let $\A=\text{div}$, the divergence operator on matrix fields in $C^\infty(\R^n, \LL(\R^n,E))$. Then for any $M\in \LL(\R^n,E)$, 
\begin{align*}
\Aa_{\text{div}(\xi)}(M)=M(\xi). 
\end{align*}
$M(\xi)=0$ if and only if $\xi$ is orthogonal to every row in $M$. Complete $\xi$ to an ON-basis for $\R^n$, so that $\R^n=\text{span}\{\xi, e_2,....,e_n\}$. Any $M$ can be written as
\begin{align*}
M=u_1\otimes \xi+\sum_{j=2}^nu_j\otimes e_j
\end{align*}
for some $u_j\in E$. Thus, for any $\xi \in S^{n-1}$
\begin{align*}
\Aa_{\text{div}(\xi)}(M)=\vert \xi\vert^2 u_1. 
\end{align*}
Thus $\text{rank}(\Aa_{\text{div}(\xi)})=1$ for all $\xi\in S^{n-1}$, and so $\text{div}$ is a constant rank operator. 
\end{ex}

\begin{ex}
Let $\A=d$, the exterior derivative on differential forms on $C^\infty(\R^n, \Lambda \R^n)$. Then for any $u\in \Lambda \R^n$, 
\begin{align*}
\Aa_{d}(\xi)u=\xi \wedge u. 
\end{align*}
For any fixed $\xi \in S^{n-1}$, we have $\R^n=V\oplus \text{span}(\xi)$, where $V=(\text{span}(\xi))^\perp$. Any $u\in \Lambda \R^n$ can be written as an orthogonal sum 
\begin{align*}
u=\xi \wedge u_1+u_2
\end{align*}
where $u_1,u_2\in \Lambda V$. Thus,
\begin{align*}
\Aa_{d}(\xi)u=\xi \wedge u_2
\end{align*}
$\text{rank}(\Aa_{\text{div}(\xi)})=2^{n-1}$ for all $\xi\in S^{n-1}$
\end{ex}

\begin{ex}\label{ex:DiracWave}
Let $\A=\mathcal{D}_M$, the hyperbolic Dirac operator acting on $C^\infty(\R^n,C\ell_M (\R^{n-1,1}))$, where $C\ell_M (\R^{n-1,1})$ is the Clifford algebra generated from the Minkowski inner product on $\R^n$. If we let $\{e_0,e_1,...,e_{n-1}\}$ be an ON-basis, then with Clifford multiplication denoted by $\gp$, we have 
\begin{align*}
e_0\gp e_0=-1, \,\,\, e_1\gp e_1=e_2\gp e_2=...=e_n\gp e_n=1. 
\end{align*}
In addition, for any $u\in C\ell_M (\R^{n-1,1})$
\begin{align*}
\Aa_{\mathcal{D}_M}(\xi)u=\xi \gp u. 
\end{align*}
Since $\xi \gp \xi=0$ whenever $\xi\in \Gamma_L$, the Lorentz cone, it follows that $\Aa_{\mathcal{D}_M}(\xi)$ is not invertible when $\xi \in S^{n-1}\cap \Gamma_L$, and so $\mathcal{D}_M$ is not elliptic (as indicated by the name hyperbolic Dirac operator). On the other hand  $\Aa_{\mathcal{D}_M}(\xi)$ is invertible for all $\xi S^{n-1}\setminus \Gamma_L$, so $\mathcal{D}_M$ is not a constant rank operator. 
\end{ex}

\begin{Thm}[Schulenberger-Wilcox coercivity theorem]
\label{thm:Coer}
Let $\A$ be a first order constant rank operator as in Definition \ref{def:ConstantR} and the coefficients of $\A$ satisfies $A_j\in C^1_b(\R^n,\LL(E,F))$ and $B\in C_b(\R^n,\LL(E,F))$, and the subscript $b$ denotes bounded. Then there exists a $\mu>0$ such that
\begin{align*}
\sum_{j=1}^n\Vert \dv_j u\Vert^2_{L^2(\Om,E)}\leq \mu(\Vert \A u\Vert^2_{L^2(\Om,F)}+\Vert u\Vert^2_{L^2(\Om,E)})
\end{align*}
for all $u\in \text{dom}_{2,\Om}(\A_0)\cap \text{ker}_{2,\Om}(\A_0)^\perp$. 
\end{Thm}

Theorem \ref{thm:Coer} has later been extended also to extensions in $L^p$ for $1<p<\infty$ and also to operators of higher order, see \cite{K}. In fact, in the case when $\A$ has constant coefficients the coercive 
inequality
\begin{align*}
\sum_{j=1}^n\Vert \dv_j u\Vert^p_{L^p(\Om,E)}\leq \mu(\Vert \A u\Vert^p_{L^p(\Om,F)}+\Vert u\Vert^p_{L^p(\Om,E)})
\end{align*}
for all $1<p<\infty$ characterizes constant rank operators as is shown in \cite{GR20}. Let $\A$ be a constant rank operator and satisfy the conditions of Theorem \ref{thm:Coer}. Let $u\in \text{dom}_{2,\Om}(\A)$. Since we are only interested in the existence of the limit \eqref{eq:FLimit2} in the interior of $\Om$, we may always multiply $u$ be a function $\phi\in C^\infty_0(\Om)$ such that $\phi=1$ in a neighbourhood of the point of interest. Moreover, since $\phi u\in \text{dom}_{2,\R^n}(\A_w)=\text{dom}_{2,\R^n}(\A_0)$, we let $\mathbf{P}_\A: \text{dom}_{2,\R^n}(\A_w)\to \text{dom}_{2,\R^n}(\A_w)$ be the orthogonal projection operator onto $\text{ker}_{2,\R^n}(\A_0)^\perp$. Then, we have an orthogonal decomposition of any $u\in \text{dom}_{2,\R^n}(\A_0)$ according to 
\begin{align*}
u=\mathbf{P}_\A u+(I-\mathbf{P}_\A)u:=u_W+u_0
\end{align*}
By Theorem \ref{thm:Coer}, $u_W\in W^{1,2}(\R^n,E)$ and $\A u_0=0$. In particular, any irregular behaviour not possible for Sobolev functions comes entirely from $u_0$. Thus,
\begin{align*}
\int_{\dv U}\Aa(x,\nu(x))u(x)d\sigma(x)&=\int_{\dv U}\Aa(x,\nu(x))u_W(x)d\sigma(x)+\int_{\dv \Om}\Aa(x,\nu(x))u_0(x)d\sigma(x)\\&=\int_{\dv U}\Aa(x,\nu(x))u_W(x)d\sigma(x)
\end{align*}
for $p$-a.e. Green domain $U\Subset \Om$ by the generalized Cauchy integral theorem for $\A$. The study of the limit \eqref{eq:FLimit2} therefore reduces to the case of for Sobolev functions.

\begin{Lem}\label{lem:Volume}
Let $1\leq p\leq \infty$ and let $u\in W^{1,p}(\Om,E)$. Then for $\eps>0$ sufficiently small so that $B_\eps(x)\Subset \Om$
\begin{align*}
\mathscr{A}_\eps u(x)&=\frac{1}{\omega_n\eps}\int_{B_\eps(x)}\mathbb{A}\bigg(\frac{y-x}{\vert y-x\vert}\bigg)Du(y)\frac{y-x}{\vert y-x\vert^n} dy
\end{align*}
\end{Lem}

\begin{proof}
First assume that $u\in C^1(\Om,E)$. Using that 
\begin{align*}
u(x+r\xi)-u(x)=\int_0^rDu(x+s\xi)\xi ds
\end{align*}
we have 
\begin{align*}
\mathscr{A}_ru(x)&=\frac{1}{\omega_n}\int_{S^{n-1}}\mathbb{A}(\xi)\frac{(u(x+r\xi)-u(x))}{r}d\sigma(\xi)=\frac{1}{\omega_n}\int_{S^{n-1}}\mathbb{A}(\xi)\frac{1}{r}\int_0^rDu(x+s\xi)\xi dsd\sigma(\xi)\\
&=\frac{1}{\omega_nr}\int_{S^{n-1}}\int_0^r\mathbb{A}(\xi)Du(x+s\xi)\xi dsd\sigma(\xi)=\frac{1}{\omega_nr}\int_{S^{n-1}}\int_0^r\mathbb{A}(\xi)Du(x+s\xi)\frac{\xi}{s^{n-1}} s^{n-1}dsd\sigma(\xi)\\
&=\frac{1}{\omega_nr}\int_{B_r(0)}\mathbb{A}\bigg(\frac{y}{\vert y\vert}\bigg)Du(x+y)\frac{y}{\vert y\vert^{n}} dy\\
&=\frac{1}{\omega_nr}\int_{B_r(x)}\mathbb{A}\bigg(\frac{y-x}{\vert y-x\vert}\bigg)Du(y)\frac{y-x}{\vert y-x\vert^n} dy. 
\end{align*}
The general case now follows from an approximation argument. 
\end{proof}

Thus if $\A$ is a homogeneous constant rank operator, then by Lemma \ref{lem:Volume} we have for $1<p<\infty$ and $u\in \text{dom}_{p,\Om}(\A)$
\begin{align*}
\mathscr{A}_\eps u(x)&=\frac{1}{\omega_n\eps}\int_{B_\eps (x)}\mathbb{A}\bigg(\frac{y-x}{\vert y-x\vert}\bigg)Du_W(y)\frac{y-x}{\vert y-x\vert^n} dy. 
\end{align*}

If we let 
\begin{align*}
\mu_\A:=\max_{\xi\in S^{n-1}}\Vert \Aa(\xi)\Vert,
\end{align*}
then by \cite[Lem. 3.9]{Kinnunen}

\begin{align*}
\vert \A_\eps u(x)\vert\leq \frac{\mu_\A}{\omega_n\eps}\int_{B_\eps (x)}\frac{\Vert Du_W(y)\Vert}{\vert y-x\vert^{n-1}} dy\leq c\mathcal{M}(\Vert Du_W\Vert)(x) 
\end{align*}
for some constant $c=c(n)>0$. This implies that for any $u\in \text{dom}_{2,\Om}(\A)$, 
\begin{align*}
 \A^\star u(x)\leq c\mathcal{M}(\Vert Du_W\Vert)(x) 
\end{align*}
where $\A^\star$ is the maximal operator in Section \ref{sec:MaxOp}.

%============NEW SUBSUBSECTION=================================================================

\subsection{\sffamily Wave cone and non-elliptic operators}

\begin{Def}[Characteristic variety and characteristic hypersurfces]
\emph{The characteristic variety} $\mathcal{V}_{\mathscr{A}}$ of $\mathscr{A}$ is defined according to 
\begin{align*}
\mathcal{V}_{\mathscr{A}}=\{(x,\xi)\in \Om\times S^{n-1}: \mathbb{A}(x,\xi)\text{ is not injective}\}. 
\end{align*}
If $E=F$, then 
\begin{align*}
\mathcal{V}_{\mathscr{A}}=\{(x,\xi)\in \Om\times S^{n-1}: \det(\mathbb{A}(x,\xi))=0\}. 
\end{align*}
Furthermore a $C^1$-surface $\Sigma$ is called \emph{characteristic with respect to} $\A$ if for all $x\in \Sigma$, $\Aa(x,\nu(x))$ is not injective. 
\end{Def}

\begin{Def}[Wave cone]
\emph{The wave cone} $\Lambda_{\mathscr{A}}(x)$ at $x$ is defined according to 
\begin{align*}
\Lambda_{\mathscr{A}}(x)=\bigcup_{\xi \in S^{n-1}}\text{ker}\,\mathbb{A}(x,\xi).
\end{align*}
\end{Def}

The wave cone was introduced by Murat and Tatar in their work on compensated compactness. The concept has later proven decisive also in other areas of PDE, in particular in the work \cite{DePR}. Loosely speaking the complement of the wave cone can be thought of as the elliptic directions of an operator $\A$. Indeed, if $\A$ is a homogeneous constant coefficient operator and
\begin{align*}
u(x+\eps\xi)-u(x)\in \Lambda_{\mathscr{A}}^c\quad \text{for all $\xi\in S^{n-1}$ and all $\eps>0$, \text{ $\eps$ sufficiently small}},
\end{align*}
then 
\begin{align*}
u(x+\eps\xi)-u(x)\notin \text{ker}\, \Aa(\xi)
\end{align*}
for any $\xi \in S^{n-1}$ and $\eps>0$. Let $\mathcal{C}\Subset \Lambda_{\mathscr{A}}^c$ be a closed cone containing $0$ and assume that for $u\in C(U,E)$, $u(U)-u(U)\subset \mathcal{C}$. Then by assumption, $(u(y)-u(x))/\eps\in \mathcal{C}$ for all $x,y\in U$.

Using Lemma \ref{lem:CancelFund}, we have 
\begin{align*}
\int_{S^{n-1}}\Aa(\xi)u(x+\eps \xi)d\sigma(\xi)=\int_{S^{n-1}}\Aa(\eta)D_{\eps,\A}u(x)(\eta)d\sigma(\eta),
\end{align*}
where 
\begin{align*}
D_{\eps,\A}u(x)=\frac{1}{\sigma_{n-1}}\int_{S^{n-1}}\Aa(\xi)^+\Aa(\xi)\frac{u(x+\eps \xi)-u(x)}{\eps}\otimes \xi d\sigma(\xi).
\end{align*}

Since $\Aa(\xi)^+\Aa(\xi)$ is the orthogonal projection onto $(\text{ker}\,\Aa(\xi))^\perp$, the assumptions on $u$ implies that 
\begin{align*}
\Aa(\xi)^+\Aa(\xi)(u(x+\eps \xi)-u(x))=(u(x+\eps \xi)-u(x))
\end{align*}
whenever $\eps>0$ is sufficiently small. Thus 
\begin{align*}
&\lim_{\eps\to 0^+}\int_{S^{n-1}}\Aa(\xi)u(x+\eps \xi)d\sigma(\xi)\\&=\lim_{\eps\to 0^+}\int_{S^{n-1}}\Aa(\eta)\frac{1}{\sigma_{n-1}}\int_{S^{n-1}}\Aa(\xi)^+\Aa(\xi)\frac{u(x+\eps \xi)-u(x)}{\eps}\otimes \xi d\sigma(\xi)(\eta)d\sigma(\eta)\\
&=\lim_{\eps\to 0^+}\int_{S^{n-1}}\Aa(\eta)\frac{1}{\sigma_{n-1}}\int_{S^{n-1}}\frac{u(x+\eps \xi)-u(x)}{\eps}\otimes \xi d\sigma(\xi)(\eta)d\sigma(\eta)
\end{align*}
which is the same expression as for an elliptic operator.

\begin{ex}[Mizhohata operators]
Consider the Mizohata operators
\begin{align*}
\mathscr{A}_{M^k}=\dv_x+ix^k\dv_y,
\end{align*}
$k\in \mathbb{N}$. As an operator $\mathscr{A}_{M^k}: C^\infty(\R^2,\R^2)\to C^\infty(\R^2,\R^2)$
\begin{align*}
\mathbb{A}_{M^k}(\xi)=\xi_1 I+x^k\xi_2J
\end{align*}
acting on $\R^2$. Note that 
\begin{align*}
(\xi_1 I-x^k\xi_2J)(\xi_1 I+x^k\xi_2J)=\xi_1^2I+x^k\xi_1\xi_2J-x^k\xi_1\xi_2J+x^{2k}\xi_2I=(\xi_1^2+x^{2k}\xi_2^2)I,
\end{align*}

and so $\mathbb{A}_{M^k}(\xi)$ is invertible if and only if $\xi_1^2+x^{2k}\xi_2^2\neq 0$ in which case 

\begin{align*}
\mathbb{A}_{M^k}^{-1}(\xi)=\frac{1}{\xi_1^2+x^{2k}\xi_2^2}(\xi_1 I-x^k\xi_2J).
\end{align*}
In the case when $x=0$, then $\mathbb{A}_{M^k}(\xi)=\xi_1$, and $\text{ker}(\mathbb{A}_{M^k}(\xi))=0$ if $\xi_1\neq 0$, and $\text{ker}(\mathbb{A}_{M^k}(\xi))=\R^2$ if $\xi_1=0$.  
Hence 
\begin{align*}
\Lambda_{\mathscr{A}_{M^k}}(x,y)=\left\{
    \begin{array}{rl}
      \{0\}& x\neq 0,\\
      \R^2 & \text{if } x = 0.
      \end{array} \right.
\end{align*}
The importance of the Mizhota operators is that like the Lewy operator, the equation $\A u=v$ for $v\in C^\infty(U,F)$ is generically (in Baire category sense) {\bf not} solvable in any neighbourhood $U$ of the origin when $k$ is odd. Furthermore,
\begin{align*}
\mathcal{V}_{\A_{M^k}}=(\{0\}\times \R)\times ((\{0\}\times \R)\cap S^{2}).
\end{align*}

\end{ex}

\begin{ex}
Let $\A=\text{curl}$ be the curl operator on $n\times n$-matrix fields. Then $\text{ker}(\Aa_{\text{curl}}(\xi))=\{a\otimes \xi: a\in \R^n\}$ and 
\begin{align*}
\Lambda_{\text{curl}}=\{a\otimes b: a,b\in \R^n\}. 
\end{align*}
Thus, $\Lambda_{\text{curl}} $ consists of all rank one maps. In addition,
\begin{align*}
\mathcal{V}_{\text{curl}}=\R^n\times S^{n-1}
\end{align*}
and so every $C^1$-hypersurface of $\R^n$ is characteristic with respect to $\text{curl}$.
\end{ex}

\begin{ex}
Let $\A=\mathcal{D}_M$ be the hyperbolic Dirac operator in Example \ref{ex:DiracWave}. Then 
\begin{align*}
\Aa_{\mathcal{D}_M}(\xi)u=\xi \gp u. 
\end{align*}
Since $\xi \gp \xi=0$ whenever $\xi\in \Gamma_L=\{\xi \in \R^n:\xi_0^2=\xi_1^2+...\xi_{n-1}^2\}$ the Lorentz cone, $\text{ker}(\Aa_{\mathcal{D}_M}(\xi))=\{0\}$ whenever $\xi \notin \Gamma_L$. On the other hand when $\xi \in \Gamma_L$, then 
\begin{align*}
\text{ker}(\Aa_{\mathcal{D}_M}(\xi))=\span{\xi}\gp C\ell(\R^{n-1}), 
\end{align*}
where $C\ell(\R^{n-1})$ is the euclidean Clifford subalgebra of $C\ell(\R^{n-1,1})$. This implies that 
\begin{align*}
\Lambda_{\mathcal{D}_M}=\Gamma_L\gp C\ell(\R^{n-1})
\end{align*}
and 
\begin{align*}
\mathcal{V}_{\mathcal{D}_M}=\R^n\times (\Gamma_L\cap S^{n-1} ). 
\end{align*}
\end{ex}

Operators may be further distinguished by considering the behaviour of the kernel and image of $\Aa$. The following notions ware introduced by Van Schaftingen in \cite{VanS} in studying endpoint estimates for elliptic operators. 

\begin{Def}[Cancelling and co-cancelling operators]
A first order constant coefficient operator is called a \emph{cancelling operator} if
\begin{align}\label{eq:Cancel}
\bigcap_{\xi\in S^{n-1}}\text{im}(\Aa(\xi))=\{0\}
\end{align}
and \emph{cocancelling operator} if
\begin{align}\label{eq:CoCancel}
\bigcap_{\xi\in S^{n-1}}\text{ker}(\Aa(\xi))=\{0\}
\end{align}
\end{Def}
Note in particular that every elliptic operator is cocancelling. Condition \eqref{eq:Cancel} already distinguish between different kinds of \emph{elliptic operators} as shown in \cite{VanS}.  Furthermore, if a constant coefficient operator $\A$ is not cocancelling, let $V\subset E$ be the linear subspace $V= \bigcap_{\xi\in S^{n-1}}\text{ker}(\Aa(\xi))$. This implies directly that any $u\in L^1(\Om,V)$ satisfy the conditions of Theorem \ref{thm:FugledeChar}, which implies that $u\in \text{dom}_{p,\Om}(\A)$ and $\A u=0$. Thus for operators that are not co-cancelling, the regularity of an element of $u\in \text{dom}_{p,\Om}(\A)$ is in general no better than $L^p$. Note that by \cite[Prop. 6.2]{VanS}, the cancelling condition is equivalent to 
\begin{align*}
\text{span}(\bigcup_{\xi\in S^{n-1}}(\text{im}(\Aa(\xi)))^\perp)=F.
\end{align*}
Similarly, the cocancelling condition is equivalent to 
\begin{align*}
\text{span}(\bigcup_{\xi\in S^{n-1}}\text{im}(\Aa(\xi)^\ast))=\text{span}(\bigcup_{\xi\in S^{n-1}}(\text{ker}(\Aa(\xi)))^\perp)=E.
\end{align*}
\begin{ex}
Let $\A=D_S$ be symmetric gradient on vector fields. Then $\text{Im}(\Aa_{D_S}(\xi))=\text{span}\{\frac{1}{2}(a\otimes \xi+\xi\otimes a:a\in \R^n)\}$ and $\bigcap_{\xi\in S^{n-1}}\text{Im}(\Aa_{D_S}(\xi))=\{0\}$ which implies that $D_S$ is cancelling. 
\end{ex}

\begin{ex}\cite[Prop. 3.1]{VanS}
Let $\A=\text{div}$ on vector fields. Then $\text{ker}(\Aa_{\text{div}}(\xi))=\xi^\perp$ and $\bigcap_{\xi\in S^{n-1}}\xi^\perp=\{0\}$ which implies that $\text{div}$ is cocancelling. 
\end{ex}

Another important condition imposed on the wave cone is the spanning condition.

\begin{Def}[Spanning condition]
The wave cone $\Lambda_{\mathscr{A}}(x)$ at $x$ is said to satisfy the \emph{spaning condition}, (SC) for short, if 
\begin{align}
\text{span}(\Lambda_{\mathscr{A}}(x))=E. 
\end{align}
\end{Def}

While we will not be using the (SC) condition we point out that it occurs naturally in the theory of Young measures and weak sequential lower semicontinuity of integrals see \cite{KR}.

\begin{Def}
A first order linear partial differential system $\mathscr{A}$ is said to be \emph{hyperbolic in direction $\nu$} if the linear map
\begin{align*}
\mathbb{A}(\nu)^{-1}\mathbb{A}(\xi)
\end{align*}
is a symmetric positive operator for every $\xi\perp \nu$.
\end{Def}

\begin{Lem}
If $\mathscr{A}$ is hyperbolic in direction $\nu$, then there exists a $n=\alpha \nu+\beta \xi$, with $\xi\perp \nu$ such that $\mathbb{A}(n)$ is not invertible. 
\end{Lem}

\begin{proof}
\begin{align*}
\mathbb{A}(n)=\alpha \mathbb{A}(\nu)+\beta \mathbb{A}(\xi)
\end{align*}
Then $\mathbb{A}(n)$ is invertible if and only if $\mathbb{A}(\nu)^{-1}\mathbb{A}(n)$ is. But
\begin{align*}
\mathbb{A}(\nu)^{-1}\mathbb{A}(n)=\alpha I+\beta  \mathbb{A}(\nu)^{-1}\mathbb{A}(\xi)
\end{align*}
and 
\begin{align*}
\det[\alpha I+\beta  \mathbb{A}(\nu)^{-1}\mathbb{A}(\xi)]=\prod_{j=1}^n(\alpha +\beta\lambda_j)
\end{align*}
where $\lambda_1\leq \lambda_2\leq...\leq \lambda_n$ are the positive eigenvalues of $\mathbb{A}(\nu)^{-1}\mathbb{A}(\xi)$. But for every $j=1,2,..,n$ there are real $\alpha,\beta$ such that $(\alpha +\beta \lambda_j)=0$. 
\end{proof}

The set of direction $\xi \in \R^n\setminus\{0\}$ such that the principal symbol of a homogeneous hyperbolic operator is not invertible is called the Lorentz cone $\Gamma_\A$ associated to $\A$. In particular,
\begin{align*}
\Aa(\xi)^+\Aa(\xi)=\text{id} 
\end{align*} 
for $\sigma$-a.e $\xi \in S^{n-1}$. In this sense hyperbolic operators seem more similar to elliptic operator than constant rank non-elliptic operators. Yet we know that hyperbolic operators behave very differently from elliptic operators. The difference is of course that for elliptic operators,
\begin{align*}
\min_{\xi\in S^{n-1}}\Vert \Aa(\xi)\Vert>0,
\end{align*}
while for hyperbolic operators we have 
\begin{align*}
\min_{\xi\in S^{n-1}}\Vert \Aa(\xi)\Vert=0. 
\end{align*}
Assume for example $u$ has one sided directional derivatives $u'(x;\xi)$ at $x$, such that $u'(x,\cdot)\in C(S^{n-1}\setminus \Gamma_{\A})$. Assume furthermore that for some $\xi_0\in \Gamma_\A$
\begin{align*}
\lim_{\substack{\xi \to \xi_0\\ \xi \in S^{n-1}\setminus \Gamma_A }}\frac{u'(x;\xi)}{\vert u'(x;\xi)\vert} \in \text{ker}(\Aa(\xi_0))
\end{align*}
exists while
\begin{align*}
\lim_{\substack{\xi \to \xi_0\\ \xi \in S^{n-1}\setminus \Gamma_\A }}\vert u'(x;\xi)\vert =+\infty
\end{align*}
and $u'(x,\cdot)\notin L^1(S^{n-1},\sigma;E)$. Then we could still have that $\Aa(\cdot)u'(x,\cdot)\in L^1(S^{n-1},\sigma;E)$, where the vanishing of $\Vert \Aa(\xi)\Vert$ as $\xi \to \xi_0$ compensates for the blow up in $u'(x,\xi)$ at $\xi_0$. Before concluding this work we would like to point out that the limit for hyperbolic operators
\begin{align*}
\lim_{\eps\to 0^+}\frac{1}{\omega_n}\int_{S^{n-1}}\Aa(\xi)u(x+\eps \xi)d\sigma(\xi)
\end{align*}
is in some sense unphysical as, for any $\eps>0$, the value of the integral depends on future points as well as past points. From a physical stand point the value of $\A u(x)$ should only depend on $u(y)$ for points in the \emph{causal past} $J_-(x)$ of the point $x$. We hope to return to these issues with a more detailed analysis in future work. 

\bibliographystyle{alpha}

\end{document}